\documentclass[a4paper, 11pt, headings = small, abstract]{scrartcl}

\usepackage[english]{babel}
\usepackage[utf8]{inputenc}
\usepackage[T1]{fontenc}
\usepackage{geometry}

\usepackage[semibold]{libertine}

\usepackage{amsmath}
\usepackage{amsfonts}
\usepackage{amssymb}
\usepackage{amsthm}
\usepackage{mathabx}
\usepackage{mathtools}
\usepackage{paralist}
\usepackage{comment}
\usepackage{graphicx}
\usepackage{caption}
\usepackage{subcaption}
\usepackage{algorithm}
\usepackage{algpseudocode}
\usepackage{paralist}
\usepackage{booktabs}
\usepackage{multirow}

\usepackage[authoryear]{natbib}
\usepackage[dvipsnames]{xcolor}

\usepackage{tikz}
\usetikzlibrary{positioning}
\usetikzlibrary{calc,tikzmark}
\usetikzlibrary{arrows,automata}
\usetikzlibrary{backgrounds}
\usetikzlibrary{arrows.meta}

\definecolor{darkblue}{rgb}{0.0, 0.0, 0.55}
\usepackage[colorlinks=true,
            linkcolor=blue,
            citecolor=blue,
            pdfborder={0 0 0}]{hyperref}
\usepackage{orcidlink}
\usepackage{url}

\usepackage{dsfont}
\usepackage{bm}

\usepackage{graphicx}
\usepackage{enumerate}
\usepackage{stackengine}
\stackMath

\newcommand{\triangleA}{\stackon[0.5pt]{A}{\kern0.2em\scalebox{0.6}{$\bm{\triangle}$}}}

\newcommand{\triangleAa}{\stackon[0.5pt]{A}{\kern0.2em\scalebox{0.6}{$\bm{\triangle}$}}\,\!^{(1)}}

\newcommand{\triangleAb}{\stackon[0.5pt]{A}{\kern0.2em\scalebox{0.6}{$\bm{\triangle}$}}\,\!^{(2)}}

\usepackage{tikz}
\usetikzlibrary{positioning}
\usetikzlibrary{backgrounds}

\newtheorem{theorem}{Theorem}[section]
\newtheorem{proposition}[theorem]{Proposition}
\newtheorem{lemma}[theorem]{Lemma}
\newtheorem{corollary}[theorem]{Corollary}

\theoremstyle{definition}
\newtheorem{definition}[theorem]{Definition}

\newtheorem*{condition*}{Condition}

\makeatletter
\newenvironment{condition}[1]
  {%
    \def\condition@tag{#1}%
    \begin{condition*}[#1]%
    \phantomsection
    \def\@currentlabel{(#1)}%
  }
  {%
    \end{condition*}%
  }
\makeatother

\theoremstyle{remark}
\newtheorem{remark}[theorem]{Remark}
\newtheorem{example}[theorem]{Example}

\allowdisplaybreaks[4]

\numberwithin{equation}{section}

\newcommand{\R}{\mathbb{R}}
\newcommand{\Q}{\mathbb{Q}}

\newcommand{\N}{\mathbb{N}}
\newcommand{\eps}{\varepsilon}

\newcommand{\Bb}{\mathbb {B}}

\newcommand{\Gb}{\mathbb {G}}
\newcommand{\Hb}{\mathbb {H}}

\newcommand{\Nc}{\mathcal{N}}
\newcommand{\Pc}{\mathcal{P}}

\newcommand{\an}{{\mathrm{an}}}
\newcommand{\pa}{{\mathrm{pa}}}

\newcommand{\An}{{\mathrm{An}}}

\newcommand{\Des}{{\mathrm{De}}}

\newcommand{\weak}{\rightsquigarrow}

\newcommand{\cweak}[1]{\underset{#1}{\overset{\mathbb P}{\rightsquigarrow}}}

\newcommand{\Prob}{\mathbb{P}}   
\newcommand{\Exp}{\mathbb E}
\newcommand{\Var}{\operatorname{Var}}
\newcommand{\Cov}{\operatorname{Cov}}

\newcommand{\diag}{\operatorname{diag}}

\newcommand{\diff}{\mathrm{d}}

\renewcommand{\subset}{\subseteq}

\newcommand{\gghat}{\widehat H_n}
\newcommand{\ggmean}{H_n}
\newcommand{\ggcmean}{G}
\newcommand{\ggprocess}{\mathbb H_n}
\newcommand{\gglimitprocess}{\mathbb H}

\newcommand{\abound}{K}

\begin{document}

\title{\fontsize{16}{19} Causal and conditional tail coefficients: asymptotics, inference and graph recovery}

\author{
Axel B\"ucher\thanks{Ruhr-Universität Bochum, Fakultät für Mathematik, 44780 Bochum, Germany. Email: \href{mailto:axel.buecher@rub.de}{axel.buecher@rub.de}}
~\orcidlink{0000-0002-1947-1617}
\and 
Mario Krali\thanks{Ruhr-Universität Bochum, Fakultät für Mathematik, 44780 Bochum, Germany. Email: \href{mailto:mario.krali@rub.de}{mario.krali@rub.de}}
~\orcidlink{0000-0002-6203-2804}
}

\date{\today}

\maketitle

\begin{abstract}
The causal tail coefficient (CTC) measures the influence of an extreme value of one variable on another variable and has been proposed as a tool for causal discovery in heavy-tailed models. Despite its popularity, the asymptotic behavior of its empirical counterpart has not been well understood. We develop an asymptotic theory for the rank-based CTC estimator with an intermediate threshold under nonparametric conditions that are satisfied by all standard copula models. A key observation is that the CTC is determined by the mean of a conditional tail limit distribution, which clarifies its interpretation and determines the asymptotic variance of the estimator. We establish joint asymptotic normality for the entire CTC matrix and prove consistency of a multiplier bootstrap for simultaneous inference. For heavy-tailed linear and recursive max-linear structural equation models, we show that the CTC equals one if and only if the first variable is an ancestor of the second, which leads to confidence graphs for the ancestral graph with asymptotic coverage guarantees. All results extend to general conditional tail coefficients, defined as further functionals of the conditional tail limit distribution. A simulation study and an application to dietary intake data illustrate the methods.
\end{abstract}



\tableofcontents

\section{Introduction}

Extreme events seldom occur in isolation. Floods arise from the conjunction of heavy precipitation and saturated soils, financial crises are triggered by the propagation of losses across markets, and, in nutrition, extreme intake of one nutrient tends to accompany extreme intake of others. In each case, the joint occurrence of the events often reflects an underlying causal mechanism. Learning this mechanism faces two main obstacles: the relevant information is limited to the few observations that fall in the tail, and controlled interventions are often costly, unethical or impossible. The mechanism must therefore be inferred from scarce observational data, with substantial uncertainty.

Extreme value theory provides the probabilistic and statistical foundations for the analysis of such events~\citep{resnick1987extreme, resnick2007heavy, beirlant2004statistics, DeHaan2006}, and a growing literature combines it with graphical modelling and causal inference~\citep{Pearl2009}. Two classes of models along a directed acyclic graph (DAG) have received particular attention: recursive max-linear models (RMLMs), introduced by \citet{Gissibl2018a}, in which each variable is the maximum of scaled parent variables and an independent innovation, and heavy-tailed linear structural equation models (linear SEMs), in which the maximum is replaced by a sum \citep{gnecco2021causal}. Learning the graph from data is an active area of research. For RMLMs, \citet{gissibl2021identifiability} propose the minimum ratio estimator, which exploits the exact parametric form of the model, \citet{buck2021recursive} extend it to propagating noise, and \citet{KluKra21} learn the structure from the angular measure of the multivariate extreme value distribution. For heavy-tailed linear models, \citet{gnecco2021causal} estimate a causal ordering by means of the causal tail coefficient (CTC), and \citet{krali2026causal} adapts the angular-measure approach to non-negative coefficients. Alternative approaches are considered in \citet{mhalla2020causal, tran2024estimating, mhalla2025causal, jiang2025separation}. More recently, \citet{engelke2025extremes} and \citet{eSCM_FangBaiWang} proposed more general extremal structural causal models, in which causal relations and interventions are defined at the extreme scale. Undirected graphical models for extremes, building on the notion of conditional independence of \citet{engelke2020graphical}, lack a notion of direction and therefore lie outside the scope of the present work.

Most of these approaches rely on the full multivariate tail or on a specific parametric model. The CTC takes a different route and has become a central tool for causal inference in extremes. For a pair $(X_1,X_2)$, it is denoted by $\Gamma_{12}$ and defined as the limiting mean 
quantile level of $X_2$ given that $X_1$ is extreme; see \eqref{eq:def-ctc} for a precise definition.
Three features explain its appeal. First, it captures the basic asymmetry of causality in extremes: in heavy-tailed structural equation models, extremes propagate from causes to effects but not conversely, so that, under suitable conditions, $\Gamma_{12} = 1$ when $X_1$ causes $X_2$, while $\Gamma_{21} < 1$. Second, it is model-free: it depends only on the copula of $(X_1,X_2)$ and can be estimated from ranks, without fitting a structural model. Third, it is pairwise, so that the matrix of all CTCs, and with it a causal ordering, remains accessible in high dimensions, where the full multivariate tail is difficult to estimate. These features have made the CTC a building block for a range of methods: it has been extended to time series
\citep{bodik2024causality}, to settings with observed confounders \citep{pasche2023causal}, to Granger causality \citep{bodik2024granger}, and to compound extreme events \citep{yin2025causal}.

Despite this central role, all these works draw causal conclusions from estimates that are based on the few observations in the tail, without any quantification of their uncertainty. Moreover, although the CTC is a model-free quantity, its existing theory is tied to the structural model that gives its causal interpretation.

\paragraph{Distributional theory for the CTC.}
\citet{gnecco2021causal} proposed a rank-based estimator of the CTC and established its consistency
within heavy-tailed linear structural equation models. Its distributional properties have only recently attracted attention. 
\citet{zhou2027efficient} state a limit theorem for the empirical CTC under the heavy-tailed linear SEM of \citet{gnecco2021causal}; however, the stated limit involves the divergent integral $\int_0^\infty R(s,1)\diff s$ of the tail copula $R$ and depends on the marginal tail index, although the estimator is rank-based.
\citet{garcin2026directional} and \citet{leimenstoll2026identification} consider related directional coefficients that differ from the CTC and are based on a fixed, rather than an intermediate, threshold. The former build on the relation between the CTC and partial derivatives of the copula and impose smoothness assumptions in the spirit of \cite{Segers2012}, which exclude important models such as max-linear copulas; the latter focus on linear structural equation models for the bulk of the distribution. An asymptotic theory for the empirical CTC itself, under conditions that cover standard models of extremal dependence and with confidence statements about the underlying causal structure, has hence been missing. More broadly, the CTC is a conditional expectation of one variable given that another one is extreme, and our results are thus related to the literature on such tail conditional expectations, for instance the marginal expected shortfall studied by \citet{cai2015estimation}, which, however, is defined on the original scale and requires different assumptions.

\paragraph{Our contribution.}
We develop an asymptotic theory for the empirical CTC with an intermediate threshold, under nonparametric conditions that are satisfied by many standard copula models, including extreme-value and max-linear ones, together with a second-order condition that controls the bias. Our starting point is the observation that the CTC is one minus the mean of a \emph{conditional tail limit distribution}, which also clarifies its interpretation: the CTC is not a measure of tail dependence, and every value in $[0,1]$ can be attained under upper-tail independence. For the estimator of \citet{gnecco2021causal}, we establish asymptotic normality, with asymptotic variance given by the variance of the conditional tail limit distribution; interestingly, this is smaller than for known margins. We then derive joint asymptotic normality for the entire CTC matrix, with an explicit covariance matrix, and prove consistency of a multiplier bootstrap for simultaneous inference. For heavy-tailed linear and recursive max-linear structural equation models, we show that all conditions are met and that both models share a common \emph{ancestral tail mixture structure}, under which the CTC matrix and the asymptotic variances are available in closed form, and the CTC equals one if and only if the first variable is an ancestor of the second. This leads to confidence graphs for the ancestral graph with asymptotic coverage guarantees, via a $\gamma$-signal graph that coincides with the ancestral graph for $\gamma$ close to one. Finally, the main results extend to general \emph{conditional tail coefficients}, defined as expectations of functions of the conditional tail limit distribution, which include the CTC as a special case. A simulation study and an application to dietary intake data illustrate the finite-sample performance.

\paragraph{Organization of the paper.}
Section~\ref{sec:bivariate-CTC} is concerned with the bivariate CTC: Section~\ref{sec:bivariate-foundations} introduces the conditional tail limit distribution and discusses the interpretation of the CTC, and Section~\ref{sec:bivariate-inference} establishes asymptotic normality of the empirical CTC. Joint asymptotics for the entire CTC matrix, including covariance estimation and a multiplier bootstrap, are derived in Section~\ref{sec:multivariate-CTC}. Section~\ref{sec:causal-models} introduces the ancestral tail mixture structure and shows that it is satisfied by heavy-tailed linear and recursive max-linear structural equation models. These results are used in Section~\ref{sec:graph-recovery} to construct confidence graphs for the ancestral graph. Section~\ref{sec:general-tail-quantities} extends the theory to general conditional tail coefficients. The finite-sample performance of the methods is investigated in a simulation study in Section~\ref{sec:simulations}, and Section~\ref{sec:data-application} presents an application to dietary intake data. All proofs are deferred to Appendix~\ref{sec:proofs}, auxiliary results are collected in Appendix~\ref{sec:auxiliary}, and additional simulation results are presented in Appendix~\ref{sec:additional-simulations}.

\section{The bivariate causal tail coefficient}
\label{sec:bivariate-CTC}

\subsection{Definition and interpretation}
\label{sec:bivariate-foundations}

Let $(X_1,X_2)$ be a random vector with continuous marginal CDFs $F_1$ and $F_2$. The \emph{causal tail coefficient} (CTC) of \citet{gnecco2021causal} measures how large $X_2$ is, after standardization to standard uniform margins, when $X_1$ is extreme:
\begin{align} \label{eq:def-ctc}
\Gamma = \Gamma_{12} = \lim_{u \uparrow 1} \Exp\big[ F_2(X_2) \mid F_1(X_1) > u\big],
\end{align}
provided the limit exists. Clearly, $\Gamma \in [0,1]$. We have $\Gamma=1/2$ if $X_1$ and $X_2$ are independent, while $\Gamma=1$ if and only if, conditionally on $F_1(X_1)>u$,  $F_2(X_2)$ converges to one in probability as $u\uparrow1$ (indeed, for bounded random variables, $L^1$-convergence is equivalent to convergence in probability).

For the formal analysis of $\Gamma$, it is convenient to work with $V_j = 1-F_j(X_j)$, which is also uniformly distributed on $(0,1)$ and small when $X_j$ is large. The joint CDF $D$ of $(V_1,V_2)$ is the survival copula of $(X_1,X_2)$, and is related to the copula $C$ of $(X_1,X_2)$ by $D(v_1,v_2) = v_1+v_2-1+C(1-v_1,1-v_2)$. In this notation,
\[
\Gamma = \lim_{h \downarrow 0} \Exp[1-V_2 \mid V_1 < h],
\]
and the conditional distribution $\mathcal L(V_2 \mid V_1 < h)$ has CDF $v \mapsto \Prob(V_2\le v \mid V_1<h) = D(h,v)/h$. The CTC is hence determined by the limiting behavior of these conditional distributions, which motivates the following condition.

\begin{condition}{W}
\label{cond:conditional-weak-limit-bivariate}
Let $D$ be a copula and let $(V_1, V_2) \sim D$.
The conditional distribution $\mathcal L(V_2 \mid V_1<h)$ converges weakly for $h\downarrow0$.
We write the CDF of the limit distribution as $H=H_{12}$, and we write $Z=Z_{12}$ for a random variable with CDF $H$. We refer to this limiting distribution as the \emph{conditional tail limit distribution}.
\end{condition}

Since $V_2$ is bounded, Condition~\ref{cond:conditional-weak-limit-bivariate} implies that the CTC exists and is given by $\Gamma = 1-\Exp[Z]$. Moreover, since the conditional CDF of $V_2$ given $V_1<h$ is $v \mapsto D(h,v)/h$, Condition~\ref{cond:conditional-weak-limit-bivariate} can be expressed in terms of the right-hand partial derivative of $D$ with respect to its first argument at the boundary point $0$, as made precise in the following lemma.

\begin{lemma} \label{lem:gamma-representation-bivariate}
Let $D$ be a copula and let $(V_1, V_2) \sim D$. Then
Condition \ref{cond:conditional-weak-limit-bivariate} holds if and only if the right-hand partial derivative of $D$ with respect to its first argument at $(0,v)$,
\begin{align}
\label{eq:pd0-bivariate}
\partial_1^+ D(0,v)  := \lim_{h \downarrow 0} \frac{D(h,v)}{h},
\end{align}
exists for Lebesgue-almost every $v \in (0,1)$. In that case, we have $H(v) = \partial_1^+ D(0,v)$ for Lebesgue-almost every $v \in (0,1)$, and it holds
\[
\Gamma  = 1- \Exp_{Z \sim H}[Z] = \int_0^1  \partial_1^+ D(0,v)  \, \diff  v.
\]
\end{lemma}

We illustrate the generality of Condition~\ref{cond:conditional-weak-limit-bivariate} by two examples. Further examples on linear and max-linear structural equation models are discussed in Section~\ref{sec:causal-models}.

\begin{example}
\label{ex:evc-conditions}
Suppose $C$ is an extreme value copula with Pickands dependence function $A$, i.e.,
\[
C(u,v) = \exp\Big\{ (\log u+\log v) A\Big(\frac{\log v}{\log u+\log v}\Big)\Big\}, \qquad (u,v)\in(0,1)^2,
\]
and let $D$ be the associated survival copula. A straightforward calculation shows that $\partial_1^+ D(0,v) = v+A'(1-)(1-v)$ for $v\in(0,1]$, where $A'(1-)$ denotes the left-hand derivative of $A$ at $1$. Hence, by Lemma~\ref{lem:gamma-representation-bivariate}, Condition~\ref{cond:conditional-weak-limit-bivariate} holds, with conditional tail limit distribution $A'(1-)\,\delta_0 + \{1-A'(1-)\}\,\mathrm{Unif}(0,1)$, where $\delta_x$ denotes the Dirac measure at $x$, and $\Gamma = \{1+A'(1-)\}/2$.

Alternatively, $\Gamma$ can be expressed in terms of the spectral measure of $A$ with respect to the $\|\cdot\|_1$-norm, that is, the unique probability measure $\Phi$ on $[0,1]$ with mean $1/2$ such that $A(t) = 2 \int_{[0,1]} \max\{tw, (1-t)(1-w)\} \, \Phi(\diff w)$ for $t\in[0,1]$; see \citet[Lemma 2.1]{Boulin2026}. A dominated convergence argument shows that $A'(1-) = 1-2\Phi(\{0\})$; see the appendix. Consequently, $\Gamma = 1-\Phi(\{0\})$, such that any extreme value copula whose spectral measure does not have an atom at zero has CTC equal to $1$.
\end{example}

\begin{example}
\label{ex:gauss-conditions}
Consider the bivariate Gaussian copula $C=C_\rho$ with correlation
parameter $\rho\in(-1,1)$. The associated survival copula satisfies
$D_\rho=C_\rho$, and, for every $v\in(0,1)$,
\[
\partial_1^+ D_\rho(0,v)
=
\begin{cases}
0, & \rho<0,\\
v, & \rho=0,\\
1, & \rho>0.
\end{cases}
\]
Hence, by Lemma~\ref{lem:gamma-representation-bivariate}, Condition~\ref{cond:conditional-weak-limit-bivariate} holds, with conditional tail limit distribution $\delta_1$ if $\rho<0$, $\mathrm{Unif}(0,1)$ if $\rho=0$, and $\delta_0$ if $\rho>0$. Consequently, $\Gamma=0$ for $\rho<0$, $\Gamma=1/2$ for $\rho=0$, and $\Gamma=1$ for $\rho>0$.
\end{example}

The previous examples show that the CTC is not a functional of the extremal dependence structure. Rather, it depends on the entire conditional distribution of $V_2$ given that $V_1$ is small, including the non-extreme part of $X_2$. In particular, $\Gamma$ is not determined by the limiting extreme value copula of a distribution in its max-domain of attraction: for $\rho\ne0$, the Gaussian copula from Example~\ref{ex:gauss-conditions} belongs to the max-domain of attraction of the independence copula, whose Pickands dependence function $A\equiv1$ would suggest $\Gamma=1/2$, while in fact $\Gamma\in\{0,1\}$. For the same reason, the CTC should not be interpreted as a measure of tail dependence. The Gaussian copula with $\rho>0$ is upper-tail independent, while $\Gamma=1$. The following lemma shows that, in fact, every value of the CTC in $[0,1]$ can be attained under upper-tail independence.

\begin{lemma} \label{lem:ctc-asymptotic-independence}
For every $p\in [0,1]$, there exists an absolutely continuous bivariate
copula $C$ that is upper-tail independent and whose causal tail coefficient
satisfies $\Gamma=p$.
\end{lemma}

In view of these examples and Lemma~\ref{lem:ctc-asymptotic-independence}, the CTC might as well be called a \emph{conditional tail coefficient}: it is a functional of the conditional tail limit distribution, and its causal interpretation arises from structural assumptions, as discussed in Section~\ref{sec:causal-models}. In this paper, we keep the established name \emph{causal tail coefficient} for $\Gamma$, but refer to the more general functionals of the conditional tail limit distribution studied in Section~\ref{sec:general-tail-quantities} as \emph{conditional tail coefficients}.

\subsection{Estimation and asymptotic normality}
\label{sec:bivariate-inference}

A natural estimator of $\Gamma$ is obtained by replacing the unknown marginal CDFs in the definition of the CTC by their empirical counterparts and the limit $u \uparrow 1$ by a high empirical quantile. This estimator was proposed by \citet{gnecco2021causal}, who established its consistency under specific model assumptions that we take up again in Section~\ref{sec:linear-scm}. In this section, we go two steps further. First, we work under general, nonparametric conditions on the survival copula $D$ rather than under specific model assumptions. Second, we go beyond consistency and show that the estimator is asymptotically normal, which provides the basis for statistical inference.

Let $\bm X_1, \dots, \bm X_n$ be a random sample from $\bm X=(X_1, X_2)^\top$. For $j \in \{1,2\}$ and $i\in[n]$, write $R_{ij}$ for the rank of $X_{ij}$ among $X_{1j}, \dots, X_{nj}$, and let
\[
\widehat U_{ij} = \frac{R_{ij}}{n+1}, \qquad \widehat V_{ij} = 1-\frac{R_{ij}}{n+1}.
\]
The sample version of the CTC is then given by
\[
\widehat\Gamma_n = \widehat \Gamma_{n,12}
= 
\frac1k \sum_{i=1}^n (1-\widehat V_{i2} )\bm 1\Big( \widehat V_{i1} \le \frac{k}{n} \Big)
= 
\frac1k \sum_{i=1}^n \widehat U_{i2}  \bm 1\Big( \widehat V_{i1} \le \frac{k}{n} \Big),
\]
where $k \in [n]$ is some threshold chosen by the statistician (up to the normalization of the ranks, this is exactly the estimator proposed in \citealp{gnecco2021causal}). Since there are almost surely no ties, exactly $k$ observations satisfy $\widehat V_{i1}\le k/n$, which explains the normalization by $k^{-1}$. Asymptotically, we will require $k=k_n$ being an intermediate sequence, that is, $k=k_n \to \infty$ such that $k=o(n)$ for $n \to \infty$.

The estimator $\widehat\Gamma_n$ does not directly estimate the limit $\Gamma$, but rather its counterpart at the finite threshold $k/n$. To make this precise, let
\begin{align*}G(z) = \Exp\big[1-V_2 \mid V_1 \le z \big]
= \int_0^1 \frac{D(z, v)}z \, \diff v, \qquad z \in (0,1],
\end{align*}
such that $\Gamma = \lim_{z \downarrow 0} G(z)$, and define the \emph{pre-asymptotic CTC} as
$
    \Gamma_n = G(k/n).
$
The name reflects that $\Gamma_n$ is the value of the CTC before passing to the limit, evaluated at the threshold used by the estimator. Indeed, if the margins were known, the empirical CTC $\widehat \Gamma_n^*$ from \eqref{eq:knownCTC} below would be an unbiased estimator of $\Gamma_n$. Accordingly, we decompose the estimation error as
\[
\widehat\Gamma_n - \Gamma = (\widehat\Gamma_n - \Gamma_n) + (\Gamma_n - \Gamma),
\]
into a stochastic part and a deterministic approximation error. We first treat the stochastic part and then the bias.

Condition~\ref{cond:conditional-weak-limit-bivariate} only concerns $D(u,v)$ at the boundary $u=0$. Asymptotic normality of $\widehat\Gamma_n$ requires a finer analysis: the estimator uses those observations with $\widehat V_{i1} \le k/n$, which corresponds to $V_{i1} \le \widehat x_n k/n$ for some random $\widehat x_n$ with $\widehat x_n - 1 = O_\Prob(k^{-1/2})$. The effect of this random threshold is of the same order as the estimation error itself and is governed by $\partial_1 D(u,v)$ for $u$ in a neighborhood of $0$. We therefore require an additional regularity condition on $D$, which ensures that this effect is asymptotically linear.
Recall that, for every $v \in [0,1]$, the function $u \mapsto D(u,v)$ is Lip\-schitz continuous, and hence the first partial derivative $\partial_1 D(u,v)$ exists for Lebesgue-almost every $u \in (0,1)$. The following condition strengthens this property near the boundary by requiring the corresponding right-hand partial derivative to exist for all sufficiently small $u$, and to be continuous at $u=0$.

\begin{condition}{D}
\label{cond:continuity2-bivariate}
Let $D$ be a copula.
For almost every $v \in (0,1)$, the right-hand partial derivative
\[\partial_1^+ D(u,v)
=
\lim_{h \downarrow 0}
\frac{D(u+h,v)-D(u,v)}{h}
\]
exists for all $u$ in some sufficiently small interval $[0,\eps_v)$. Moreover,
$\lim_{u \downarrow 0} \partial_1^+ D(u,v)
 = \partial_1^+ D(0,v)$.
\end{condition}

Note that $\partial_1^+D(u,v)$ admits the interpretation
\[
\partial_1^+D(u,v)
=
\lim_{h\downarrow0}
\Pr(V_2\le v \mid V_1\in [u, u+h]).
\]
It further follows from Lemma~\ref{lem:gamma-representation-bivariate} that Condition~\ref{cond:continuity2-bivariate} implies Condition~\ref{cond:conditional-weak-limit-bivariate}, since Condition~\ref{cond:continuity2-bivariate} ensures the existence of $\partial_1^+D(0,v)$ for almost every $v\in(0,1)$.
We briefly address Condition~\ref{cond:continuity2-bivariate} in the context of the the examples of Section~\ref{sec:bivariate-foundations}.

\begin{example}[Example~\ref{ex:evc-conditions} continued]
\label{ex:evc-D}
The survival copula $D$ of every bivariate extreme value copula satisfies Condition~\ref{cond:continuity2-bivariate}; see the appendix.
\end{example}

\begin{example}[Example~\ref{ex:gauss-conditions} continued]
\label{ex:gauss-D}
The survival copula $D_\rho=C_\rho$ of the bivariate Gaussian copula from Example~\ref{ex:gauss-conditions} satisfies Condition~\ref{cond:continuity2-bivariate} for every $\rho\in(-1,1)$; see the appendix.
\end{example}

The following theorem treats the stochastic part of the estimation error. It only requires Condition~\ref{cond:continuity2-bivariate} and no rate condition on $k$ beyond $k\to\infty$ and $k=o(n)$.

\begin{theorem}
\label{theo:CTC-bivariate} Let $k \in [n]$ with $k \to \infty$ and $k=o(n)$ as $n \to \infty$.
Under Condition~\ref{cond:continuity2-bivariate}, we have, 
    \begin{align} \label{eq:linear-bivariate}
    \sqrt{k}(\widehat \Gamma_n - \Gamma_n)
    &=
    \frac1{\sqrt n} \sum_{i=1}^n  \phi_n(V_{i1},V_{i2})  + o_\Prob(1)
    \end{align}
    for $n \to \infty$, where
    \[
    \phi_n(v_1,v_2)
    =  \sqrt{\frac{n}{k}} \big( 1 - v_2 - \Gamma_n \big) \bm 1\Big(v_{1} \le  \frac{k}{n} \Big).
    \]
    Moreover, for $Z \sim H$ with $H$ from Condition~\ref{cond:conditional-weak-limit-bivariate}, 
    \begin{align} \label{eq:weak-bivariate} 
    \sqrt{k}(\widehat \Gamma_n - \Gamma_n) \weak \Nc(0, \sigma^2), \qquad \sigma^2 = \Var(Z).
    \end{align}
\end{theorem}

The asymptotic linearization in \eqref{eq:linear-bivariate} expresses the stochastic part of the estimation error as a sum of independent terms, up to a negligible remainder. Beyond the asymptotic normality in \eqref{eq:weak-bivariate}, it is the key ingredient for the joint asymptotics of several CTCs and for the validity of the bootstrap in Section~\ref{sec:multivariate-CTC}.
Theorem~\ref{theo:CTC-bivariate} is a special case of more general results on empirical conditional tail coefficients, that is, on empirical versions of $\Exp[g(V_2) \mid V_1 \le h]$, which are developed in Section~\ref{sec:general-tail-quantities} and which also include a functional version with respect to the threshold and a multiplier bootstrap. Bootstrap approximations, which are particularly useful for simultaneous inference on several CTCs, are discussed in Sections~\ref{sec:multivariate-CTC} and~\ref{sec:general-tail-quantities}.

Interestingly, the use of ranks in the definition of the empirical CTC leads to a reduction in asymptotic variance. Indeed, if the marginal CDFs $F_1$ and $F_2$ were known, one could alternatively consider the known-margin empirical CTC
\begin{align}\label{eq:knownCTC}
    \widehat \Gamma_n^* = \frac1k \sum_{i=1}^n (1-V_{i2}) \bm 1(V_{i1} \le k/n).
\end{align}
It follows from Proposition~\ref{prop:known-margins-bivariate-2param} that the asymptotic variance of $\widehat \Gamma_n^*$ is
$
\sigma^2_{\mathrm{known}} = \Exp[(1-Z)^2].
$
In contrast, the rank-based estimator has asymptotic variance
\[
\sigma^2 = \Var(Z) = \Var(1-Z) = 
\Exp[(1-Z)^2]-\{\Exp[1-Z]\}^2
=
\sigma^2_{\mathrm{known}} - \Gamma^2.
\]
Consequently, $\sigma^2\le \sigma^2_{\mathrm{known}}$, with equality if and only if $\Gamma=0$. Thus, estimating the marginal distributions rather than treating them as known can reduce the asymptotic variance of the empirical CTC; see also \cite{GenestSegers2010, Buecher2014} for related phenomena.

For statistical inference on $\Gamma$ itself, the deterministic approximation error $\Gamma_n-\Gamma$ must be negligible compared to the stochastic part, which is of order $k^{-1/2}$. This reflects the usual bias--variance trade-off in extreme value statistics: a larger $k$ reduces the variance, but increases the bias, since the threshold $k/n$ is then less extreme. Formally, we need the convergence $G(z) \to \Gamma$ for $z \downarrow 0$ to be sufficiently fast relative to the growth of $k$, which is a second-order condition on $D$, in the spirit of second-order conditions in classical extreme value theory \citep[see, e.g.,][]{DeHaan2006}.

\begin{condition}{B}
\label{cond:bias-bivariate}
The sequence $k=k_n$ and the copula $D$ satisfy $\sqrt k B(n/k) = o(1)$ as $n \to \infty$, where $B:(0,\infty) \to \R$ is defined as
\[
B(t)
= G(1/t) - \Gamma = \int_0^1 \big\{ tD(1/t,v) - \partial_1^+D(0,v) \big\}  \, \diff v
\]
\end{condition}

\begin{example}[Example~\ref{ex:evc-conditions} continued]
\label{ex:evc-bias}
Let $C$ be an extreme value copula with Pickands dependence function $A$, and suppose additionally that $A$ admits an expansion of the form 
\[
A(1-h) = 1 - ah + bh^2   + o(h^2)
\]
for $h \downarrow 0$, for some $a,b\ge 0$. We then have
$a=A'(1-)$ and $b= A''(1-)/2$, with $A''(1-)$ the left-hand derivative at $1$ of $x \mapsto A'(x -)$. 
We show in the appendix that 
\[
B(t) =
\int_0^1 \big\{ tD(1/t,v) - \partial_1^+D(0,v) \big\}  \, \diff v
=
- b \frac{\log t}t + o\Big(\frac{\log t}t \Big), \qquad t \to \infty.
\]
The required rate on $k$ from Condition~\ref{cond:bias-bivariate} is therefore
$k^{3/2} \log(n/k) = o(n)$ as $n \to \infty$.
If $A$ is even affine in a left neighbourhood of $1$, that is, $A(t)=1-a(1-t)$ for all $t \in [t_0,1]$ and some $t_0 \in [0,1)$, then the rate improves to $B(t)=O(1/t)$ as $t \to \infty$, such that Condition~\ref{cond:bias-bivariate} holds whenever $k=o(n^{2/3})$; see the appendix. This is the case, for instance, for the max-linear models from Section~\ref{sec:max-linear-scm}.
\end{example}

\begin{corollary}
\label{cor:CTC-bivariate-bias}
In the setting of Theorem~\ref{theo:CTC-bivariate}, assume additionally that Condition \ref{cond:bias-bivariate} holds. Then, for $n \to \infty$,
\[
\sqrt k ( \widehat \Gamma_n - \Gamma) \weak \Nc(0,\sigma^2).
\]
\end{corollary}

\begin{proof} 
Condition \ref{cond:bias-bivariate} implies
$
\sqrt k \big(\Gamma_n - \Gamma \big)  
=
\sqrt k \big( G(k/n) - \Gamma \big)
=
\sqrt k B(n/k) = o(1).
$
The claim then follows from Theorem~\ref{theo:CTC-bivariate} and Slutsky's lemma.
\end{proof}

Provided $\sigma^2>0$, the corollary yields asymptotic confidence intervals for $\Gamma$ of the form $\widehat\Gamma_n \pm z_{1-\alpha/2}\,\widehat\sigma_n/\sqrt k$, with $z_{1-\alpha/2}$ the standard normal quantile and $\widehat\sigma_n$ a consistent estimator of $\sigma$. A natural choice is the empirical variance of the $k$ selected values $\widehat U_{i2}$,
\[
\widehat\sigma_n^2 = \frac1k \sum_{i=1}^n \big(\widehat U_{i2} - \widehat\Gamma_n\big)^2 \bm 1\Big( \widehat V_{i1} \le \frac{k}{n} \Big),
\]
whose consistency under Condition~\ref{cond:continuity2-bivariate} follows from Corollary~\ref{cor:unknown-margins-bivariate} applied with $g(v)=(1-v)^2$; indeed, $\widehat\sigma_n^2 = k^{-1}\sum_{i=1}^n \widehat U_{i2}^2\, \bm 1(\widehat V_{i1} \le k/n) - \widehat\Gamma_n^2$, where the first term is the empirical conditional tail coefficient from Section~\ref{sec:general-tail-quantities} with that choice of $g$. Condition~\ref{cond:bias-bivariate} is further verified for the structural equation models in Section~\ref{sec:causal-models}.

\section{The causal tail coefficient matrix}
\label{sec:multivariate-CTC}

Throughout, $\bm X=(X_1, \dots, X_d)^\top$ is a random vector with continuous marginal CDFs $F_1, \dots, F_d$. As in Section~\ref{sec:bivariate-CTC}, let $V_j = 1-F_j(X_j)$, write $\bm V=(V_1, \dots, V_d)^\top$, and let $D_{j\ell}$ denote the survival copula of $(X_j,X_\ell)$, that is, the CDF of $(V_j,V_\ell)$. For each ordered pair $(j,\ell)$ in
$\Pc_2 = \{(j,\ell) \in [d]^2: j \ne \ell\}$, the CTC $\Gamma_{j\ell}$ from $X_j$ to $X_\ell$ is the bivariate CTC of Section~\ref{sec:bivariate-CTC} applied to $(X_j,X_\ell)$, that is, $\Gamma_{j\ell} = \lim_{h \downarrow 0} \Exp[1-V_\ell \mid V_j < h]$. We collect these coefficients in the \emph{CTC matrix} $(\Gamma_{j\ell})_{(j,\ell)\in\Pc_2}$~\citep{gnecco2021causal}.

Let $\bm X_1, \dots, \bm X_n$ be a random sample from $\bm X$, and let $\widehat V_{ij} = 1-R_{ij}/(n+1)$, where $R_{ij}$ is the rank of $X_{ij}$ among $X_{1j}, \dots, X_{nj}$. Applying the estimator from Section~\ref{sec:bivariate-inference} to each pair of coordinates yields
\[
\widehat\Gamma_{j\ell} = \frac1k \sum_{i=1}^n (1-\widehat V_{i\ell})\, \bm 1\Big( \widehat V_{ij} \le \frac{k}{n} \Big), \qquad (j,\ell)\in\Pc_2,
\]
with corresponding pre-asymptotic CTC $\Gamma_{j\ell}^{(n)} = (n/k)\,\Exp[(1-V_\ell)\, \bm 1(V_j \le k/n)]$, the analogue of $\Gamma_n$ from Section~\ref{sec:bivariate-inference}. Throughout, $k=k_n \to \infty$ with $k=o(n)$.

The marginal asymptotics of each $\widehat\Gamma_{j\ell}$ follow from Theorem~\ref{theo:CTC-bivariate}. The joint asymptotics additionally require controlling the dependence between estimators with different conditioning variables, which we capture via the following condition on the joint tail of $\bm V$.

\begin{condition}{M}
\label{cond:multivariate}
For every nonempty $I\subset[d]$ with $|I|\le2$,
there exists a finite measure $\nu_I$ on $[0,1]^d$
such that
\[
    \nu_{I,h} := h^{-1}\Prob\big(
\bm V\in \cdot,\ V_j\le h\text{ for all }j\in I
\big)\weak\nu_I,
    \qquad h\downarrow0.
\]
\end{condition}

If $I=\{j\}$ is a singleton, the measure $\nu_j:=\nu_{\{j\}}$ is a probability measure. In that case, let $\bm Z_j=(Z_{j1},\dots,Z_{jd})^\top$ denote a random vector with distribution $\nu_j$. Condition~\ref{cond:multivariate} implies that each bivariate copula $D_{j\ell}$ satisfies Condition~\ref{cond:conditional-weak-limit-bivariate}, with the distribution of $Z_{j\ell}$ defined therein coinciding with the corresponding marginal distribution of $\bm Z_j$.

If $I=\{j,s\}$ has cardinality two, the measure $\nu_I$ is a subprobability measure on $[0,1]^d$, with total mass equal to the classical tail dependence coefficient of $(X_j,X_s)$, that is, 
\[
    \nu_{\{j,s\}}([0,1]^d)
    =
    \chi_{js}
    :=
    \lim_{h\downarrow0}
    \Prob(V_j\le h\mid V_s\le h).
\]

\begin{theorem}
\label{theo:CTC-multivariate}
Let $k \in [n]$ with $k \to \infty$ and $k=o(n)$ as $n \to \infty$.
Suppose that, for every $(j,\ell)\in\Pc_2$, the copula $D_{j\ell}$ satisfies Condition~\ref{cond:continuity2-bivariate}.
Then, for every $(j,\ell)\in\Pc_2$,
    \begin{align*}    \sqrt{k}(\widehat \Gamma_{j\ell} - \Gamma_{j\ell}^{(n)}) 
    &=
    \frac1{\sqrt n} \sum_{i=1}^n \phi_{j\ell}^{(n)}(V_{ij},V_{i\ell})  + o_\Prob(1)
    \end{align*}
    for $n \to \infty$, where
    \begin{align} 
    \label{eq:phin-multivariate}
    \phi_{j\ell}^{(n)}(v_j,v_\ell)
    =  \sqrt{\frac{n}{k}} \big( 1 - v_\ell - \Gamma_{j\ell}^{(n)} \big) \bm 1\Big(v_{j} \le  \frac{k}{n} \Big).
    \end{align}
    Moreover, with $Z_{j\ell} \sim H_{j\ell}$ and $H_{j\ell}$ from Condition~\ref{cond:conditional-weak-limit-bivariate} applied to $D_{j\ell}$,
    \begin{align} \label{eq:weak-CTC-multivariate}
    \Gb_{n,j\ell} := \sqrt{k}(\widehat \Gamma_{j\ell} - \Gamma_{j\ell}^{(n)}) \weak \Nc(0, \sigma_{j\ell}^2), 
    \qquad \sigma_{j\ell}^2
    = \Var(Z_{j\ell}).
    \end{align}
    If, additionally, Condition~\ref{cond:multivariate} is met, then
    \begin{align*}
    \Gb_n := \big(\Gb_{n,j\ell} \big)_{(j,\ell) \in \mathcal P_2}
    \weak
    \mathcal N_{d(d-1)} (0, \Sigma)
    \end{align*}
    where $\Sigma$ has entries 
    \begin{align} \label{eq:covariance-identity}
\sigma_{j\ell;st}
=
\int_{[0,1]^d}
(1-v_\ell-\Gamma_{j\ell})
(1-v_t-\Gamma_{st})
\,\nu_{\{j,s\}}(\mathrm d\bm v).
    \end{align}
\end{theorem}

\begin{remark}
The covariance formula in \eqref{eq:covariance-identity}
simplifies in several interesting special cases.
First, for $j=s$, we obtain
\[
    \sigma_{j\ell;jt}
    = \Cov(Z_{j\ell},Z_{jt}),
    \qquad \ell,t\ne j.
\]
Thus, the asymptotic covariance matrix of the CTC estimators
with the same conditioning variable $X_j$ is given by
the covariance matrix of the corresponding components of $\bm Z_j$.

Second, for $j\ne s$ with $\chi_{js}=0$, we have
\[
    \sigma_{j\ell;st}=0,
    \qquad \ell\ne j,\ t\ne s.
\]
Hence, asymptotic independence of the conditioning variables
$X_j$ and $X_s$ implies asymptotic independence of the
corresponding rows of the estimated CTC matrix.

Finally, for $j\ne s$, the covariance between the estimators
of the CTC in opposite directions takes the particularly
simple form
\[
    \sigma_{js;sj}
    = \chi_{js}(1-\Gamma_{js})(1-\Gamma_{sj}).
\]
Indeed, this follows from \eqref{eq:covariance-identity}
since $\nu_{\{j,s\}}$ is concentrated on
$\{\bm v\in[0,1]^d:v_j=v_s=0\}$.
In particular, this covariance is always nonnegative.
\end{remark}

Statistical inference on the entire CTC matrix, for instance by simultaneous confidence intervals or by the confidence graphs in Section~\ref{sec:graph-recovery}, requires two ingredients: an estimator of the limiting covariance matrix $\Sigma$, and an approximation of the distribution of non-linear functionals of $\Gb_n$, such as $\max_{(j,\ell)\in\Pc_2} |\Gb_{n,j\ell}|$. Both can be built from the asymptotic linearization in Theorem~\ref{theo:CTC-multivariate}, by replacing the unknown influence functions $\phi_{j\ell}^{(n)}$ by plug-in estimates. For the second ingredient, we use a multiplier bootstrap, in which the estimated influence functions are multiplied by IID random weights that are independent of the data. Conditionally on the data, this mimics the covariance structure of the linearization, without the need to resample the data or to recompute any ranks. Multiplier bootstrap methods are well established in empirical process theory, see, e.g., \citet[Chapters~10 and~11]{Kosorok2008}, and have been used for tail dependence functionals in \citet{BuecherDette2011}. Below, $\cweak{\xi}$ denotes weak convergence conditionally on the data, in probability, with respect to the multipliers $\xi_1, \xi_2, \dots$; see \citet[Section~2.2.3]{Kosorok2008} for a formal definition.

\begin{theorem}[Covariance estimation and multiplier bootstrap for the CTC matrix]
\label{theo:bootstrap-multivariate}
Let $k \in [n]$ with $k\to\infty$ and $k=o(n)$
as $n\to\infty$.
Suppose that Condition~\ref{cond:multivariate} holds and
that, for every $(j,\ell)\in\Pc_2$, the survival copula $D_{j\ell}$
satisfies Condition~\ref{cond:continuity2-bivariate}.
Then the estimator $\widehat\Sigma_n$ with entries
\[
    \widehat\sigma_{j\ell;st}^{(n)}
    =
    \frac1n\sum_{i=1}^n
    \widehat\phi_{j\ell}^{(n)}
    (\widehat V_{ij},\widehat V_{i\ell})
    \widehat\phi_{st}^{(n)}
    (\widehat V_{is},\widehat V_{it}),
    \qquad (j,\ell),(s,t)\in\Pc_2,
\]
where
\[
    \widehat\phi_{j\ell}^{(n)}(v_j,v_\ell)
    =
    \sqrt{\frac nk}
    (1-v_\ell-\widehat\Gamma_{j\ell})
    \bm 1\Big(v_j\le\frac kn\Big),
\]
is consistent for $\Sigma$ defined
in \eqref{eq:covariance-identity}, that is, $\widehat\Sigma_n =\Sigma + o_\Prob(1)$.
Moreover, if $\xi_1,\xi_2,\dots$ are IID centered
random variables with variance $1$, independent
of the observations, then
\[
    \widehat \Gb_{n,j\ell}
    =
    \frac1{\sqrt n}\sum_{i=1}^n
    \xi_i\,
    \widehat\phi_{j\ell}^{(n)}
    (\widehat V_{ij},\widehat V_{i\ell}),
    \qquad (j,\ell)\in\Pc_2,
\]
satisfies
\[
    \widehat \Gb_n := \big(\widehat \Gb_{n,j\ell}\big)_{(j,\ell)\in\Pc_2}
    \cweak{\xi}
    \mathcal N_{d(d-1)}(0,\Sigma).
\]
\end{theorem}

It is worthwhile to mention that the estimated influence functions are empirically centered, i.e.,
\[
\frac1n\sum_{i=1}^n
\widehat\phi_{j\ell}^{(n)}
(\widehat V_{ij},\widehat V_{i\ell})
=0,
\]
since exactly $k$ observations satisfy $\widehat V_{ij} \le k/n$. Hence, no additional centering of the bootstrap replicates is required. Moreover, the diagonal entries $\widehat\sigma_{j\ell;j\ell}^{(n)}$ of $\widehat\Sigma_n$ coincide with the plug-in variance estimator from Section~\ref{sec:bivariate-inference}, applied to the pair $(X_j,X_\ell)$.

Theorem~\ref{theo:bootstrap-multivariate} provides the basis for various simultaneous inference procedures. For instance, if Condition~\ref{cond:bias-bivariate} holds for every $D_{j\ell}$ and $\Sigma \ne 0$, then the conditional $(1-\alpha)$-quantile $\widehat q_{1-\alpha}$ of $\max_{(j,\ell)\in\Pc_2}|\widehat\Gb_{n,j\ell}|$ yields simultaneous confidence intervals $\widehat\Gamma_{j\ell} \pm \widehat q_{1-\alpha}/\sqrt k$ for all entries of the CTC matrix. Likewise, by the continuous mapping theorem, Theorem~\ref{theo:CTC-multivariate} implies asymptotic normality of the difference $\widehat\Gamma_{js}-\widehat\Gamma_{sj}$ between the two causal directions, with asymptotic variance $\sigma_{js}^2+\sigma_{sj}^2-2\chi_{js}(1-\Gamma_{js})(1-\Gamma_{sj})$, which can be used to test for asymmetry of the CTC. A particularly relevant application is the recovery of causal graphs from the CTC matrix, which we address in Section~\ref{sec:graph-recovery}. There, we use a studentized version of the maximum statistic, which accounts for the fact that some of the variances $\sigma_{j\ell}^2$ may be zero, for instance when $\Gamma_{j\ell}=1$. The finite-sample performance of the covariance estimator and of the bootstrap is investigated in Section~\ref{sec:simulations}.

\section{Structural equation models}
\label{sec:causal-models}

We now discuss two classes of structural equation models (SEMs) on directed acyclic graphs for which the causal tail coefficient admits a transparent causal interpretation: heavy-tailed linear SEMs (Section~\ref{sec:linear-scm}) and recursive max-linear models (RMLMs) (Section~\ref{sec:max-linear-scm}). 
The two models are known to share a similar extremal dependence structure which is characterised by discrete atoms. Interestingly, they also share the same conditional tail structure: conditionally on $X_j$ being extreme, the vector $\bm V$ converges to a mixture of distributions indexed by the ancestors of $j$, with mixture weights given by normalized ancestral coefficients. We formalize this structure in Section~\ref{sec:ancestral-mixture} and derive its consequences for the CTC matrix and for the limiting covariance matrix of its estimator. In particular, $\Gamma_{j\ell}=1$ if and only if $j$ is an ancestor of $\ell$, which forms the basis for the graph recovery procedure in Section~\ref{sec:graph-recovery}. For both model classes, it then suffices to verify the tail structure, together with Conditions~\ref{cond:continuity2-bivariate} and~\ref{cond:bias-bivariate}.

\subsection{Graph notation and a common tail structure}
\label{sec:ancestral-mixture}

Let $\mathcal G$ be a directed acyclic graph (DAG) on node set $[d]$. For $j \in [d]$, we write $\pa(j)$ for the set of parents of $j$, $\an(j)$ for the set of ancestors of $j$, that is, the set of nodes $m \ne j$ from which there is a directed path to $j$, and $\An(j)=\an(j)\cup\{j\}$. Moreover, for $r\in[d]$, we write
$
    \Des(r)=\{m\in[d]:r\in\An(m)\}
$
for the set of descendants of $r$, including $r$ itself.

For $r\in[d]$, define the random vector
$\bm W^{(r)}=(W_1^{(r)},\dots,W_d^{(r)})^\top$ by
\[
    W_m^{(r)}
    =
    \begin{cases}
        0, & m\in\Des(r),\\
        V_m, & m\notin\Des(r),
    \end{cases}
    \qquad m\in[d],
\]
and let $\mu_r$ denote its distribution on $[0,1]^d$. Thus, $\bm W^{(r)}$ is obtained from $\bm V$ by setting all coordinates that are descendants of $r$ to zero.

\begin{definition}[Ancestral tail mixture]
\label{def:ancestral-mixture}
Let $w=(w_{jr})_{j,r\in[d]}$ be non-negative weights such that $w_{jr}>0$ if and only if $r \in \An(j)$, and such that $\sum_{r\in[d]} w_{jr}=1$ for every $j \in [d]$. We say that $\bm X$ has an \emph{ancestral tail mixture structure} with weights $w$ if Condition~\ref{cond:multivariate} holds with
\begin{align}
    \nu_j
    &=
    \sum_{r\in\An(j)}
    w_{jr}\mu_r,
    \qquad
    \nu_{\{j,s\}}
    =
    \sum_{r\in\An(j)\cap\An(s)}
    (w_{jr}\wedge w_{sr})\mu_r,
    \qquad j,s \in [d], j \ne s.
    \label{eq:ancestral-mixture}
\end{align}
\end{definition}

The measure $\nu_j$ admits the following alternative
representation. Let $R_j$ be a random variable taking
values in $\An(j)$, independent of $\bm V$, with
$\Prob(R_j=r)=w_{jr}$, and define
\[
    \widetilde V_m^{(j)}
    =
    W_m^{(R_j)}
    =
    V_m\bm 1\{m\notin\Des(R_j)\},
    \qquad m\in[d].
\]
Then
$
    \nu_j=\mathcal L(\widetilde{\bm V}^{(j)}).
$
This representation has a natural interpretation.
Conditionally on $X_j$ being extreme, the extreme event asymptotically originates at a randomly chosen ancestor $R_j\in\An(j)$, which equals $r$ with probability $w_{jr}$.
All descendants of $R_j$ then become extreme as well, such that
their survival probability transforms converge to zero,
while the remaining coordinates retain their joint
distribution. In the structural equation models below, in which $\bm X$ is a function of independent innovations $\varepsilon_1, \dots, \varepsilon_d$ and $X_j$ depends on $\varepsilon_r$ if and only if $r \in \An(j)$, the event $\{R_j=r\}$ corresponds to a single large innovation $\varepsilon_r$, in the spirit of the single big jump principle for heavy-tailed random variables.

\begin{proposition}
\label{prop:ancestral-mixture}
Suppose that $\bm X$ has an ancestral tail mixture structure with weights $w$, and, for $(j,\ell) \in \Pc_2$, let
\[
    p_{j\ell}
    =
    \sum_{r\in\An(j)\cap\An(\ell)}w_{jr}.
\]
Then, for every $(j,\ell)\in\Pc_2$, the bivariate survival copula $D_{j\ell}$ satisfies Condition~\ref{cond:conditional-weak-limit-bivariate}, with conditional tail limit distribution
\[
    \mathcal L(V_\ell\mid V_j < h)
    \weak
    p_{j\ell}\delta_0
    +(1-p_{j\ell})\mathrm{Unif}(0,1),
    \qquad h\downarrow0.
\]
Consequently,
\begin{align*}
    \Gamma_{j\ell}
    =
    \frac{1+p_{j\ell}}2,
    \qquad
    \sigma_{j\ell}^2
    =
    \Var(Z_{j\ell})
    =
    \frac{(1-p_{j\ell})(1+3p_{j\ell})}{12}.
\end{align*}
In particular, $\Gamma_{j\ell}=1$ if and only if $j\in\An(\ell)$, whereas $\Gamma_{j\ell}=1/2$ if and only if $j$ and $\ell$ have no common ancestors. Finally, for $j \ne s$, the upper-tail dependence coefficient between $X_j$ and $X_s$ is given by
\[
    \chi_{js}
    =
    \nu_{\{j,s\}}([0,1]^d)
    =
    \sum_{r\in\An(j)\cap\An(s)}
    (w_{jr}\wedge w_{sr}).
\]
\end{proposition}

In view of \eqref{eq:covariance-identity}, the entire limiting covariance matrix $\Sigma$ in Theorem~\ref{theo:CTC-multivariate} can hence be expressed in terms of the weights $w$ and the distributions $\mu_r$. Note that $\sigma_{j\ell}^2=0$ if and only if $p_{j\ell}=1$, that is, if and only if $\Gamma_{j\ell}=1$.

\subsection{Heavy-tailed linear structural equation models}
\label{sec:linear-scm}

A random vector $\bm X=(X_1, \dots, X_d) \in [0,\infty)^d$ is said to satisfy a linear structural equation model with innovations $\varepsilon_1, \dots, \varepsilon_d$ and underlying DAG $\mathcal G$ if
\begin{align*}
    X_j=\sum_{\ell\in \pa(j)} b_{j\ell} X_\ell + b_{jj}\varepsilon_j, \qquad j \in [d],
\end{align*}
for coefficients $b_{j\ell} \ge 0$ such that $b_{j\ell}>0$ if and only if $\ell \in \pa(j) \cup \{j\}$. 
Suppose that the innovations
$\varepsilon_1,\dots,\varepsilon_d$ are IID non-negative
random variables with continuous distribution and
regularly varying survival function with index
$-\alpha$, for some $\alpha>0$.
The model admits the ancestral representation
\begin{align}
    X_j
    =
    \sum_{r\in\An(j)}
    \beta_{r\to j}\varepsilon_r,
    \qquad j\in[d],
    \label{eq:linear-scm-ancestral}
\end{align}
where the coefficients $\beta_{r\to j}$ are determined by the directed paths from $r$ to $j$ and satisfy $\beta_{r\to j}>0$ if $r\in\An(j)$ and $\beta_{r\to j}=0$ if $r\notin\An(j)$. Specifically, we have $\beta := (\beta_{r \to j})_{jr} = (I-B')^{-1} \diag(b_{11}, \dots, b_{dd})$, where $B'=(b_{j\ell} \bm 1(j \ne \ell))_{j\ell}$ and where $I-B'$ is invertible because $B'$ is nilpotent as a consequence of $\mathcal G$ being acyclic. We refer to \citet{gnecco2021causal} for further details on this model and its causal tail coefficients. 

\begin{lemma}[Verification of Conditions \ref{cond:conditional-weak-limit-bivariate} and \ref{cond:multivariate} for linear structural equation models]
\label{lem:linear-scm-multivariate}
Consider the linear structural equation model in \eqref{eq:linear-scm-ancestral}.
Then $\bm X$ has an ancestral tail mixture structure in the sense of Definition~\ref{def:ancestral-mixture}, with weights
\begin{align}
\label{eq:weights-lsem}
    w_{jr}
    = c_j^{-1}\beta_{r\to j}^{\alpha},
    \qquad c_j=\sum_{q\in \An(j)}\beta_{q\to j}^{\alpha}.
\end{align}
In particular, Condition~\ref{cond:multivariate} holds, and every bivariate survival copula $D_{j\ell}$, $(j,\ell)\in\Pc_2$, satisfies Condition~\ref{cond:conditional-weak-limit-bivariate}.
\end{lemma}

As a consequence of Lemma~\ref{lem:linear-scm-multivariate} and Proposition~\ref{prop:ancestral-mixture}, we have
\begin{align*}
    \Gamma_{j\ell}
    &=
    \frac12
    +
    \frac12\frac{
    \sum_{r\in\An(j)\cap\An(\ell)}
    \beta_{r\to j}^{\alpha}
    }{
    \sum_{r\in\An(j)}
    \beta_{r\to j}^{\alpha}
    },
\end{align*}
which recovers Lemma~1 of \citet{gnecco2021causal}.

We next verify Condition~\ref{cond:continuity2-bivariate}, which concerns the behavior of the first partial derivative of the survival copula near the boundary.
For this, we impose an additional regularity assumption on the innovation density.

\begin{lemma}[Verification of Condition~\ref{cond:continuity2-bivariate} for linear structural equation models]
\label{lem:linear-scm-condition-D}
Consider the linear structural equation model
\eqref{eq:linear-scm-ancestral}.
Suppose additionally that the innovations admit a density $f_\varepsilon$ that is
eventually positive and continuous and satisfies
\[
    f_\varepsilon\in\mathrm{RV}_{-(\alpha+1)}.
\]
Then, for every $(j,\ell)\in\Pc_2$, the bivariate survival copula $D_{j\ell}$ satisfies Condition~\ref{cond:continuity2-bivariate}.
\end{lemma}

Finally, we provide a sufficient condition for Condition~\ref{cond:bias-bivariate}, which requires a suitable rate of convergence of the pre-asymptotic
causal tail coefficient to its limit. Such a rate can be obtained under a Hall--Welsh-type second-order assumption on the innovation tails \citep{HallWelch1984}.

\begin{lemma}[Verification of Condition~\ref{cond:bias-bivariate} for linear structural equation models]
\label{lem:linear-scm-bias}
Consider the linear structural equation model \eqref{eq:linear-scm-ancestral}.
Suppose additionally that, for some $C,\alpha,\beta>0$ and $A<\infty$,
\[
    \bar F_\varepsilon(x)
    =
    Cx^{-\alpha}\{1+r(x)\},
    \qquad
    |r(x)|\le Ax^{-\beta},
\]
for all sufficiently large $x$. Then, for every $(j,\ell)\in\Pc_2$ and every
$\delta\in(0,\alpha)$ with $\delta\le\min\{\beta,1\}$,
\[
    B_{j\ell}(t)
    =
    t\int_0^1 D_{j\ell}(1/t,v)  \, \diff v - \Gamma_{j\ell}
    =
    O(t^{-\delta/\alpha}),
    \qquad t\to\infty.
\]
In particular, whenever $k = o(n^{2\delta/(\alpha+2\delta)})$ for some such $\delta$, Condition~\ref{cond:bias-bivariate} holds for every copula $D_{j\ell}$.
\end{lemma}

\subsection{Recursive max-linear models}
\label{sec:max-linear-scm}

A random vector $\bm X=(X_1, \dots, X_d) \in [0,\infty)^d$ is said to satisfy a recursive max-linear model with innovations $\varepsilon_1, \dots, \varepsilon_d$ and underlying DAG $\mathcal G$ if
\begin{align*}
    X_j=\bigvee_{\ell\in \pa(j)} b_{j\ell} X_\ell \vee b_{jj}\varepsilon_j, \qquad j \in [d],
\end{align*}
for coefficients $b_{j\ell} \ge 0$ such that $b_{j\ell}>0$ if and only if $\ell \in \pa(j) \cup \{j\}$. Stochastic properties related to the tail behavior of recursive max-linear models have been studied in \citet{Gissibl2018a, KluKra21}, among others.  A key structural feature is their ancestral representation
\begin{align}
X_j
=
\bigvee_{r\in\An(j)} \beta_{r\to j}\varepsilon_r,
\qquad j \in [d],
\label{eq:max-linear-ancestral}
\end{align}
for suitable non-negative coefficients $\beta_{r\to j}$ depending on the DAG $\mathcal G$ and the coefficients $b_{j\ell}$, with $\beta_{r\to j}>0$ if and only if $r \in \An(j)$; see \citet{Gissibl2018a}. Note the similarity to \eqref{eq:linear-scm-ancestral}.

Throughout, we assume that the innovations are standard $\alpha$-Fr\'echet distributed for some $\alpha>0$, that is, $\Prob(\varepsilon_r \le z) = \exp(-z^{-\alpha})$ for $z>0$.
Since $X_j^\alpha = \bigvee_{r \in \An(j)} \beta_{r\to j}^\alpha \varepsilon_r^\alpha$ and $\varepsilon_r^\alpha$ is standard Fr\'echet, the copula of $\bm X$ depends on $\alpha$ and the coefficients $\beta_{r\to j}$ only through $\beta_{r\to j}^\alpha$.
More precisely, for $\bm x \in (0,\infty)^d$,
\[
\Prob(\bm X \le \bm x)
=
\Prob\Big(\varepsilon_r \le \min_{j \in [d]: \beta_{r\to j}>0} \frac{x_j}{\beta_{r\to j}}  , \ \forall r \in [d]\Big)
=
\exp\Big( - \sum_{r \in [d]} \bigvee_{j \in [d]} \beta_{r\to j}^\alpha x_j^{-\alpha}\Big).
\]
In particular, $F_j(x) = \exp(-c_j x^{-\alpha})$ with $c_j = \sum_{r \in [d]} \beta_{r\to j}^\alpha$, and the copula of $\bm X$ is the extreme value copula $C(\bm u) = \exp\{-L(-\log u_1, \dots, -\log u_d)\}$ with stable tail dependence function
\begin{align}
L(\bm x)
=
\sum_{r\in[d]}
\bigvee_{j\in[d]}
w_{jr}\, x_j,
\label{eq:max-linear-stdf}
\end{align}
where $w_{jr} = c_j^{-1}\beta_{r\to j}^{\alpha}$, with the coefficients $\beta_{r\to j}$ from \eqref{eq:max-linear-ancestral}, has the same form as the weights in the linear SEM in \eqref{eq:weights-lsem}. Further, for each fixed $(j,\ell)$ with $j \ne \ell$, the copula $C_{j\ell}$ of the subvector $(X_j,X_\ell)$ is a bivariate extreme value copula with stable tail dependence function
\[
L_{j\ell}(x_j,x_\ell)
=
\sum_{r\in[d]} \max(w_{jr}x_j, w_{\ell r} x_\ell),
\]
and, in the parametrization of Example~\ref{ex:evc-conditions}, with Pickands dependence function
\begin{align}
\label{eq:max-linear-pickands}
A_{j\ell}(t) = L_{j\ell}(1-t,t) = \sum_{r \in [d]} \max\big\{w_{jr}(1-t), w_{\ell r} t\big\}, \qquad t \in [0,1].
\end{align}
As a consequence, the bivariate Conditions~\ref{cond:continuity2-bivariate} and \ref{cond:bias-bivariate} can be verified by the results for bivariate extreme value copulas from Section~\ref{sec:bivariate-CTC}. The only additional structural input is that $A_{j\ell}$ is piecewise linear and therefore affine in a left neighbourhood of $1$. The ancestral tail mixture structure, which concerns the joint behavior of all $d$ coordinates, follows from a direct calculation based on the explicit form of the joint CDF.

\begin{lemma}[Verification of Conditions \ref{cond:conditional-weak-limit-bivariate} and \ref{cond:multivariate} for recursive max-linear models]
\label{lem:max-linear-multivariate}
Consider the recursive max-linear model in \eqref{eq:max-linear-ancestral} with IID standard $\alpha$-Fr\'echet innovations.
Then $\bm X$ has an ancestral tail mixture structure in the sense of Definition~\ref{def:ancestral-mixture}, with weights $w_{jr}$ given by \eqref{eq:weights-lsem}, where $\beta_{r\to j}$ now denote the coefficients from \eqref{eq:max-linear-ancestral}.
In particular, Condition~\ref{cond:multivariate} holds, and every bivariate survival copula $D_{j\ell}$, $(j,\ell)\in\Pc_2$, satisfies Condition~\ref{cond:conditional-weak-limit-bivariate}.
\end{lemma}

As a consequence of Lemma~\ref{lem:max-linear-multivariate} and Proposition~\ref{prop:ancestral-mixture},
\begin{align}
    \Gamma_{j\ell}
    =
    \frac12
    +
    \frac12
    \sum_{r\in\An(j)\cap\An(\ell)}
    w_{jr},
    \label{eq:max-linear-ctc}
\end{align}
that is, $\Gamma_{j\ell}$ only depends on the normalized weights of $X_j$ associated with common ancestors of $j$ and $\ell$. Alternatively, \eqref{eq:max-linear-ctc} follows from Example~\ref{ex:evc-conditions}, since $A_{j\ell}'(1-)=p_{j\ell}$; see \eqref{eq:max-linear-pickands-affine} in the appendix.

Condition~\ref{cond:continuity2-bivariate} is satisfied by every bivariate extreme value copula, see Example~\ref{ex:evc-D}. Hence, contrary to the linear structural equation model, no additional assumption is needed. 

\begin{lemma}[Verification of Condition~\ref{cond:continuity2-bivariate} for recursive max-linear models]
\label{lem:max-linear-condition-D}
Consider the recursive max-linear model in \eqref{eq:max-linear-ancestral} with IID standard $\alpha$-Fr\'echet innovations.
Then, for every $(j,\ell)\in\Pc_2$, the bivariate survival copula $D_{j\ell}$ satisfies Condition~\ref{cond:continuity2-bivariate}.
\end{lemma}

Finally, since the Pickands dependence functions $A_{j\ell}$ are affine in a left neighbourhood of $1$, the pre-asymptotic CTCs converge to their limits at the parametric rate $1/t$, by Example~\ref{ex:evc-bias}.

\begin{lemma}[Verification of Condition~\ref{cond:bias-bivariate} for recursive max-linear models]
\label{lem:max-linear-bias}
Consider the recursive max-linear model in \eqref{eq:max-linear-ancestral} with IID standard $\alpha$-Fr\'echet innovations.
Then, for every $(j,\ell) \in \Pc_2$, the Pickands dependence function $A_{j\ell}$ is affine in a left neighbourhood of $1$, and hence $B_{j\ell}(t)=O(t^{-1})$ as $t \to \infty$. In particular, Condition~\ref{cond:bias-bivariate} holds for every copula $D_{j\ell}$ whenever $k = o(n^{2/3})$.
\end{lemma}

\section{Graph recovery}
\label{sec:graph-recovery}

In the structural equation models of Section~\ref{sec:causal-models}, the population CTC matrix carries causal information:  we have $\Gamma_{j\ell}=1$ if and only if $j$ is an ancestor of $\ell$, whereas $\Gamma_{j\ell}\in[1/2,1)$ otherwise; see Proposition~\ref{prop:ancestral-mixture}. Consequently, the  CTC matrix determines the \emph{ancestral graph} $\mathcal A=([d],E)$ of the underlying DAG $\mathcal G$, whose edge set is defined as
\[
E=\big\{(j,\ell)\in\Pc_2: j\in\an(\ell)\big\}.
\]
Note that each edge in the underlying DAG is also contained in $E$.
The ancestral graph encodes all ancestor--descendant relations and thus is a natural target for causal discovery. \citet{gnecco2021causal} use the CTC to estimate a causal ordering, i.e., an ordering on $[d]$ such that $j<\ell$ if $(j,\ell) \in E$, 
without quantifying the statistical uncertainty of the result.

From a statistical perspective, the characterization $\Gamma_{j\ell}=1$ is delicate for two reasons. First, the value $1$ lies on the boundary of the parameter space, and estimates of $\Gamma_{j\ell}$ that equal one are difficult to distinguish from values that are close to one. Second, the asymptotic variance $\sigma_{j\ell}^2$ in Theorem~\ref{theo:CTC-multivariate} vanishes exactly when $\Gamma_{j\ell}=1$, such that the normal approximation degenerates. We therefore consider, for a threshold $\gamma\in(0,1)$, the \emph{$\gamma$-signal graph} $\mathcal A_\gamma=([d],E_\gamma)$ with edge set
\[
    E_\gamma
    =
    \big\{(j,\ell)\in\Pc_2:\Gamma_{j\ell}>\gamma\big\}.
\]
In the models of Section~\ref{sec:causal-models}, we have $E\subset E_\gamma$ for every $\gamma<1$, with equality if and only if $\gamma\ge\gamma^*:=\max_{(j,\ell)\in\Pc_2\setminus E}\Gamma_{j\ell}$. Since $\gamma^*<1$, the $\gamma$-signal graph coincides with the ancestral graph for all $\gamma$ sufficiently close to one. Beyond these models, $\mathcal A_\gamma$ remains well defined and collects all pairs for which extremes of $X_j$ are accompanied by large values of $X_\ell$, measured on the scale of the CTC.
Note that $\mathcal A_\gamma$ need not be acyclic: for $\gamma<\gamma^*$, it may
contain both $(j,\ell)$ and $(\ell,j)$, namely when $X_j$ and $X_\ell$ share
ancestors carrying most of their tail weight. For the same reason, the gap
$1-\gamma^*$ can be small in large and dense DAGs, making it difficult to
separate ancestral from non-ancestral~pairs.

Below, we construct a lower and an upper confidence graph with edge sets $\widehat E^L_{\gamma,\alpha}$ and $\widehat E^U_{\gamma,\alpha}$ such that $\widehat E^L_{\gamma,\alpha}\subset E_\gamma\subset\widehat E^U_{\gamma,\alpha}$ with asymptotic probability at least $1-\alpha$. The construction is based on simultaneous bootstrap confidence intervals for the CTC matrix, as discussed after Theorem~\ref{theo:bootstrap-multivariate}. Since some of the variances $\sigma_{j\ell}^2$ may vanish, we use a studentized maximum statistic, where each estimated standard deviation is bounded from below by a fraction $c$ of the largest one.

For $\alpha\in(0,1)$ and $c>0$, let
$\widehat q_{1-\alpha}$ denote the conditional
$(1-\alpha)$-quantile of
\[
    \max_{(j,\ell)\in\Pc_2}
    \frac{|\widehat\Gb_{n,j\ell}|}
    {\widehat s_{j\ell}},
    \qquad
    \widehat s_{j\ell}
    =
    \max\{\widehat\sigma_{j\ell},\widehat\delta_n\}
    \qquad 
    \widehat \delta_n = 
c\max_{(j,\ell)\in\mathcal P_2}\widehat\sigma_{j\ell},
\]
where $\widehat\Gb_{n,j\ell}$ and
$\widehat\sigma_{j\ell}^2
=\widehat\sigma_{j\ell;j\ell}^{(n)}$
are defined in Theorem~\ref{theo:bootstrap-multivariate}. For completeness, we define $\widehat q_{1-\alpha}=0$ if $\hat \delta_n =0$; the assumptions below imply that the probability of this event converges to zero.  We define the lower and upper confidence graphs through their respective edge sets
\begin{align*}
    \widehat E_{\gamma,\alpha}^L
    &=
    \Big\{
    (j,\ell)\in\Pc_2:
    \widehat\Gamma_{j\ell}
    -\frac{\widehat q_{1-\alpha}
    \widehat s_{j\ell}}{\sqrt{k}}
    >\gamma
    \Big\},\\
    \widehat E_{\gamma,\alpha}^U
    &=
    \Big\{
    (j,\ell)\in\Pc_2:
    \widehat\Gamma_{j\ell}
    +\frac{\widehat q_{1-\alpha}
    \widehat s_{j\ell}}{\sqrt{k}}
    >\gamma
    \Big\}.
\end{align*}

\begin{proposition}[Confidence graphs and graph recovery]
\label{prop:graph-recovery-studentized}
Suppose that the assumptions of Theorem~\ref{theo:bootstrap-multivariate} hold.
Additionally, suppose that
\[
    \sqrt{k}
    \max_{(j,\ell)\in\Pc_2}
    |\Gamma_{j\ell}^{(n)}-\Gamma_{j\ell}|
    \to 0
\]
and that $\sigma_{j\ell}^2>0$ for at least one
$(j,\ell)\in\Pc_2$.
Then, for every $\gamma,\alpha\in(0,1)$ and $c>0$,
\begin{align}
\label{eq:graphrecovery-studentized}
    \liminf_{n\to\infty}
    \Prob\big(
        \widehat E_{\gamma,\alpha}^L
        \subset E_\gamma
        \subset \widehat E_{\gamma,\alpha}^U
    \big)
    \ge 1-\alpha.
\end{align}
Moreover,
\[
    \lim_{n\to\infty}
    \Prob\big(
        E_\gamma\subset\widehat E_{\gamma,\alpha}^L
    \big)=1,
\]
and consequently,
\[
    \liminf_{n\to\infty}
    \Prob\big(
        \widehat E_{\gamma,\alpha}^L=E_\gamma
    \big)
    \ge 1-\alpha.
\]
If, additionally, $\Gamma_{j\ell}\ne\gamma$
for all $(j,\ell)\in\Pc_2$, then
\[
    \lim_{n\to\infty}
    \Prob\big(
        \widehat E_{\gamma,\alpha}^L
        =E_\gamma
        =\widehat E_{\gamma,\alpha}^U
    \big)=1.
\]
\end{proposition}

\section{Beyond the CTC: conditional tail coefficients}
\label{sec:general-tail-quantities}

Condition~\ref{cond:conditional-weak-limit-bivariate} frames the CTC as a
functional of the conditional tail limit distribution: $\Gamma=\Exp[1-Z]$.
This suggests replacing $z\mapsto 1-z$ by other functions. For
$g:[0,1]\to\R$, we refer to
\[
\Gamma_g = \Exp[g(Z)] = \lim_{h \downarrow 0} \Exp[g(V_2) \mid V_1 < h]
\]
as the \emph{conditional tail coefficient} associated with $g$, where the
second equality holds for continuous $g$. In particular, $\Gamma=\Gamma_g$
for $g(z)=1-z$. Such coefficients appear naturally: for instance, by
Theorem~\ref{theo:CTC-multivariate}, the asymptotic variance of the empirical
CTC equals $\Var(Z)=\Gamma_{g}-\Gamma^2$ with $g(z)=(1-z)^2$. Below, we
discuss further examples and show that the asymptotic theory of
Section~\ref{sec:multivariate-CTC} extends to $\Gamma_g$.

\begin{example}[Generalized power CTCs]
From a causal perspective, one is often interested in situations where the CTC equals 1. This happens if and only if the conditional tail limit distribution is
degenerate at 0. A natural generalization of the CTC  that emphasizes the point mass at 0  is given by the \emph{generalized power CTCs}
\[
    \Gamma(q)=\Exp[(1-Z)^q], \qquad q\in\N.
\]
The results below show that the empirical counterpart has asymptotic variance $\Var((1-Z)^q)$, which may be smaller than $\Var(Z)$ and thus motivates efficiency
considerations. Moreover, by dominated convergence,
\[
    \lim_{q\to\infty}\Gamma(q)=\Prob(Z=0),
\]
which is of particular interest from a causal perspective. Finally, note that the sequence $(\Gamma(q))_{q\in\N}$ uniquely determines the
distribution of $Z$ \citep[Theorem~30.1]{Billingsley1995}.
\end{example}

\begin{example}[Legendre CTCs]
The classical Legendre polynomials $P_q$ form an orthogonal basis of $L^2([-1,1])$ and satisfy
\[
    \int_{-1}^1 P_q(x)P_m(x)\,\diff x
    =\frac{2}{2q+1}\bm 1(m=q), \qquad m,q \in \N.
\]
The \emph{shifted Legendre polynomials}, defined by $\phi_q(z)=P_q(1-2z)$, are
orthogonal on $[0,1]$. The first four are
\begin{align*}
    \phi_0(z)&=1,\\
    \phi_1(z)&=1-2z,\\
    \phi_2(z)&=1-6z+6z^2,\\
    \phi_3(z)&=1-12z+30z^2-20z^3.
\end{align*}
In general, they satisfy
$
    \Exp[\phi_q(U)\phi_m(U)]
    =\bm 1(m=q)/(2q+1)$ for $U\sim\mathrm{Unif}(0,1)$.
We define the $q$th Legendre CTC by
\[
    \theta_q=\Exp[\phi_q(Z)],
    \qquad q\in\N.
\]
Note that the classical CTC satisfies $\Gamma=(1+\theta_1)/2$.
Moreover, since $\phi_q(0)=1$ and
$\Exp[\phi_q(U)]=0$ for $q\ge1$, we have $\theta_q=p$ for all $q\in\N$,
whenever
$Z\sim p\delta_0+(1-p)\mathrm{Unif}(0,1)$,
as in the linear SCM from Section~\ref{sec:linear-scm}.
Consequently, the identities $\theta_1=\theta_2=\cdots$ provide moment restrictions that may be useful for assessing the conditional tail limit distribution implied by the linear SCM.
\end{example}

A detailed discussion of these coefficients is beyond the scope of this paper. However, the following results provide a basis for statistical inference on general conditional tail coefficients. Let $g:[0,1] \to \R$ be a continuous function, and define
\begin{align}\label{eq:est_ghat}
\gghat (x) 
&= \frac1k \sum_{i=1}^n g(\widehat V_{i2}) \bm 1\Big(\widehat V_{i1} \le \frac{kx}n \Big),\\
\ggmean (x) 
&=  \label{eq:gn-g}
\frac{n}k \Exp \Big[ g(V_2) \bm 1\Big(V_1 \le \frac{kx}n\Big) \Big].
\end{align}
For $g(v)=1-v$, these quantities reduce to
$\gghat(1)=\widehat\Gamma_n$ and
$\ggmean(1)=\Gamma_n$.
Further, for any $x>0$,  we may write
\begin{align} \label{eq:gn-Ag}
\ggmean(x) =
x \ggcmean \Big(\frac{kx}n \Big), \qquad \ggcmean(v) = \Exp\Big[ g(V_2) \mid V_1 \le v \Big].
\end{align}
Under Condition~\ref{cond:conditional-weak-limit-bivariate},
continuity of $g$ implies that
\[
\ggcmean(0) := \lim_{v \to 0} \ggcmean(v) = \Exp_{Z \sim H}[g(Z)]
\]
exists, and the equation $\ggmean(x) = x \ggcmean(kx/n)$ then also holds for $x=0$.

\begin{theorem}[Linearization of empirical conditional tail coefficients]
\label{theo:unknown-margins-bivariate}
Let $D$ be a bivariate copula, and suppose that $(V_1, V_2) \sim D$ satisfies \ref{cond:continuity2-bivariate}. 
Let $(X_{i1}, X_{i2})_{i=1}^n$ be $n$ IID observations with continuous marginal CDFs and survival copula $D$, let $g:[0,1] \to \R$ be Lipschitz-continuous,  let $k \in [n]$ with $k \to \infty$ and $k=o(n)$ as $n \to \infty$, 
and let
\begin{align*}
\ggprocess(x) := 
\sqrt{k} \big\{ \gghat(x) - \ggmean(x) \big\}.
\end{align*}
Then, for any $T >0$ and as $n \to \infty$,
\begin{align} \label{eq:linear-bivariate-process-g}
    \sup_{x \in [0,T]} \Big| \ggprocess(x) 
    -
    \frac1{\sqrt n} \sum_{i=1}^n  \phi_{n,x}(V_{i1},V_{i2}) \Big|  = o_\Prob(1)
    \end{align}
    for $n \to \infty$, where
    \[
    \phi_{n,x}(v_1,v_2)
    =  \sqrt{\frac{n}{k}} \Big\{ g(v_2) - \ggcmean\Big(\frac{kx}n \Big) \Big\} \bm 1\Big(v_{1} \le  \frac{kx}{n} \Big).
\]
\end{theorem}

\begin{corollary}[Functional CLT for empirical conditional tail coefficients]
\label{cor:unknown-margins-bivariate}
Under the conditions of Theorem~\ref{theo:unknown-margins-bivariate}, we have
\begin{align*}    \big( \ggprocess(x)\big)_{x \in [0,T]} \weak 
    \big( \gglimitprocess(x)\big)_{x \in [0,T]}
    \qquad \text{ in } \quad \ell^{\infty}([0,T]),
\end{align*}
where $\gglimitprocess$ is a tight, centered Gaussian process on $[0,\infty)$ with almost sure continuous sample paths and with covariance function 
\[
\Cov\big[\gglimitprocess(x),\gglimitprocess(y)\big]
=
(x\wedge y) \Var_{Z \sim H}[g(Z)].
\]
\end{corollary}

\begin{theorem}[Multiplier bootstrap for empirical conditional tail coefficients]
\label{theo:bootstrap}
    Consider the setting of Theorem~\ref{theo:unknown-margins-bivariate}. Let $\xi_1, \xi_2, \dots$ be IID centered random variables with variance $1$, independent of the $\bm X_i$. Let
    \[
    \widehat \Hb_n(x) := \frac1{\sqrt n} \sum_{i=1}^n \xi_i \, \widehat \phi_{n,x}(\widehat V_{i1}, \widehat V_{i2}),
    \]
    where, using the convention $\gghat(0)/0=0$, 
    \[
    \widehat \phi_{n,x}(v_1, v_2) = \sqrt{\frac{n}{k}} \Big\{ g(v_2) - \frac{\gghat(x)}{x} \Big\} \bm 1\Big(v_{1} \le  \frac{kx}{n} \Big).
    \]
    Then, for every $T >0$, 
    \[
    \big( \widehat \Hb_n(x)\big)_{x \in [0,T]} \cweak{\xi}
    \big( \gglimitprocess(x)\big)_{x \in [0,T]}
    \qquad \text{ in } \quad \ell^{\infty}([0,T]).
    \]
\end{theorem}

\section{Simulation study} 
\label{sec:simulations}

In this section, we investigate the finite-sample performance of the CTC estimators through simulations in bivariate and multivariate settings. We first consider bivariate max-linear copulas to assess the performance of $\hat\Gamma_n$ and compare the plug-in and multiplier-bootstrap estimators of its asymptotic variance. We then turn to the multivariate linear structural equation framework of Section~\ref{sec:linear-scm}, focusing on the joint coverage probabilities of the estimated lower and upper confidence graphs. Throughout, the bootstrap multipliers $\xi_i$ are IID standard normal random variables.

\subsection{Finite sample properties of the CTC estimators}

We begin by investigating the normal approximation for the CTC estimators and the variance reduction achieved by the rank-based estimator $\widehat\Gamma_n$ relative to the estimator $\widehat\Gamma_n^*$ in~\eqref{eq:knownCTC}, which assumes known margins and is referred to as the \emph{oracle estimator}.

We consider a bivariate RMLM defined recursively by $X_1=b_{12}X_2\vee b_{11}\varepsilon_1$ and $X_2=b_{22}\varepsilon_2$, where
$\varepsilon_1,\varepsilon_2$ are independent standard Fr\'echet innovations and
$b_{11}=\theta$, $b_{12}=1-\theta$, and $b_{22}=1$ for some $\theta\in(0,1)$.
The underlying DAG consists of the single edge $2\to1$, and the ancestral
representation \eqref{eq:max-linear-ancestral} reads
$X_1=\theta\varepsilon_1\vee(1-\theta)\varepsilon_2$ and $X_2=\varepsilon_2$,
that is, $\beta_{1\to1}=\theta$, $\beta_{2\to1}=1-\theta$ and $\beta_{2\to2}=1$.
In particular, both margins are standard Fr\'echet.
By Proposition~\ref{prop:ancestral-mixture}, $\Gamma_{21}=1$,
$\Gamma_{12}=1-\theta/2$ and $\chi_{12}=1-\theta$, so that $\theta$
controls the strength of extremal dependence between $X_1$ and $X_2$, with
smaller values corresponding to stronger dependence.

\begin{figure}[thb!]
    \hspace*{2.5cm}%
    \includegraphics[width=.875\linewidth, height=.425\linewidth]{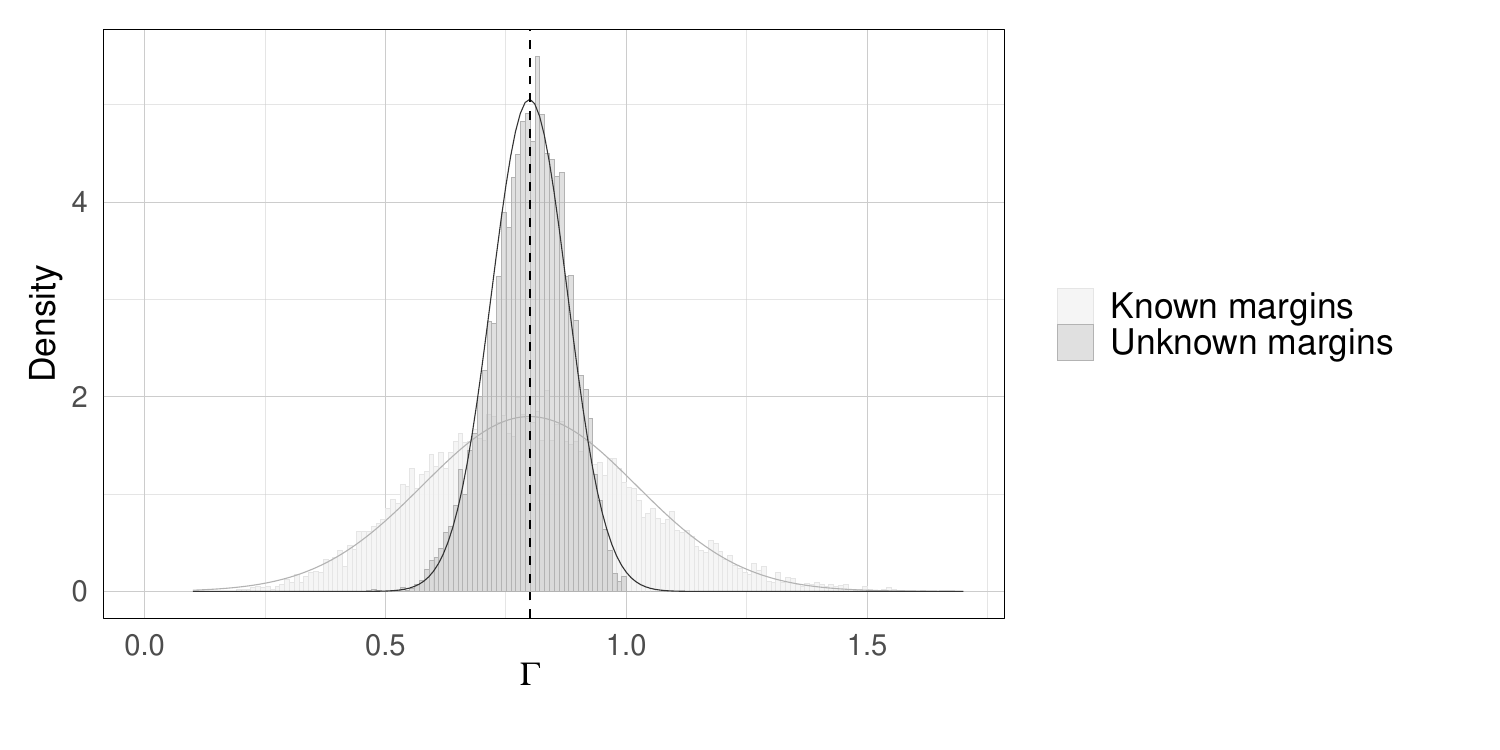} 
    \caption{Histograms of the rank-based and oracle CTC estimators,
    together with fitted normal densities, based on
    $R=10{,}000$ replications with $n=5000$ and
    $k=\lfloor 0.5 \, n^{0.4}\rfloor$.
    The population CTC is $\Gamma_{12}=0.8$.}
   \label{fig:densityplots}
\end{figure}

As a first illustration, Figure~\ref{fig:densityplots} displays histograms of $\widehat\Gamma_n$ and $\widehat\Gamma_n^*$, together with fitted normal densities,
based on $R=10{,}000$ independent replications with $\theta=0.4$, $n=5000$, and $k=\lfloor 0.5 \, n^{0.4}\rfloor=15$. The fitted densities closely match the empirical distributions, supporting the normal approximation even for the relatively small number of exceedances considered here. Moreover, the rank-based estimator exhibits substantially less variability than the oracle estimator, in line with the theoretical variance reduction.

We next examine the finite-sample behaviour of both estimators
over a wider range of dependence strengths, sample sizes,
and thresholds. Specifically, we consider $\theta\in\{0.25, 0.4, 0.6, 0.85\}$, $n\in \{500, 5000, 50000\}$ and choose the number of exceedances as $k=\max( \lfloor c n^{0.4} \rfloor, 10)$ with $c \in\{0.25, 0.5, 0.75, 1\}$. For each parameter configuration, we generate $R=1000$ independent replications. The exponent $0.4$ follows the choice $k=\lfloor n^{0.4}\rfloor$ recommended by \citet{gnecco2021causal}, and is further supported by an additional simulation study in Appendix~\ref{sec:additional-simulations_doublegrid}, which shows that it provides a reasonable trade-off between the mean squared errors (MSEs) of $\widehat \Gamma_{12}$ and $\widehat \Gamma_{21}$.

Using the results of Section~\ref{sec:max-linear-scm} and the asymptotic variance formulas derived above, we obtain
\begin{align*}
    \Gamma_{12} &= 1-\frac{\theta}{2},
    &\sigma_{12}^2 &= \frac{\theta(4-3\theta)}{12},
    &\sigma_{12,\mathrm{known}}^2 &= 1-\frac{2\theta}{3},\\
    \Gamma_{21} &= 1,
    &\sigma_{21}^2 &= 0,
    &\sigma_{21,\mathrm{known}}^2 &= 1.
\end{align*}
Thus, the rank-based estimator has a strictly smaller asymptotic variance in both directions,~with~a degenerate limiting distribution in the direction $(2,1)$. In what follows, we omit the directional~subscripts from the asymptotic variances whenever the corresponding CTC is clear from the context.

\begin{figure}[thb!]
    \centering
    \includegraphics[width=.85\linewidth]{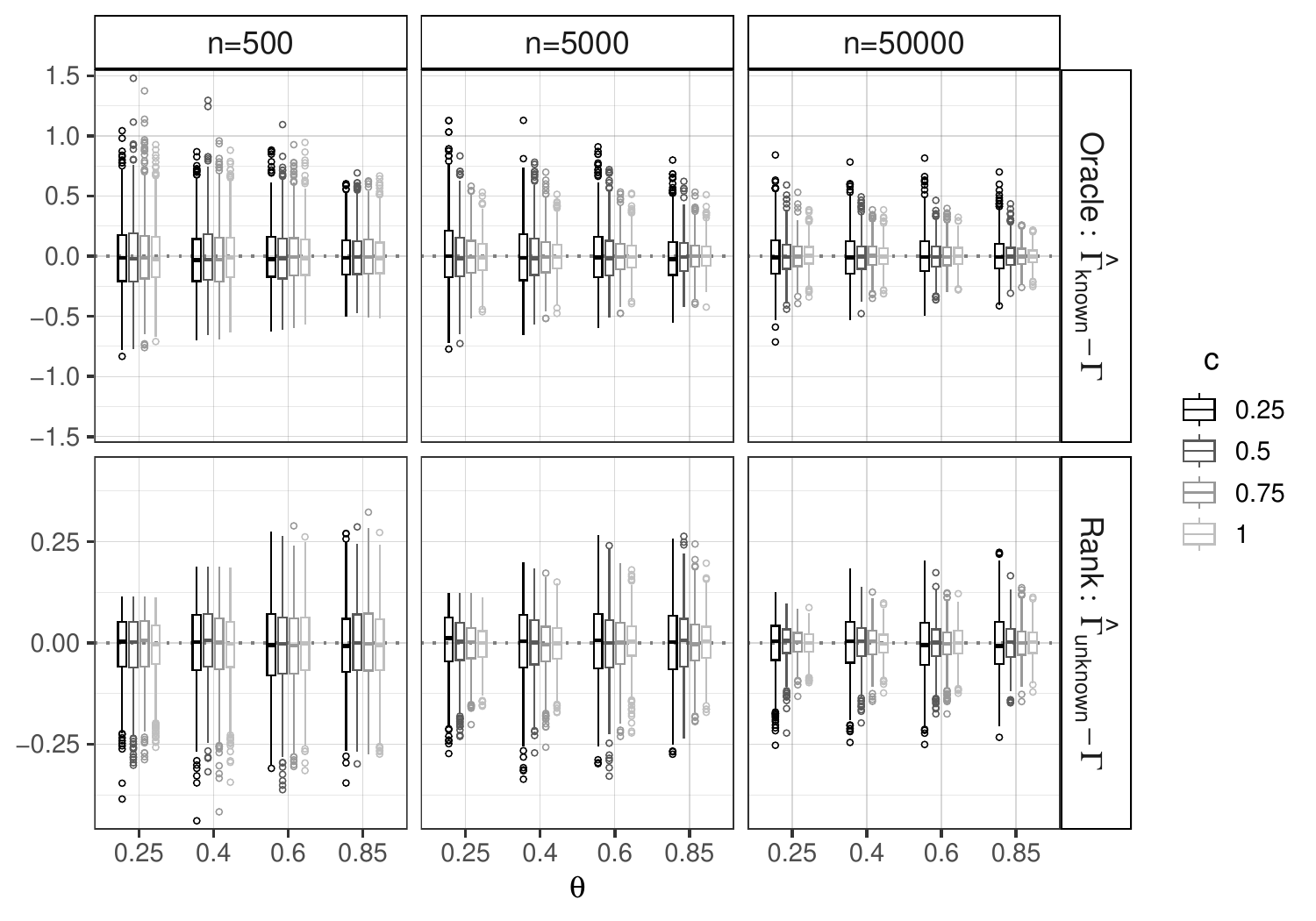} 
    \caption{Bivariate RMLM: Boxplots of the estimation
    errors $\widehat\Gamma_{12}^*-\Gamma_{12}$ (top) and $\widehat\Gamma_{12}-\Gamma_{12}$ (bottom) based on
    $R=1000$ replications for different values of
    $\theta$, $n$, and $k$.}
   \label{fig:MLgamma12}
\end{figure}

Figure~\ref{fig:MLgamma12} presents boxplots of the estimation errors for $\Gamma_{12}$ across the considered parameter configurations. The rank-based estimator consistently exhibits substantially less variability than the oracle estimator, confirming the theoretical variance reduction across different dependence strengths. For both estimators, the variability generally decreases
with increasing sample size and number of exceedances. However, larger choices of $k$ also lead to a noticeable negative bias for both estimators in some settings, 
illustrating the usual bias-variance trade-off in threshold selection. Corresponding results for $\Gamma_{21}$, where the asymptotic variance of the rank-based estimator vanishes, are provided in Figure~\ref{fig:MLgamma21} in the Appendix~\ref{sec:additional-simulations}.

\subsection{Finite-sample properties of the variance estimators}
\label{sec:variance-simulations}

We next investigate the finite-sample performance of the estimators of the asymptotic variances $\sigma_{12}^2$ and $\sigma_{21}^2$. Recall from
Corollary~\ref{cor:unknown-margins-bivariate} and Theorem~\ref{theo:bootstrap} that these variances coincide with the limiting variances of the respective multiplier bootstrap processes. We consider two approaches to variance estimation under unknown margins. The first is a plug-in estimator $\widehat\sigma^2$,
obtained by subtracting $\widehat\Gamma_n^2$ from the estimator in~\eqref{eq:est_ghat} with $g(v)=(1-v)^2$. The second is the bootstrap estimator
$\widehat\sigma_{\mathrm{boot}}^2$, defined as the empirical variance of $1000$ multiplier bootstrap replicates. We use the same simulation settings as in the preceding subsection.

\begin{figure}[thb!]
    \centering
    \includegraphics[width=.85\linewidth]
        {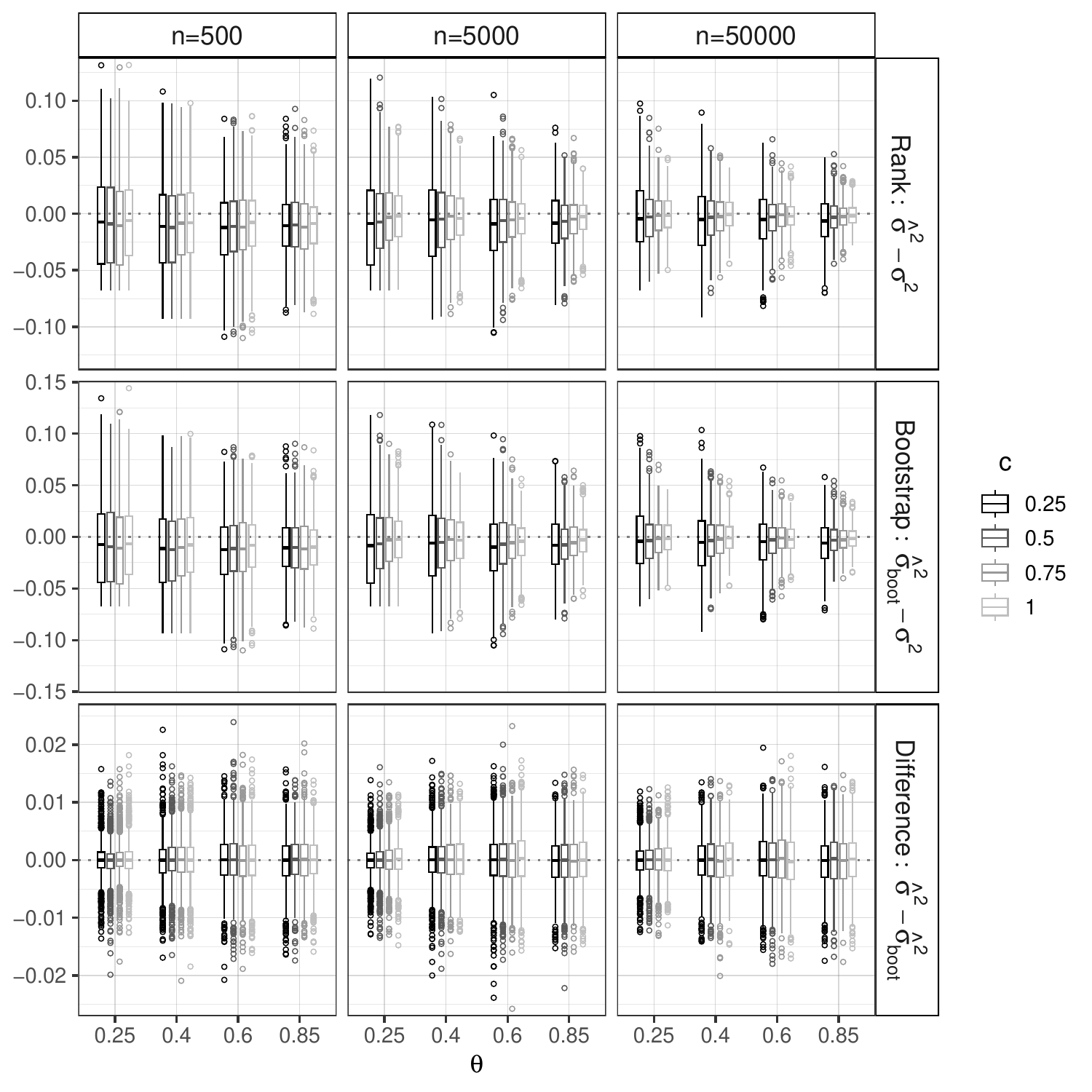}
    \caption{Bivariate RMLM, direction $(1,2)$:
    Boxplots of the estimation errors of the plug-in
    variance estimator (top) and the multiplier-bootstrap
    variance estimator (middle), together with their
    differences (bottom).
    Results are based on $R=1000$ simulation replications
    with $1000$ multiplier replicates each.}
    \label{fig:MLvars}
\end{figure}

Figure~\ref{fig:MLvars} presents boxplots of $\widehat\sigma^2-\sigma^2$ (top),
$\widehat\sigma_{\mathrm{boot}}^2-\sigma^2$ (middle), and $\widehat\sigma^2-\widehat\sigma_{\mathrm{boot}}^2$ (bottom) for the direction $(1,2)$. Both variance estimators perform well overall, with estimation errors generally decreasing as the sample size and number of exceedances increase.
Some negative bias is visible for small samples and large values of $k$, particularly under weaker dependence (large $\theta$). The differences between the two variance estimators are generally small relative to their estimation errors,
indicating close agreement between the plug-in and bootstrap approaches. Corresponding results for the direction $(2,1)$, where the asymptotic variance vanishes, are reported in Appendix~\ref{app:variance21}.

\subsection{Multivariate inference and graph recovery}
\label{subsec:sims-graph}

We finally investigate the finite-sample performance of the lower and upper
confidence graphs from Section~\ref{sec:graph-recovery} in linear SEMs as in
Section~\ref{sec:linear-scm}. Throughout, we choose $\gamma=0.95$ and
significance level $\alpha=0.1$. The target is thus the $\gamma$-signal graph
$\mathcal A_\gamma$, which coincides with the ancestral graph $\mathcal A$
whenever $\gamma\ge\gamma^*$; see Section~\ref{sec:graph-recovery}.

We consider dimensions $d\in\{5,10,20\}$, edge probabilities
$p\in\{0.1,0.2\}$ and sample sizes $n\in\{500,5000,50000\}$.
For each configuration, we generate $R=1000$ random DAGs, each with a
corresponding data set, as follows. For $j<\ell$, the edge $\ell\to j$ is
included independently with probability $p$, and its coefficient $b_{j\ell}$
is drawn from $\mathcal U[0.2,1.2]$; all other off-diagonal coefficients are
zero, so that $B'$ is strictly upper triangular. The diagonal coefficients
$b_{11},\dots,b_{dd}$ are drawn independently from $\mathcal U[0.5,1.5]$.
The innovations $\varepsilon_1,\dots,\varepsilon_d$ are IID copies of $|T|$,
where $T$ follows a Student-$t$ distribution with one degree of freedom; they
are thus regularly varying with index $1$. The data are then generated
from the ancestral representation $\bm X=\beta\bm\varepsilon$, with
$\beta=(\beta_{r \to j})_{jr}=(I-B')^{-1}\operatorname{diag}(b_{11},\dots,b_{dd})$ as in
Section~\ref{sec:linear-scm}.

For each replication $r$, we compute the CTC matrix $\widehat\Gamma=(\hat \Gamma_{r,j\ell})$ using $k=\max(\lfloor 0.5\,n^{0.4}\rfloor, 10)$ exceedances. As in Section~\ref{sec:simulations}, the exponent $0.4$ follows the recommendation of \citet{gnecco2021causal}; the choice of both $c$ and the exponent is further supported by a preliminary simulation study on $k=\max(\lfloor cn^\kappa\rfloor, 10)$ over a two-dimensional $(c, \kappa)$-grid in Appendix~\ref{app:k_graph}. 
We then construct the lower and upper confidence edge sets $\widehat E_{r,\gamma,\alpha}^L$ and $\widehat E_{r,\gamma,\alpha}^U$ as described in Section~\ref{sec:graph-recovery}, using $1000$ multiplier bootstrap replications to estimate the required quantiles. 

We evaluate the confidence graphs in terms of both simultaneous and pairwise coverage. Let $E_{r,\gamma} =\{(j,\ell) \in \mathcal P_2: \Gamma_{r,j\ell}>\gamma\}$ denote the population edge set in replication $r$, and let $e_{r,\gamma}(j, \ell)$ and $\widehat e_{r,\gamma,\alpha}^T(j, \ell)$, $T\in\{L,U\}$,
denote the corresponding edge indicators. The empirical joint and pairwise coverage probabilities are defined as
\begin{align*}
\widehat{\mathrm{Cov}}_{\mathrm{joint}}
&=
\frac1R\sum_{r=1}^R
\bm 1\big\{
\widehat E_{r,\gamma,\alpha}^L
\subseteq E_{r,\gamma}
\subseteq \widehat E_{r,\gamma,\alpha}^U
\big\},\\
\widehat{\mathrm{Cov}}_{\mathrm{pair}}
&=
\frac1{R|\mathcal P_2|}
\sum_{r=1}^R\sum_{(j,\ell)\in\mathcal P_2}
\bm 1\big\{
\widehat e_{r,\gamma,\alpha}^L(j, \ell)
\le e_{r,\gamma}(j, \ell)
\le \widehat e_{r,\gamma,\alpha}^U(j, \ell)
\big\}.
\end{align*}

To assess the recovery performance of the individual confidence graphs, we compare, for each replication $r$ and $T\in\{L,U\}$, the estimated edge set
$\widehat E_{r,\gamma,\alpha}^T$ with the population edge set $E_{r,\gamma}$.
We denote by $\mathrm{TP}_r^T$, $\mathrm{FP}_r^T$, $\mathrm{TN}_r^T$, and $\mathrm{FN}_r^T$ the numbers of ordered pairs $(j, \ell)\in\mathcal P_2$
belonging to both edge sets (true positives), only to the estimated edge set (false positives), to neither edge set (true negatives), and only to the population
edge set (false negatives), respectively. We then define
\begin{align*}
\mathrm{Prec}_r^T
&=\frac{\mathrm{TP}_r^T}
        {\mathrm{TP}_r^T+\mathrm{FP}_r^T},
&
\mathrm{TPR}_r^T
&=\frac{\mathrm{TP}_r^T}
        {\mathrm{TP}_r^T+\mathrm{FN}_r^T},\\
\mathrm{FPR}_r^T
&=\frac{\mathrm{FP}_r^T}
        {\mathrm{FP}_r^T+\mathrm{TN}_r^T},
&
\mathrm{FNR}_r^T
&=\frac{\mathrm{FN}_r^T}
        {\mathrm{FN}_r^T+\mathrm{TP}_r^T}.
\end{align*}
The reported metrics are obtained by averaging over the $R$ replications, excluding any replication for which the corresponding denominator vanishes.
More precisely, letting $\mathcal R_M^T$ denote the set of replications for which $M_r^T$ is defined, we report
\[
    \widehat M_T
    =\frac{1}{|\mathcal R_M^T|}
    \sum_{r\in\mathcal R_M^T} M_r^T,
    \qquad
    M\in\{\mathrm{Prec},\mathrm{TPR},
             \mathrm{FPR},\mathrm{FNR}\}.
\]

Finally, we assess how closely the population graph $E_{r,\gamma}$ approximates the ancestral graph. Let $\mathcal A_r$ denote the ancestral graph in
replication $r$, with edge set $E_r=\{(j,\ell)\in \Pc_2: \Gamma_{r,j\ell}=1\}$.
Since $E_r\subseteq E_{r,\gamma}$, the two graphs differ only through additional edges in the $\gamma$-signal graph. We therefore define the normalised structural
Hamming distance by
\[
    \mathrm{nSHD}_r
    =
    \frac{|E_{r,\gamma}\setminus E_r|}
    {|E_r|+|E_{r,\gamma}|},
\]
with the ratio set to zero if both edge sets are empty. We report 
$\widehat{\mathrm{nSHD}}=R^{-1}\sum_{r=1}^R
\mathrm{nSHD}_r$.
Note that this quantity depends only on the population graphs and thus measures the approximation error associated with the choice of $\gamma$, independently
of statistical estimation.

Table~\ref{tab:graph-recovery} summarises the results, with all metrics reported as percentages. Pairwise coverage is consistently high, exceeding $96\%$ across all configurations and approaching $100\%$ for larger sample sizes. Joint coverage, by contrast, depends strongly on the dimension, sparsity, and sample size. For $d=5$, it is close to or above the nominal level of $90\%$ even for $n=500$, and approaches $100\%$ as $n$ increases. For larger and denser graphs, however, simultaneous coverage can be substantially below the nominal level.
The most challenging configuration, $d=20$ and $p=0.2$, yields joint coverage probabilities of only $1\%$, $30.9\%$, and $83.1\%$ for the three sample sizes.
Thus, considerably larger samples are needed to achieve reliable simultaneous inference in more complex graphs.

The lower and upper confidence graphs exhibit the expected complementary behaviour.
The lower graphs have high precision and very low false positive rates, but may miss a substantial proportion of true edges in small samples. Conversely, the upper graphs achieve high true positive rates~at the expense of including additional edges. Both bounds become more accurate as the sample size
increases, although the upper graphs offer a relatively lower degree of precision in the higher-dimensional settings.

The normalised structural Hamming distances between the ancestral and $\gamma$-signal graphs are very small, never exceeding $1.3\%$ in the considered configurations. This indicates that the choice $\gamma=0.95$ provides a close approximation to the underlying ancestral graph while allowing for the confidence
graph construction.

Finally, additional boxplots in Appendix~\ref{app:graph} illustrate the variability
of the recovery metrics~across the replications, which is particularly pronounced for small sample sizes and larger~dimensions.

\begin{table}[t]
\centering
\caption{Graph recovery under multivariate linear SEMs:
Empirical coverage probabilities, recovery metrics for the lower and upper confidence graphs, and normalised structural Hamming distance between the ancestral and $\gamma$-signal graphs. All entries are percentages based on $R=1000$
simulation replications.}
\label{tab:graph-recovery}
\footnotesize
\setlength{\tabcolsep}{4pt}
\begin{tabular}{@{}cc*{11}{r}@{}}
\toprule
& & \multicolumn{2}{c}{Coverage}
    & \multicolumn{2}{c}{Precision}
    & \multicolumn{2}{c}{TPR}
    & \multicolumn{2}{c}{FPR}
    & \multicolumn{2}{c}{FNR}
    & \multicolumn{1}{c}{nSHD} \\
\cmidrule(lr){3-4}
\cmidrule(lr){5-6}
\cmidrule(lr){7-8}
\cmidrule(lr){9-10}
\cmidrule(lr){11-12}
\cmidrule(l){13-13}
$d$ & $n$
& ${\mathrm{Cov}}_{\mathrm{joint}}$
& ${\mathrm{Cov}}_{\mathrm{pair}}$
& $\mathrm{Prec}_L$
& $\mathrm{Prec}_U$
& $\mathrm{TPR}_L$
& $\mathrm{TPR}_U$
& $\mathrm{FPR}_L$
& $\mathrm{FPR}_U$
& $\mathrm{FNR}_L$
& $\mathrm{FNR}_U$
& $\mathrm{nSHD}$ \\
\midrule
\multicolumn{13}{@{}l}{{Panel A: $p=0.1$}} \\
\midrule
\multirow{3}{*}{$5$}
  & $500$ & 96.4 & 99.8 & 98.4 & 34.0 & 85.3 & 99.5 & 0.2 & 10.8 & 14.7 & 0.5 & 0.0 \\
 & $5000$ & 99.8 & 100.0 & 99.9 & 77.7 & 99.9 & 100.0 & 0.0 & 2.7 & 0.1 & 0.0 & 0.0 \\
 & $50000$ & 100.0 & 100.0 & 100.0 & 95.1 & 100.0 & 100.0 & 0.0 & 0.6 & 0.0 & 0.0 & 0.0 \\
\addlinespace
\multirow{3}{*}{$10$}
 & $500$ & 80.9 & 99.7 & 97.4 & 27.2 & 76.2 & 98.6 & 0.2 & 16.8 & 23.8 & 1.4 & 0.0 \\
 & $5000$ & 98.5 & 100.0 & 99.9 & 65.3 & 99.1 & 99.9 & 0.0 & 4.6 & 0.9 & 0.1 & 0.0 \\
 & $50000$ & 100.0 & 100.0 & 100.0 & 86.9 & 100.0 & 100.0 & 0.0 & 1.5 & 0.0 & 0.0 & 0.0 \\
\addlinespace
\multirow{3}{*}{$20$}
 & $500$ & 20.6 & 99.2 & 94.6 & 26.5 & 55.9 & 95.3 & 0.4 & 24.5 & 44.1 & 4.7 & 0.1 \\
 & $5000$ & 79.7 & 99.9 & 99.6 & 51.6 & 94.7 & 99.6 & 0.0 & 9.6 & 5.3 & 0.4 & 0.1 \\
 & $50000$ & 97.9 & 100.0 & 100.0 & 75.6 & 99.7 & 100.0 & 0.0 & 3.6 & 0.3 & 0.0 & 0.0 \\

\midrule
\multicolumn{13}{@{}l}{{Panel B: $p=0.2$}} \\
\midrule
\multirow{3}{*}{$5$}
 & $500$ & 91.5 & 99.5 & 98.3 & 44.4 & 83.2 & 98.8 & 0.4 & 16.7 & 16.8 & 1.2 & 0.0 \\
 & $5000$ & 99.5 & 100.0 & 99.9 & 76.9 & 99.7 & 100.0 & 0.0 & 6.4 & 0.3 & 0.0 & 0.0 \\
 & $50000$ & 100.0 & 100.0 & 100.0 & 93.9 & 100.0 & 100.0 & 0.0 & 1.7 & 0.0 & 0.0 & 0.0 \\
\addlinespace
\multirow{3}{*}{$10$}
 & $500$ & 51.1 & 98.7 & 95.2 & 37.2 & 66.2 & 96.9 & 0.9 & 29.0 & 33.8 & 3.1 & 0.1 \\
 & $5000$ & 91.3 & 99.8 & 99.6 & 60.0 & 97.4 & 99.7 & 0.1 & 13.6 & 2.6 & 0.3 & 0.1 \\
 & $50000$ & 99.3 & 100.0 & 100.0 & 81.8 & 99.9 & 100.0 & 0.0 & 5.4 & 0.1 & 0.0 & 0.0 \\
\addlinespace
\multirow{3}{*}{$20$}
 & $500$ & 1.0 & 96.5 & 86.3 & 35.5 & 49.8 & 94.2 & 3.0 & 49.1 & 50.2 & 5.8 & 1.3 \\
 & $5000$ & 30.9 & 99.3 & 97.9 & 47.4 & 90.1 & 99.2 & 0.7 & 33.3 & 9.9 & 0.8 & 1.3 \\
 & $50000$ & 83.1 & 99.9 & 99.8 & 63.6 & 97.7 & 99.9 & 0.1 & 18.3 & 2.3 & 0.1 & 1.2 \\
\bottomrule
\end{tabular} 
\end{table}

\section{Data application}
\label{sec:data-application}

We illustrate the proposed graph inference procedure using dietary intake data from the 2015--2016 National Health and Nutrition Examination Survey (NHANES)\footnote{Available at 
\url{https://wwwn.cdc.gov/Nchs/Data/Nhanes/Public/2015/DataFiles/DR1TOT_I.xpt}}.
We consider the six nutrients selected by \citet{KraliDavisonKluppelberg2026} as the largest subset whose extremal dependence structure could be reasonably
approximated by an RMLM: retinol (RET, DR1TRET), $\alpha$-carotene (AC, DR1TACAR),
$\beta$-carotene (BC, DR1TBCAR), vitamin A (VA, DR1TVARA), lutein and zeaxanthin
(LZ, DR1TLZ), and vitamin K (VK, DR1TVK). The data consist of estimated nutrient intakes from a 24-hour dietary recall. After retaining only complete observations, the effective sample size is $n=8327$. For the purposes of our analysis,
we treat the observations as independent and identically
distributed.

Our working assumption is that the joint distribution of the nutrient intakes has an ancestral tail mixture structure in the sense of Definition~\ref{def:ancestral-mixture}, which holds, for instance, if the intakes follow a heavy-tailed linear SEM or an RMLM after possibly applying strictly increasing marginal transformations. Under this assumption, $\Gamma_{j\ell}=1$ if and only if $j$ is an ancestor of $\ell$, and the proposed method allows us to quantify the uncertainty concerning the underlying ancestral graph.
The right-hand panel of Figure~\ref{fig:food6} shows the
ancestral graph estimated by
\citet{KraliDavisonKluppelberg2026} under an RMLM.
While this graph reveals several interesting directed
relationships between nutrient intakes, it comes without an assessment of the statistical
uncertainty of its edges.

\paragraph{Choice of $k$.} For $k \in [n]$, let $\hat{\mathcal A}_{\gamma,k}=([d],\widehat E_{\gamma,k})$ denote the plug-in estimate of the $\gamma$-signal graph from Section~\ref{sec:graph-recovery}, with edge set $\widehat E_{\gamma,k}=\{(j,\ell)\in\Pc_2:\widehat\Gamma_{j\ell}>\gamma\}$, where $\widehat\Gamma_{j\ell}$ is computed with $k$ exceedances. As in the simulation study, we use $\gamma=0.95$. We choose $k$ in a region where $\hat{\mathcal A}_{\gamma,k}$ is stable as a function of $k$; a conceptually related stability-based choice of $k$ was proposed by \citet{tran2024estimating} for the estimation of root-directed spanning trees. To measure the distance between two graphs, we identify a directed graph on $[d]$ with its adjacency matrix $A\in\{0,1\}^{d\times d}$, where $A_{j\ell}=1$ if and only if $(j,\ell)$ is an edge. For two adjacency matrices $A$ and $B$, let
\[
D(A,B) = \sum_{j, \ell \in [d]}\bm 1(A_{j\ell}\neq B_{j\ell})
\]
denote the number of entries in which they differ, and let $R(A,B)$ denote the number of strict reversals, that is, the number of unordered pairs $\{j,\ell\}$ with $(A_{j\ell},A_{\ell j})=(1,0)$ and $(B_{j\ell},B_{\ell j})=(0,1)$, or vice versa. The modified normalized structural Hamming distance
\[
\operatorname{nSHD}_{\mathrm{mod}}(A,B)
=\frac{D(A,B)-R(A,B)}
     {\sum_{j,\ell \in [d]}A_{j\ell}+\sum_{j,\ell \in [d]}B_{j\ell}},
\]
with $\operatorname{nSHD}_{\mathrm{mod}}(A,B)=0$ if neither graph contains an edge, then counts a strict reversal as one edit rather than two.

We consider the grid $\mathcal{K}=\{k_1,\dots,k_{21}\}=\{10,20,\ldots,210\}$, which ranges from approximately $0.1\%$ to $2.5\%$ of the sample size, and compute the distances between consecutive estimated graphs, $S_i = \operatorname{nSHD}_{\mathrm{mod}}(\hat{\mathcal A}_{\gamma,k_i}, \hat{\mathcal A}_{\gamma,k_{i+1}})$ for $i\in[20]$. Their rolling means and rolling standard deviations over five consecutive pairs,
\[
\bar S_i = \frac{1}{5}\sum_{m=0}^{4} S_{i+m},
\qquad
\hat\tau_i = \operatorname{sd}(S_i,S_{i+1},\ldots,S_{i+4}),
\qquad i\in[16],
\]
measure the average size and the variability of the changes of the estimated graph in the respective window. We select the window with the smallest rolling mean $\bar S_i$; if the minimum is attained at several consecutive indices, we take the longest such run of windows. Finally, we set $k^*=k_{i^*}$, where $i^*$ is the central index of the selected run (the smaller one if the run has even length), while the rolling standard deviations serve as a descriptive complement.

\begin{figure}
    \centering
    \includegraphics[width=0.87\linewidth]{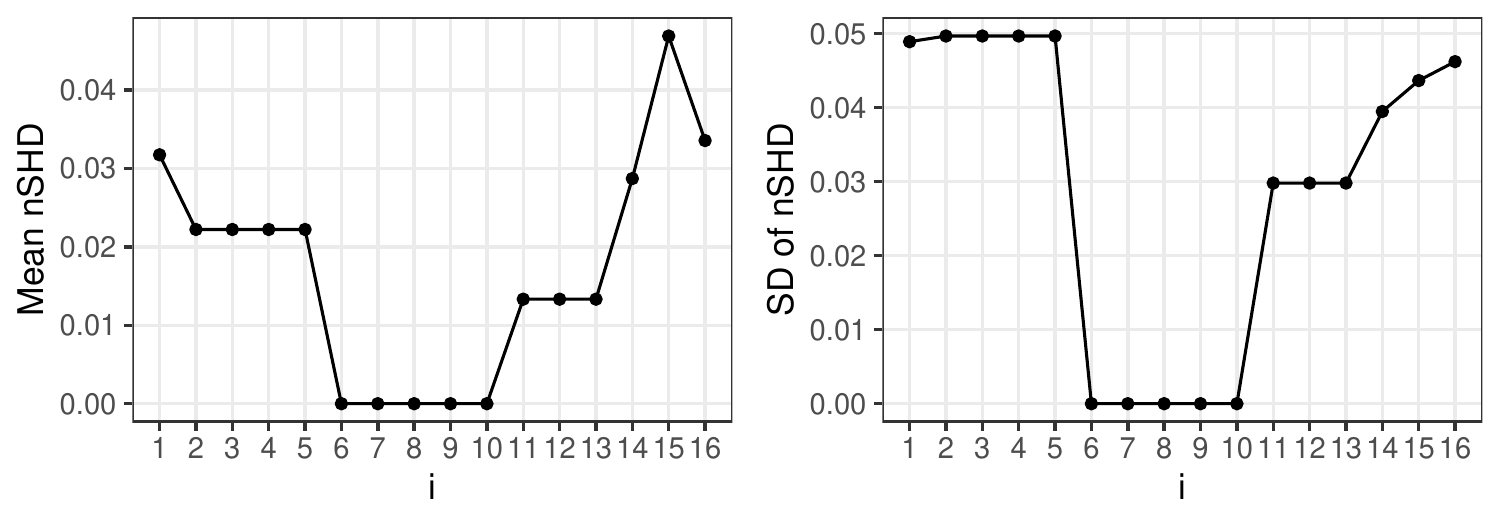}
    \caption{Five-pair rolling mean (left) and rolling standard deviation (right) of the distances $\operatorname{nSHD}_{\mathrm{mod}}$ between consecutive estimated graphs $\hat{\mathcal A}_{\gamma,k}$, $k\in\mathcal K$. Lower values of the mean indicate smaller average changes of the estimated graph.}
    \label{fig:k_stability}
\end{figure}

The left-hand panel of Figure~\ref{fig:k_stability} shows that the rolling mean is smallest for $k\in\{60, 70, 80, 90, \allowbreak 100\}$, a region in which the estimated graphs $\hat{\mathcal A}_{\gamma, k}$ coincide. The right-hand panel shows that the rolling standard deviation is smallest in this region as well. We therefore choose $k^*=80$, the midpoint of this region.

\paragraph{Confidence graphs.} For $k^*=80$, we compute the lower and upper confidence graphs from Section~\ref{sec:graph-recovery} at level $\alpha=0.1$, based on $1000$ multiplier bootstrap replications; see Figure~\ref{fig:food6}. The lower confidence graph coincides with $\hat{\mathcal A}_{\gamma,k^*}$ except for the edge VK$\to$BC. All remaining edges in the $\gamma$-signal graph are supported by the data at the chosen confidence level. The upper confidence graph contains two additional edges, LZ$\to$VA and VK$\to$VA, which our procedure neither establishes nor excludes, and it coincides with the graph estimated by \citet{KraliDavisonKluppelberg2026}. Hence, the CTC-based procedure recovers the structure found by the RMLM-based approach and, in addition, quantifies which of its edges are supported by the data. Since ancestral graphs are transitively closed, the two additional edges would in fact be implied by the paths LZ$\to$BC$\to$VA and VK$\to$BC$\to$VA if the $\gamma$-signal graph coincided with the ancestral graph.

Both directed edges between LZ and VK belong to the lower confidence graph, that is, both $\Gamma_{\mathrm{LZ},\mathrm{VK}}$ and $\Gamma_{\mathrm{VK},\mathrm{LZ}}$ significantly exceed $\gamma$. Since an ancestral graph cannot contain both directions, this does not indicate mutual causation; under the ancestral tail mixture structure, it rather suggests that the strong extremal dependence of the pair is driven by a common, possibly hidden, ancestor. This is in line with \citet{KraliDavisonKluppelberg2026}, who classified LZ and VK as indistinguishable, and, more generally, the inferred relationships are consistent with the groups~of nutrients $\{\mathrm{RET},\mathrm{VA}\}$, $\{\mathrm{AC},\mathrm{BC},\mathrm{VA}\}$ and $\{\mathrm{LZ},\mathrm{VK},\mathrm{BC},\mathrm{VA}\}$ identified by \citet{Boulin2026}.

\begin{figure}[t]
\centering
\begin{tikzpicture}[
  >=Stealth,
  vertex/.style={
    circle, draw=black, fill=white, thick,
    minimum size=8mm, inner sep=1pt
  },
  edge/.style={
    ->, shorten <=2pt, shorten >=2pt
  },
  f21/.style={edge, draw=black,     line width=1.5pt},
  f18/.style={edge, draw=black!65,  line width=1.10pt},
  f15/.style={edge, draw=black!45,  line width=0.90pt},
  f5/.style ={edge, draw=black!25,  line width=0.75pt},
  f1/.style ={edge, draw=black!10,  line width=0.65pt}
]

\node[vertex] (ret) at (-2,  0) {\footnotesize RET};
\node[vertex] (ac)  at ( 0,  0) {\footnotesize AC};
\node[vertex] (lz)  at ( 2,  0) {\footnotesize LZ};
\node[vertex] (va)  at (-2, -2) {\footnotesize VA};
\node[vertex] (bc)  at ( 0, -2) {\footnotesize BC};
\node[vertex] (vk)  at ( 2, -2) {\footnotesize VK};

\draw[f15]     (ac)  -- (va);
\draw[f21]     (ac)  -- (bc);
\draw[f5]      (lz)  -- (va);
\draw[f18]     (lz)  -- (bc);
\draw[f21]     (bc)  -- (va);
\draw[f18]     (vk)  -- (bc);

\draw[f21] ([xshift=-4pt]lz.south) -- ([xshift=-4pt]vk.north);
\draw[f21] ([xshift= 4pt]vk.north) -- ([xshift= 4pt]lz.south);

\draw[f21] ([xshift=-4pt]ret.south) -- ([xshift=-4pt]va.north);
\draw[f1] ([xshift= 4pt]va.north) -- ([xshift= 4pt]ret.south);

\draw[f5]
  (vk) to[out=230,in=315] (va);

\begin{scope}[shift={(4.2,-0.4)}]
 \node[draw=none, anchor=west] at (0,.5)
   {Edge frequency};
  \draw[f21] (0,0) -- (0.65,0)
    node[right,black] {\footnotesize 21 / 21};
  \draw[f18] (0,-0.4) -- (0.65,-0.4)
    node[right,black] {\footnotesize 18--19 / 21};
  \draw[f15] (0,-0.8) -- (0.65,-0.8)
    node[right,black] {\footnotesize 15 / 21};
  \draw[f5] (0,-1.2) -- (0.65,-1.2)
    node[right,black] {\footnotesize 5 / 21};
  \draw[f1] (0,-1.6) -- (0.65,-1.6)
    node[right,black] {\footnotesize 1 / 21};
\end{scope}
\end{tikzpicture}
\caption{Summary of the estimated $\gamma$-signal graphs $\hat{\mathcal A}_{\gamma,k}$, $k \in \mathcal K$, with edges colored according to the number of grid points for which they appear. Since edges in opposite directions can appear with different frequencies, we use a directed edge for each direction.}
  \label{fig:summary_graph}\vspace{-.2cm}
\end{figure}
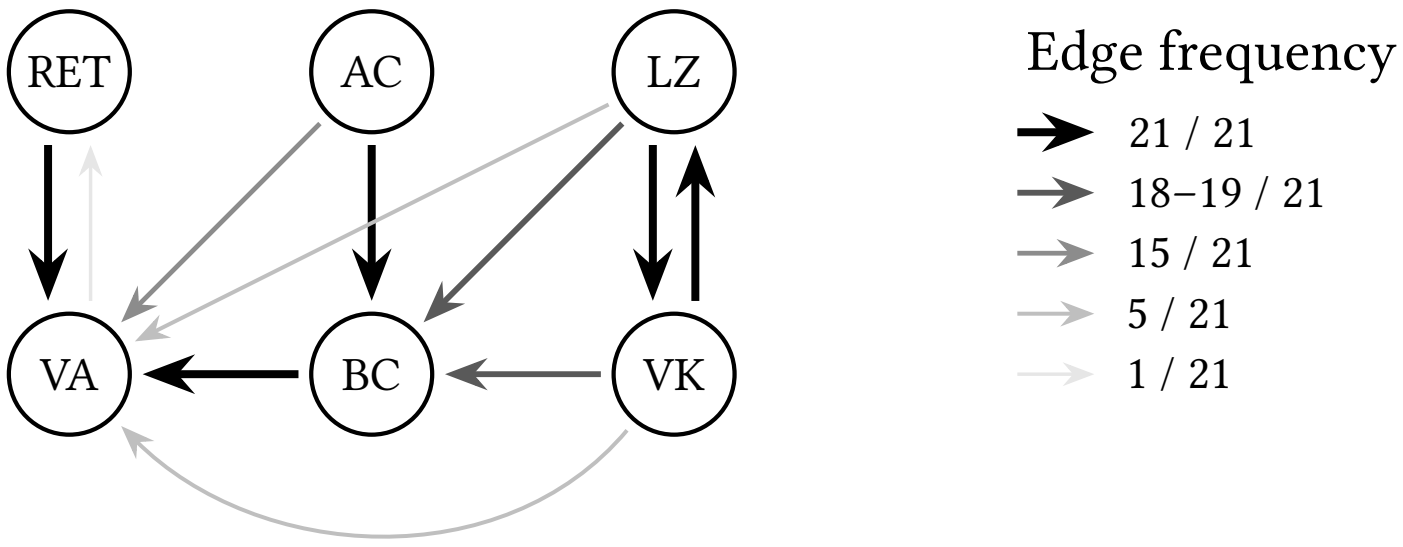

\paragraph{Sensitivity analysis.} 
Finally, we assess the sensitivity of the estimated graph with respect to the choice of $k$. Figure~\ref{fig:summary_graph} summarizes the estimated $\gamma$-signal graphs $\hat{\mathcal A}_{\gamma,k}$ for all grid points $k \in \mathcal K$, with edges colored according to the number of grid points for which they appear. The edges AC$\to$BC, BC$\to$VA, RET$\to$VA and both edges between LZ and VK are estimated for all $21$ grid points, and LZ$\to$BC and VK$\to$BC for $19$ and $18$ grid points, respectively. Hence, the main structure of the graph of \citet{KraliDavisonKluppelberg2026} is robust with respect to the choice of $k$. The edge AC$\to$VA appears for $15$ grid points, while the edges LZ$\to$VA and VK$\to$VA appear for only $5$ grid points, and VA$\to$RET for a single one. Remarkably, the edges that appear for only few grid points are precisely those whose presence is not established by the confidence graphs for $k^*$: both LZ$\to$VA and VK$\to$VA only belong to the upper confidence graph, and VA$\to$RET does not even belong to the latter, so that it can be excluded from the $\gamma$-signal graph at the chosen confidence level. The edge VK$\to$BC, absent for three grid points, is estimated but not confirmed by the lower confidence graph. The descriptive sensitivity analysis and the formal inference thus lead to consistent conclusions.

\begin{figure}[t]
        \centering
\begin{subfigure}{0.3\textwidth}
        \centering
\begin{tikzpicture}[
  >=Stealth,
  vertex/.style={
    circle, draw=black, fill=white, thick,
    minimum size=8mm, inner sep=1pt
  },
  edge/.style={
    shorten <=2pt,
    shorten >=2pt,
    line width=0.9pt
  },
  directed/.style={edge, ->},
  undirected/.style={edge, -, dashed}
]

\node[vertex] (ret) at (-2,  0) {\footnotesize RET};
\node[vertex] (ac)  at ( 0,  0) {\footnotesize AC};
\node[vertex] (lz)  at ( 2,  0) {\footnotesize LZ};
\node[vertex] (va)  at (-2, -2) {\footnotesize VA};
\node[vertex] (bc)  at ( 0, -2) {\footnotesize BC};
\node[vertex] (vk)  at ( 2, -2) {\footnotesize VK};

\draw[directed] (ret) -- (va);
\draw[directed] (ac)  -- (bc);
\draw[directed] (ac)  -- (va);
\draw[directed] (lz)  -- (bc);
\draw[directed] (vk)  -- (bc);
\draw[directed] (bc)  -- (va);

\draw[undirected] (lz) -- (vk);

\end{tikzpicture}
\end{subfigure}%
\hfill
\begin{subfigure}{0.3\textwidth}
        \centering
\begin{tikzpicture}[
  >=Stealth,
  vertex/.style={
    circle, draw=black, fill=white, thick,
    minimum size=8mm, inner sep=1pt
  },
  edge/.style={
    shorten <=2pt,
    shorten >=2pt,
    line width=0.9pt
  },
  directed/.style={edge, ->},
  undirected/.style={edge, -, dashed}
]

\node[vertex] (ret) at (-2,  0) {\footnotesize RET};
\node[vertex] (ac)  at ( 0,  0) {\footnotesize AC};
\node[vertex] (lz)  at ( 2,  0) {\footnotesize LZ};
\node[vertex] (va)  at (-2, -2) {\footnotesize VA};
\node[vertex] (bc)  at ( 0, -2) {\footnotesize BC};
\node[vertex] (vk)  at ( 2, -2) {\footnotesize VK};

\draw[directed] (ret) -- (va);
\draw[directed] (ac)  -- (bc);
\draw[directed] (ac)  -- (va);
\draw[directed] (lz)  -- (bc);
\draw[directed] (bc)  -- (va);

\draw[undirected] (lz) -- (vk);

\end{tikzpicture}
\end{subfigure}%
\hfill
\raisebox{-.385cm}{%
\begin{subfigure}{0.3\textwidth}
        \centering
\begin{tikzpicture}[
  >=Stealth,
  vertex/.style={
    circle, draw=black, fill=white, thick,
    minimum size=8mm, inner sep=1pt
  },
  edge/.style={
    shorten <=2pt,
    shorten >=2pt,
    line width=0.9pt
  },
  directed/.style={edge, ->},
  undirected/.style={edge, -, dashed}
]
\begin{scope}[yshift=-2cm]
\node[vertex] (ret) at (-2,  0) {\footnotesize RET};
\node[vertex] (ac)  at ( 0,  0) {\footnotesize AC};
\node[vertex] (lz)  at ( 2,  0) {\footnotesize LZ};
\node[vertex] (va)  at (-2, -2) {\footnotesize VA};
\node[vertex] (bc)  at ( 0, -2) {\footnotesize BC};
\node[vertex] (vk)  at ( 2, -2) {\footnotesize VK};

\draw[directed] (ret) -- (va);
\draw[directed] (ac)  -- (bc);
\draw[directed] (ac)  -- (va);
\draw[directed] (lz)  -- (bc);
\draw[directed] (lz)  -- (va);
\draw[directed] (vk)  -- (bc);
\draw[directed] (bc)  -- (va);
\draw[directed] (vk) to[bend left=28] (va);

\draw[undirected] (lz) -- (vk);
\end{scope}
\end{tikzpicture}
\end{subfigure}
}
\caption{Graph structures for the NHANES dietary intake data.
Left: Estimated ancestral graph $\hat{\mathcal A}_{\gamma,k^*}$ for $k^*=80$.
Middle and right: Lower and upper confidence graphs,
respectively, for $\gamma=0.95$ and $\alpha=0.1$. The upper confidence graph coincides with the graph estimated in \citet{KraliDavisonKluppelberg2026}.
The dashed edge indicates that both directed edges are included, whereas in \citet{KraliDavisonKluppelberg2026} the pair (LZ,VK) is classified as \emph{indistinguishable}.}
  \label{fig:food6}\vspace{-.0cm}
\end{figure}
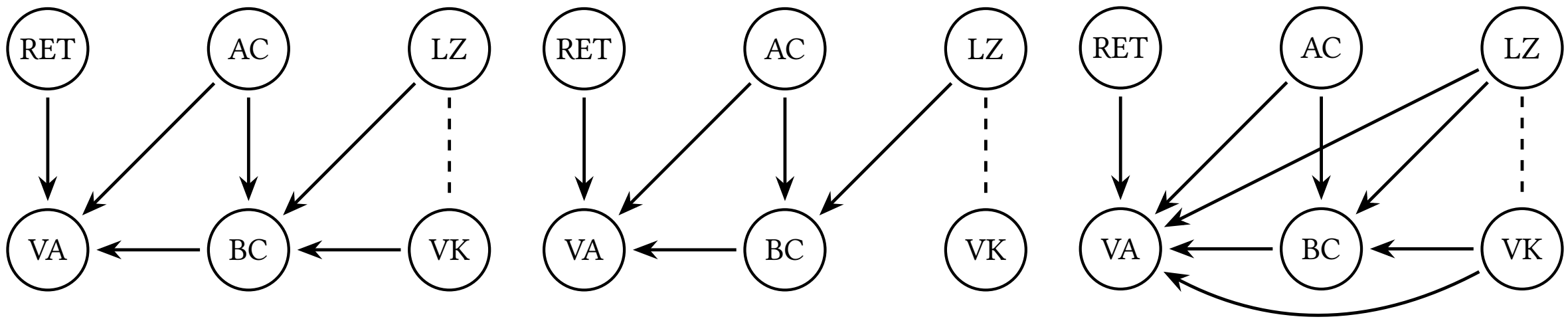

\section*{Declarations}

\begin{itemize}
\item\textbf{Ethical Approval.} Not applicable.
\item\textbf{Availability of supporting data. } Code for the case study is available in this Github repository \url{https://github.com/mariokrali/CTC_asymptotics}. The dataset used in Section~\ref{sec:data-application} is also publicly available on the website mentioned in the main text.
\item\textbf{Competing interests.} The authors declare that they have no conflict of interest.
\item\textbf{Funding.} Not applicable.
\item \textbf{Use of AI.} 
The authors used the large language models Claude Opus~5.5 (Anthropic) and ChatGPT~6 (OpenAI) as auxiliary tools for discussing technical aspects of the research, restructuring the paper, checking and streamlining proofs, improving code and language editing. In particular, the tools contributed to developing and refining the proofs of Lemmas~\ref{lem:linear-scm-condition-D} and~\ref{lem:linear-scm-bias}. The authors carefully reviewed and revised all outputs and take full responsibility for the content of the publication.
\end{itemize}

\appendix

\section{Proofs}
\label{sec:proofs}

\subsection{Proofs for Section~\ref{sec:bivariate-CTC}}

\begin{proof}[Proof of Lemma~\ref{lem:gamma-representation-bivariate}]
Throughout, for $h \in (0,1)$, we write $F_{V_2 \mid V_1 < h}(v) = \Prob(V_2 \le v \mid V_1 < h) = D(h,v)/h$, $v \in [0,1]$, for the CDF of $\mathcal L(V_2 \mid V_1 < h)$.
Suppose first that Condition~\ref{cond:conditional-weak-limit-bivariate} holds, and
let $H$ denote the CDF of the limiting distribution. By the
Portmanteau theorem,
\[
\frac{D(h,v)}{h} = F_{V_2 \mid V_1 < h}(v) \longrightarrow H(v),
\qquad h\downarrow0,
\]
at every continuity point $v$ of $H$. Since a distribution
function has at most countably many discontinuities, it follows that 
\begin{align*}
\partial_1^+ D(0,v)  = \lim_{h \downarrow 0} \frac{D(h,v)}{h}
=
H(v)
\end{align*}
exists for Lebesgue-almost every $v\in(0,1)$.

Conversely, suppose that the limit in \eqref{eq:pd0-bivariate} exists for Lebesgue-almost every $v\in(0,1)$. 
Let $(h_n)_{n\in\mathbb N}$ be an arbitrary sequence with $h_n\downarrow0$, and write $H_n = F_{V_2 \mid V_1 < h_n}$. Since each $H_n$ is supported on the compact
set $[0,1]$, Prohorov's theorem yields a subsequence $(h_{n_k})_{k\in\mathbb N}$ along which $H_{n_k}$ converges weakly to some CDF $H$.
Hence,
$
H_{n_k}(v)\longrightarrow H(v)
$
at every continuity point of $H$. On the other hand, by assumption,
$
H_{n_k}(v)\longrightarrow \partial_1^+ D(0,v) 
$
for Lebesgue-almost every $v\in(0,1)$. Consequently,
$
H(v)=\partial_1^+ D(0,v) 
$
for Lebesgue-almost every $v\in(0,1)$.
In particular, any two subsequential weak limits have CDFs that agree Lebesgue-almost everywhere. Since CDFs are nondecreasing and right-continuous, equality
Lebesgue-almost everywhere implies equality everywhere. 
Thus, all subsequential weak limits coincide. Since $(h_n){n\in\mathbb N}$ was arbitrary, it follows that $\mathcal L(V_2\mid V_1<h)$ converges weakly as $h\downarrow0$, with limiting CDF $H$ satisfying
\[
H(v)
=
\partial_1^+D(0,v)
\]
for Lebesgue-almost every $v\in(0,1)$.

It remains to establish the representation of $\Gamma$. By
Condition~\ref{cond:conditional-weak-limit-bivariate},
and 
since $x\mapsto1-x$ is bounded and continuous on $[0,1]$,
\[
\Gamma
=
\lim_{h\downarrow0}
\Exp[1-V_2\mid V_1<h]
=
1-\Exp[Z].
\]
Finally, since $Z$ takes values in $[0,1]$,
\[
1-\Exp[Z]
=
\int_0^1 H(v)\,\diff v = \int_0^1
\partial_1^+D(0,v)\,\diff v,
\]
which completes the proof.
\end{proof}

\begin{proof}[Proof of the claims in Examples~\ref{ex:evc-conditions} and~\ref{ex:evc-D}] 
Since $A$ is convex, its left derivative $A'(t-)$ exists for every $t\in(0,1]$, and $t \mapsto A'(t-)$ is left-continuous on $(0,1]$. Let
\[
L(x,y)=(x+y)\,A\Big(\frac{y}{x+y}\Big), \qquad x\ge0,\ y>0,
\]
denote the stable tail dependence function, such that $C(u,v)=\exp\{-L(-\log u,-\log v)\}$ for $u\in(0,1]$ and $v\in(0,1)$.

We first show that, for $x\ge0$ and $y>0$, the right-hand partial derivative of $L$ with respect to its first argument is given by
\begin{align}
\label{eq:evc-L-deriv}
\partial_1^+L(x,y) = A(t)-tA'(t-), \qquad t=\frac{y}{x+y}.
\end{align}
For $h>0$, let $t_h=y/(x+h+y)$, such that $t_h<t$, $t_h\uparrow t$ as $h\downarrow0$, and $(t_h-t)/h=-t/(x+h+y)$. Then
\[
\frac{L(x+h,y)-L(x,y)}{h}
=
A(t_h)+(x+y)\,\frac{A(t_h)-A(t)}{t_h-t}\,\frac{t_h-t}{h}
\longrightarrow
A(t)-tA'(t-),
\]
which is \eqref{eq:evc-L-deriv}.

Next, fix $u\in[0,1)$ and $v\in(0,1)$, and let $x=-\log(1-u)$, $y=-\log(1-v)$ and $x_h=-\log(1-u-h)$, such that $x_h\downarrow x$ and $(x_h-x)/h\to1/(1-u)$ as $h\downarrow0$. Since $D(u,v)=u+v-1+C(1-u,1-v)$, we obtain from \eqref{eq:evc-L-deriv}, continuity of $L$ and the mean value theorem applied to the exponential function that
\begin{align*}
\partial_1^+D(u,v)
&=
1-\lim_{h\downarrow0}\frac{\exp\{-L(x,y)\}-\exp\{-L(x_h,y)\}}{h}
\\&=
1-\exp\{-L(x,y)\}\,\partial_1^+L(x,y)\lim_{h\downarrow0}\frac{x_h-x}{h}
\\&=
1-\frac{C(1-u,1-v)}{1-u}\big\{A(t_{uv})-t_{uv}A'(t_{uv}-)\big\},
\qquad
t_{uv}=\frac{\log(1-v)}{\log(1-u)+\log(1-v)}.
\end{align*}
For $u=0$, we have $t_{0v}=1$, $A(1)=1$ and $C(1,1-v)=1-v$, whence
\[
\partial_1^+D(0,v)=v+A'(1-)(1-v), \qquad v\in(0,1).
\]
Condition~\ref{cond:conditional-weak-limit-bivariate} hence follows from Lemma~\ref{lem:gamma-representation-bivariate}. Moreover, for fixed $v\in(0,1)$, we have $t_{uv}\uparrow1$ as $u\downarrow0$, such that left-continuity of $t\mapsto A'(t-)$ at $1$, together with continuity of $A$ and $C$, implies $\lim_{u\downarrow0}\partial_1^+D(u,v)=\partial_1^+D(0,v)$. Hence, Condition~\ref{cond:continuity2-bivariate} is met.

It remains to show that $A'(1-) = 1-2\Phi(\{0\})$, with $\Phi$ the spectral measure of $A$. By the spectral representation of $A$, we have, for $h \in (0,1)$,
\[
\frac{A(1-h) - A(1)}{-h} = 2\int_{[0,1]} \frac{w - \max\{ (1-h)w,h(1-w) \}}h \, \Phi(\diff w).
\]
The integrand on the right-hand side is bounded in absolute value by 1 and, as $h \downarrow 0$, converges to $w$ for $w>0$ and to $-1$ for $w=0$. Hence, by dominated convergence,
\[
A'(1-) = \lim_{h \downarrow 0} \frac{A(1-h) - A(1)}{-h}
= 
2 \int_{(0,1]} w \, \Phi(\diff w) - 2 \Phi(\{0\})
= 1 - 2 \Phi(\{0\}),
\]
where the last equality follows from $\int_{[0,1]}w\,\Phi(\diff w)=1/2$.
\end{proof}

\begin{proof}[Proof of the claims in Examples~\ref{ex:gauss-conditions} and~\ref{ex:gauss-D}]
Recall that $D_\rho(u,v)=\Phi_\rho(\Phi^{-1}(u),\Phi^{-1}(v))$,
where $\Phi_\rho$ denotes the bivariate
normal CDF with mean vector zero and covariance matrix with rows $(1,\rho)$ and $(\rho,1)$.
A straightforward calculation shows that
\begin{align*}
\partial_1^+ D_\rho(u,v) = \partial_1 D_\rho(u,v)
&=
\Phi\left(
\frac{\Phi^{-1}(v)-\rho\Phi^{-1}(u)}
{\sqrt{1-\rho^2}}
\right), \qquad (u,v) \in (0,1)^2;
\end{align*}
see also \citet[Page 3056]{Omelka2009}.
As a consequence, for every $v \in (0,1)$, since $D_\rho(0,v)=0$,
\begin{align*}
\frac{D_\rho(h,v)}{h}
=
\frac{1}{h}\int_0^h \partial_1 D_\rho(s,v)\,\diff s.
\end{align*}
Since $\Phi^{-1}(s)\to-\infty$ as $s\downarrow0$, for every fixed
$v\in(0,1)$,
\begin{align}
\label{eq:convergence-gauss-pd}
\partial_1 D_\rho(s,v) = \Phi\left(
\frac{\Phi^{-1}(v)-\rho\Phi^{-1}(s)}
{\sqrt{1-\rho^2}}
\right)
\longrightarrow
\begin{cases}
0, & \rho<0,\\
v, & \rho=0,\\
1, & \rho>0.
\end{cases}
\end{align}
The corresponding Ces\`aro average has the same limit, i.e., 
\begin{align*}
\partial_1^+ D_\rho(0,v)
=
\lim_{h \downarrow 0} \frac{D_\rho(h,v)}{h}
=
\begin{cases}
0, & \rho<0,\\
v, & \rho=0,\\
1, & \rho>0,
\end{cases}
\qquad v\in(0,1).
\end{align*}
Thus, by Lemma~\ref{lem:gamma-representation-bivariate}, Condition~\ref{cond:conditional-weak-limit-bivariate} holds.

Finally, \eqref{eq:convergence-gauss-pd} also yields $\lim_{u \downarrow 0}\partial_1^+ D_\rho(u,v) = \partial_1^+ D_\rho(0,v)$ for every $v \in (0,1)$, and hence
Condition~\ref{cond:continuity2-bivariate} is met.
\end{proof}

\begin{proof}[Proof of Lemma~\ref{lem:ctc-asymptotic-independence}]
Let $D$ denote the survival copula corresponding to $C$, such  that
\[
\chi_U(C) = \lim_{u \downarrow 0} \frac{D(u,u)}u
\]
is the upper tail-dependence coefficient of $C$. It suffices to construct $D$, since $C$ is then absolutely continuous as well. We distinguish three cases with respect to $p$.

First, for $p\in\{0,1/2,1\}$, we can take $D=D_\rho$ as the Gaussian copula with correlation $\rho \in \{-1/2, 0,1/2\}$, respectively; see Example~\ref{ex:gauss-conditions}.

Next, let $p \in (0,1/2)$, and define $q=2p \in (0,1)$. Let $D$ be the checkerboard copula with bivariate density 
\[
d(u,v) = \frac1q \bm 1((u,v) \in I_1) + \frac1{1-q} \bm 1((u,v) \in I_2),
\]
where $I_1 = (0,q) \times (1-q,1)$ and $I_2 = (q,1) \times (0,1-q)$. 
It is immediate that both marginal densities of $d$ are equal to one,
so that $d$ is indeed a copula density. For every $u \in[0, q)$, we have 
\[
D(u,v) = 
\frac{1}q u (v-1+q) \bm 1(v \ge 1-q),
\]
which implies 
\[
\partial_1^+D(0,v)
=
\frac{1}q (v-1+q) \bm 1(v \ge 1-q).
\]
Hence, by Lemma~\ref{lem:gamma-representation-bivariate},
\[
\Gamma = \int_0^1 \partial_1^+D(0,v)\, \diff v = 
\frac{1}q
\int_{1-q}^1 (v-1+q) \, \diff v
= \frac q2 = p.
\]
Moreover, $D(u,u)=0$ for all sufficiently small $u$, and therefore $\chi_U(C) = \lim_{u \downarrow 0} D(u,u) / u = 0$.

Next, let $p \in (1/2,1)$, and define $q=2(1-p) \in (0,1)$. Let $D$ be the checkerboard copula with bivariate density 
\[
d(u,v) = \frac1q \bm 1((u,v) \in I_3) + \frac1{1-q} \bm 1((u,v) \in I_4),
\]
where $I_3 = (0,q) \times (0,q)$ and $I_4 = (q,1) \times (q,1)$.  Again, both marginal densities of $d$ are equal to one. For every $u\in[0,q)$, we have
\[
D(u,v) = 
\frac{ u v}q \bm 1(v  \le q) + u \bm 1(v > q),
\]
which implies
\[
\partial_1^+D(0,v)
=
\frac{v}q  \bm 1(v  \le q) +  \bm 1(v > q).
\]
Hence, by Lemma~\ref{lem:gamma-representation-bivariate},
\[
\Gamma = \int_0^1 \partial_1^+D(0,v)\, \diff v = 
\int_{0}^q \frac{v}q \, \diff v + \int_{q}^1  1 \, \diff  v 
= 1-\frac q2 = p.
\]
Moreover, $\chi_U(C) = \lim_{u \downarrow 0} D(u,u) / u = \lim_{u \downarrow 0} u/q = 0$.
\end{proof}

\begin{proof}[Proof of Theorem~\ref{theo:CTC-bivariate}]
The result is an immediate consequence of Theorem~\ref{theo:unknown-margins-bivariate} and Corollary~\ref{cor:unknown-margins-bivariate}, applied with $g(v)=1-v$ and $x=1$.
\end{proof}

\begin{proof}[Proof of Example~\ref{ex:evc-bias}]
Throughout, let $u\in(0,1)$ and $x=-\log(1-u)$, such that $x=u+O(u^2)$ as $u\downarrow0$. By definition of the survival copula,
\[
\int_0^1D(u,v)\,\diff v=u-\frac12+I(u), \qquad
I(u) = \int_0^1 C(1-u, 1-v)\, \diff v.
\]
Substituting $y=-\log(1-v)$ and then $z=x+y$, we obtain
\[
I(u) = \int_0^\infty \exp \Big\{ - y - (x+y) A\Big(\frac{y}{x+y}\Big) \Big\} \, \diff y
=
e^x \int_x^\infty  \exp\Big\{ - z - z A\Big( 1 - \frac{x}z\Big)\Big\} \, \diff z.
\]
For $b \ge 0$ to be specified below, let $R(h) = 1 - A(1-h) - a h + bh^2$ for $h \in [0,1]$, which is a bounded function. Since $zA(1-x/z) = z - ax + bx^2/z - zR(x/z)$, we have
\begin{align}
\label{eq:evc-bias-I}
I(u) = e^{(a+1)x} J(x),
\qquad
J(x)=\int_x^\infty \exp\Big\{ - 2z - b \frac{x^2}z + z R\Big( \frac{x}z\Big)\Big\} \, \diff z .
\end{align}
In both cases considered below, we show that $J(x) = e^{-2x}/2 - bx^2\log(1/x) + r(x)$, with a remainder $r(x)$ to be specified. Then \eqref{eq:evc-bias-I}, $e^{(a-1)x}=(1-u)^{1-a}=1-(1-a)u+O(u^2)$ and $e^{(a+1)x}x^2\log(1/x)=u^2\log(1/u)\{1+o(1)\}$ imply (for the latter, note that $x/u=1+O(u)$, such that $e^{(a+1)x}=1+O(u)$, $x^2=u^2\{1+O(u)\}$ and $\log(1/x)=\log(1/u)-\log(x/u)=\log(1/u)+O(u)$)
\begin{align}
\label{eq:evc-bias-D}
\int_0^1 D(u,v) \, \diff v
=
\frac{1+a}{2}\, u - b\,u^2\log(1/u)\{1+o(1)\} + O(u^2) + O(|r(x)|).
\end{align}
Since $\Gamma=(1+a)/2$ by Example~\ref{ex:evc-conditions}, the function $B(t) = t\int_0^1 D(1/t,v) \, \diff v - \Gamma$ is obtained by setting $u=1/t$ in \eqref{eq:evc-bias-D}, multiplying by $t$ and subtracting $\Gamma$.

\emph{Case 1: $A(1-h) = 1-ah+bh^2+o(h^2)$.} Then $R(h)=o(h^2)$ as $h\downarrow 0$, and $b \ge 0$, since $A(1-h)\ge 1-ah$ by convexity. First, for $z \ge x$, write $\exp(-bx^2/z) = 1 - bx^2/z + \rho_x(z)$, where $0 \le \rho_x(z) \le b^2x^4/(2z^2)$ by the elementary inequality $0\le e^{-w}-1+w\le w^2/2$ for $w\ge0$. Then
\[
\int_x^\infty \exp\Big( - 2z - b \frac{x^2}z \Big) \, \diff z
=
\int_x^\infty e^{-2z}\,\diff z
- bx^2 \int_x^\infty \frac{e^{-2z}}{z}\,\diff z
+ \int_x^\infty e^{-2z}\rho_x(z)\,\diff z.
\]
The first integral on the right-hand side equals $e^{-2x}/2$. For the second one, we have, for $x<1$ (which is no restriction, since $x\downarrow0$ for $u \downarrow 0$),
\[
\int_x^\infty \frac{e^{-2z}}{z}\,\diff z
=
\int_x^1 \frac{\diff z}{z}
+ \int_x^1 \frac{e^{-2z}-1}{z}\,\diff z
+ \int_1^\infty \frac{e^{-2z}}{z}\,\diff z
=
\log(1/x) + O(1),
\]
since $|e^{-2z}-1| \le 2z$. The third integral is bounded by $({b^2x^4}/{2})\int_x^\infty z^{-2}\,\diff z = b^2x^3/{2}$. Altogether,
\[
\int_x^\infty \exp\Big( - 2z - b \frac{x^2}z \Big) \, \diff z
=
\frac12 e^{-2x} - b x^2 \log(1/x) + O(x^2).
\]
Second, we show that the term involving $R$ in $J(x)$ is negligible. Let $\eps>0$, and choose $\delta \in (0,1]$ such that $|R(h)|\le \eps h^2$ for $h\le\delta$. Using $|e^w-1|\le |w|e^{|w|}$ and $e^{-bx^2/z} \le 1$, we obtain
\[
\Big| J(x) - \int_x^\infty \exp\Big( - 2z - b \frac{x^2}z \Big) \, \diff z \Big|
\le
\int_x^{x/\delta} z\|R\|_\infty e^{z\|R\|_\infty} \, \diff z
+
e\int_{x/\delta}^\infty \eps\frac{x^2}{z}e^{-2z} \, \diff z
\]
for $x$ small enough that $\eps x\delta\le1$. Since $\delta$ does not depend on $x$, the first integral on the right-hand side is bounded by
\[
\frac{x}{\delta}\cdot\frac{x}{\delta}\|R\|_\infty e^{\|R\|_\infty} = O(x^2)
\]
for $x\le\delta$, as the integrand is bounded by $(x/\delta)\|R\|_\infty e^{\|R\|_\infty}$ on an interval of length at most $x/\delta$. The second integral is bounded by $e\eps x^2\int_x^\infty z^{-1}e^{-2z}\,\diff z = e\eps x^2\{\log(1/x)+O(1)\}$, by the bound for the second integral in the first step. Combining both steps, we obtain
\[
\limsup_{x\downarrow0}\frac{|r(x)|}{x^2\log(1/x)}
=
\limsup_{x\downarrow0}\frac{|J(x)-e^{-2x}/2+bx^2\log(1/x)|}{x^2\log(1/x)}
\le e\eps.
\]
Since $\eps>0$ was arbitrary, $r(x)=o(x^2\log(1/x))$, which is $o(u^2\log(1/u))$ by the same arguments as for \eqref{eq:evc-bias-D}. Plugging this into \eqref{eq:evc-bias-D} and setting $u=1/t$, we obtain
\begin{align*}
B(t)
&=
t\int_0^1 D(1/t,v)\,\diff v-\frac{1+a}2
\\&=
-b\,\frac{\log t}{t}\{1+o(1)\} + O\Big(\frac1t\Big) + o\Big(\frac{\log t}{t}\Big)
=
-b\,\frac{\log t}{t} + o\Big(\frac{\log t}{t}\Big)
\end{align*}
as $t\to\infty$.

\emph{Case 2: $A(t)=1-a(1-t)$ for $t\in[t_0,1]$.} Choose $b=0$, such that $R(h)=0$ for $h\le h_0:=1-t_0$. Then, as in Case 1,
\[
|r(x)| = \Big|J(x)-\frac12e^{-2x}\Big|
\le
\int_x^{x/h_0} z\|R\|_\infty e^{z\|R\|_\infty} \, \diff z
= O(x^2) = O(u^2).
\]
Hence, $B(t) = O(1/t)$ by \eqref{eq:evc-bias-D}.
\end{proof}

\subsection{Proofs for Section~\ref{sec:multivariate-CTC}}

\begin{proof}[Proof of Theorem~\ref{theo:CTC-multivariate}]
Write $t_n=k/n$. The linearization and the marginal weak convergence are a
direct consequence of Theorem~\ref{theo:CTC-bivariate}. For the joint weak
convergence, note that the summands $n^{-1/2}\phi_{j\ell}^{(n)}(V_{ij},V_{i\ell})$
are centred by definition of $\Gamma_{j\ell}^{(n)}$ and bounded in absolute
value by $k^{-1/2}\to0$. By the Cram\'er--Wold device and the Lindeberg CLT,
it hence suffices to show that 
\[
\sigma_{j\ell;st}^{(n)}
=
\Exp\big[\phi_{j\ell}^{(n)}(V_j,V_\ell)\,\phi_{st}^{(n)}(V_s,V_t)\big]
=
\int_{[0,1]^d} f_n\,\diff\nu_{\{j,s\},t_n},
\]
converges to $\sigma_{j\ell;st}$, where $f_n(\bm v)=(1-v_\ell-\Gamma^{(n)}_{j\ell})(1-v_t-\Gamma^{(n)}_{st})$.
Condition~\ref{cond:multivariate} implies that
\[
\Gamma^{(n)}_{j\ell}=\int(1-v_\ell)\,\nu_{\{j\},t_n}(\diff\bm v)\to\Gamma_{j\ell},
\]
and likewise for $(s,t)$. Hence, $f_n\to f$ uniformly on $[0,1]^d$, where
$f(\bm v)=(1-v_\ell-\Gamma_{j\ell})(1-v_t-\Gamma_{st})$ is bounded and
continuous. Since the total masses of $\nu_{\{j,s\},t_n}$ are bounded by one,
the convergence to be shown follows from Condition~\ref{cond:multivariate}.
\end{proof}

\begin{proof}[Proof of Theorem~\ref{theo:bootstrap-multivariate}]
Write $t_n=k/n$.
Since $d$ is fixed, consistency of $\widehat\Sigma_n$
follows once we show that, for each
$(j,\ell),(s,t)\in\Pc_2$,
\begin{align} \label{eq:bootstrap-variance-consistent}
    \widehat\sigma_{j\ell;st}^{(n)}
    = \sigma_{j\ell;st} + o_\Prob(1).
\end{align}
To this end, recalling $\phi_{j\ell}^{(n)}$ from \eqref{eq:phin-multivariate}, consider the corresponding covariance estimator
based on the unobservable variables,
\[
    \widetilde\sigma_{j\ell;st}^{(n)}
    =
    \frac1n\sum_{i=1}^n
    \phi_{j\ell}^{(n)}(V_{ij},V_{i\ell})
    \phi_{st}^{(n)}(V_{is},V_{it}).
\]
The proof of \eqref{eq:bootstrap-variance-consistent} is finished once we show that
\begin{align} \label{eq:bootstrap-variance-consistent-sufficient}
    \widetilde\sigma_{j\ell;st}^{(n)}
    = \sigma_{j\ell;st} + o_\Prob(1),
    \qquad
    \widehat\sigma_{j\ell;st}^{(n)}
    -
    \widetilde\sigma_{j\ell;st}^{(n)}
    =o_\Prob(1).
\end{align}

We first establish the first convergence in \eqref{eq:bootstrap-variance-consistent-sufficient}.
By Condition~\ref{cond:multivariate} and the argument
in the proof of Theorem~\ref{theo:CTC-multivariate}, we have
$\lim_{n \to \infty} \Exp[\widetilde\sigma_{j\ell;st}^{(n)}]
    = \sigma_{j\ell;st}$.
Moreover, 
\[
    \Var(\widetilde\sigma_{j\ell;st}^{(n)})
    \le
    \frac1{nt_n^2} 
    \Prob(V_j\le t_n,V_s\le t_n)
    \le \frac1{nt_n}
    = \frac1k
    \to0,
\]
where we have used that $\phi_{j\ell}^{(n)}$ is bounded by $t_n^{-1/2}$.
Consequently, the first convergence in \eqref{eq:bootstrap-variance-consistent-sufficient} holds.

Next, for proving the second convergence in \eqref{eq:bootstrap-variance-consistent-sufficient}, it is sufficient to show that
\begin{align}
\label{eq:multivariate-bootstrap-l2}
    \frac1n\sum_{i=1}^n
    \Big|
    \widehat\phi_{j\ell}^{(n)}
    (\widehat V_{ij},\widehat V_{i\ell})
    -
    \phi_{j\ell}^{(n)}(V_{ij},V_{i\ell})
    \Big|^2
    =o_\Prob(1).
\end{align}
Indeed, the convergence already established on the left-hand side of \eqref{eq:bootstrap-variance-consistent-sufficient} implies that
\[
    \frac1n\sum_{i=1}^n
    \big\{\phi_{j\ell}^{(n)}(V_{ij},V_{i\ell})\big\}^2
    =O_\Prob(1).
\]
By \eqref{eq:multivariate-bootstrap-l2}, the same holds with 
$\phi_{j\ell}^{(n)}(V_{ij},V_{i\ell})$ replaced by
$\widehat\phi_{j\ell}^{(n)}
(\widehat V_{ij},\widehat V_{i\ell})$.
The claimed convergence then follows from the Cauchy--Schwarz inequality.

To establish \eqref{eq:multivariate-bootstrap-l2}, we apply \eqref{eq:two-param-L2}, established in the proof of Theorem~\ref{theo:bootstrap} 
under Condition~\ref{cond:continuity2-bivariate}, with $g(v)=1-v$, $T=1$ and the pair $(V_j,V_\ell)$.
Noting that
\[
    \frac1{\sqrt n}
    \phi_{j\ell}^{(n)}(V_{ij},V_{i\ell})
    =
    f_{ni}(1,\Gamma_{j\ell}^{(n)}),
    \qquad
    \frac1{\sqrt n}
    \widehat \phi_{j\ell}^{(n)}(\widehat V_{ij}, \widehat V_{i\ell})
    =
    \widehat f_{ni}(1,\widehat \Gamma_{j\ell}),
\]
we obtain
\begin{align*}
\frac1n\sum_{i=1}^n
    \Big|
    \widehat\phi_{j\ell}^{(n)}
    (\widehat V_{ij},\widehat V_{i\ell})
    -
    \phi_{j\ell}^{(n)}(V_{ij},V_{i\ell})
    \Big|^2
&\le
    2\sum_{i=1}^n
    \big\{
        \widehat f_{ni}(1,\Gamma_{j\ell}^{(n)})
        -f_{ni}(1,\Gamma_{j\ell}^{(n)})
    \big\}^2\\
&\hspace{3cm}+
    \frac{2}{k}
    (\widehat\Gamma_{j\ell}-\Gamma_{j\ell}^{(n)})^2
    \sum_{i=1}^n
    \bm 1(\widehat V_{ij}\le t_n).
\end{align*}
The first term is $o_\Prob(1)$ by
\eqref{eq:two-param-L2}, while the second is
$o_\Prob(1)$ since $\sum_{i=1}^n
    \bm 1(\widehat V_{ij}\le t_n) = k$ with probability 1
and $\widehat\Gamma_{j\ell}-\Gamma_{j\ell}^{(n)}
=o_\Prob(1)$.
This proves \eqref{eq:multivariate-bootstrap-l2}.

Finally, conditionally on the observations, the vectors
\[
    \xi_i\Big(
    \widehat\phi_{j\ell}^{(n)}
    (\widehat V_{ij},\widehat V_{i\ell})
    \Big)_{(j,\ell)\in\Pc_2},
    \qquad i\in[n],
\]
are independent and centered, with their normalized covariance matrix converging in probability to $\Sigma$ by \eqref{eq:bootstrap-variance-consistent}.
Moreover, the conditional Lindeberg condition is satisfied. Indeed, for any fixed linear combination of the components, the corresponding summands are uniformly bounded by $C\sqrt{n/k}$ for some $C>0$, while their empirical second moments are $O_\Prob(1)$. Since $k\to\infty$ and $\Exp[\xi_1^2]<\infty$, the conditional Lindeberg condition follows. The assertion is therefore a consequence of the
multivariate Lindeberg--Feller theorem, applied conditionally on the observations. More precisely, since the convergence of the conditional covariance matrix and the conditional Lindeberg condition hold in probability, every subsequence has a further subsequence along which both hold almost surely; along the latter, the Lindeberg--Feller theorem applies conditionally on the observations, which yields conditional weak convergence in probability in the sense of \citealp[Section~2.2.3]{Kosorok2008}.
\end{proof}

\subsection{Proofs for Section~\ref{sec:causal-models}}

\begin{proof}[Proof of Proposition~\ref{prop:ancestral-mixture}]
Fix $(j,\ell)\in\Pc_2$. Since $\{V_j<h\}$ and $\{V_j \le h\}$ coincide almost surely, the weak convergence $\nu_{\{j\},h}\weak\nu_j$ from Condition~\ref{cond:multivariate} and the continuous mapping theorem applied to the projection onto the $\ell$-th coordinate imply
\[
    \mathcal L(V_\ell\mid V_j < h)
    \weak
    \sum_{r\in\An(j)}
    w_{jr}\mathcal L(W_\ell^{(r)}).
\]
By definition of $\bm W^{(r)}$, and since $V_\ell$ is uniformly distributed on $(0,1)$,
\[
    \mathcal L(W_\ell^{(r)})
    =
    \begin{cases}
        \delta_0, & \ell\in\Des(r),\\
        \mathrm{Unif}(0,1), & \ell\notin\Des(r).
    \end{cases}
\]
Since $\ell\in\Des(r)$ if and only if $r\in\An(\ell)$, we obtain
\[
    \mathcal L(V_\ell\mid V_j< h)
    \weak
    p_{j\ell}\delta_0
    +(1-p_{j\ell})\mathrm{Unif}(0,1),
\]
which proves Condition~\ref{cond:conditional-weak-limit-bivariate}. For $Z_{j\ell}$ with this distribution, we have $\Exp[Z_{j\ell}] = (1-p_{j\ell})/2$ and $\Exp[Z_{j\ell}^2] = (1-p_{j\ell})/3$. Hence, by Lemma~\ref{lem:gamma-representation-bivariate}, $\Gamma_{j\ell} = 1-\Exp[Z_{j\ell}] = (1+p_{j\ell})/2$, and
\[
\Var(Z_{j\ell}) = \frac{1-p_{j\ell}}3 - \frac{(1-p_{j\ell})^2}{4} = \frac{(1-p_{j\ell})(1+3p_{j\ell})}{12}.
\]
Next, since $w_{jr}>0$ if and only if $r\in\An(j)$ and $\sum_{r \in [d]}w_{jr}=1$, we have $p_{j\ell}=1$ if and only if $\An(j)\subset\An(\ell)$, and $p_{j\ell}=0$ if and only if $\An(j)\cap\An(\ell)=\emptyset$. Since $j \in \An(j)$, the inclusion $\An(j)\subset\An(\ell)$ holds if and only if $j\in\An(\ell)$. This proves the characterizations of $\Gamma_{j\ell}=1$ and $\Gamma_{j\ell}=1/2$. Finally, the formula for $\chi_{js}$ follows from \eqref{eq:ancestral-mixture}, since each $\mu_r$ is a probability measure, and since $\nu_{\{j,s\}}([0,1]^d)=\lim_{h\downarrow 0} h^{-1}\Prob(V_j\le h, V_s\le h) = \chi_{js}$.
\end{proof}

\begin{proof}[Proof of Lemma~\ref{lem:linear-scm-multivariate}]
Let $\bar F_\varepsilon$ denote the common survival function of the innovations. For $h\in(0,1)$, $j\in[d]$ and $r\in\An(j)$, define
\begin{align}
\label{eq:xj-bj}
    x_j(h)=F_j^{-1}(1-h),
    \qquad
    B_j(h)=\{X_j>x_j(h)\},
    \qquad
    A_{jr}(h)=\{\beta_{r\to j}\varepsilon_r>x_j(h)\},
\end{align}
and let $B_j'(h)=\bigcup_{r\in\An(j)}A_{jr}(h)$. Such events are akin to those considered in the proof of Lemma~1 of~\cite{gnecco2021causal}. In what follows, we omit the argument $h$ and write, for instance, $x_j=x_j(h)$, $B_j=B_j(h)$, $A_{jr}=A_{jr}(h)$ and $B_j'=B_j'(h)$. Since $F_j$ is continuous, $\Prob(B_j)=h$, and $B_j$ coincides almost surely with $\{V_j\le h\}$. Moreover, $B_j'\subset B_j$, since $X_j\ge\beta_{r\to j}\varepsilon_r$ for every $r\in\An(j)$.

We treat both parts of \eqref{eq:ancestral-mixture} simultaneously. Let $I=\{j,s\}$ with $j,s\in[d]$, where possibly $j=s$, let $\mathcal C_I=\An(j)\cap\An(s)$ denote the set of common ancestors of $j$ and $s$, and write $B_I=B_j\cap B_s$, $w_{Ir}=w_{jr}\wedge w_{sr}$ and $A_{Ir}=A_{jr}\cap A_{sr}$ for $r\in\mathcal C_I$. With this notation, the claimed representation in \eqref{eq:ancestral-mixture} reads
\[
\nu_I=\sum_{r\in\mathcal C_I}w_{Ir}\,\mu_r,
\]
and, since $B_I$ coincides almost surely with $\{V_j\le h,V_s\le h\}$, the measure $\nu_{I,h}$ from Condition~\ref{cond:multivariate} can be written as $\nu_{I,h}=h^{-1}\Prob(\{\bm V\in\cdot\,\}\cap B_I)$. By the characterization of weak convergence of finite measures on $[0,1]^d$ through bounded Lipschitz continuous test functions, it suffices to show that, for every Lipschitz continuous $f:[0,1]^d\to\R$,
\begin{align}
\label{eq:pt-conv1}
    \lim_{h\downarrow0}\frac1h\Exp[f(\bm V)\bm 1(B_I)]
    =
    \sum_{r\in\mathcal C_I} w_{Ir}\Exp[f(\bm W^{(r)})].
\end{align}

We first show that, for all $j,s\in[d]$, $r\in\An(j)$ and $t\in\An(s)$ with $t\ne r$, as $h\downarrow0$,
\begin{align}
\label{eq:linear-scm-sbj}
    h^{-1}\Prob(A_{jr}) = w_{jr}+o(1),
    \qquad
    \Prob(A_{jr}\cap A_{st}) = O(h^2),
    \qquad
    \Prob(B_j\setminus B_j') = o(h).
\end{align}
By Lemma~S.1 of \citet{gnecco2021causal} and regular variation, $\Prob(X_j>x)\sim c_j\bar F_\varepsilon(x)$ as $x\to\infty$, with $c_j$ from \eqref{eq:weights-lsem}. Since $x_j\to\infty$ as $h\downarrow0$ and $\Prob(X_j>x_j)=h$, we obtain
\[
    \frac{\Prob(A_{jr})}{h}
    =
    \frac{\bar F_\varepsilon(x_j/\beta_{r\to j})}{\Prob(X_j>x_j)}
    \longrightarrow \frac{\beta_{r\to j}^\alpha}{c_j}=w_{jr},
\]
which is the first assertion in \eqref{eq:linear-scm-sbj}. For $r\ne t$, the events $A_{jr}$ and $A_{st}$ are independent, which together with the first assertion yields the second one. Consequently, by inclusion--exclusion and $\sum_{r\in\An(j)}w_{jr}=1$, we have $\Prob(B_j')=h+o(h)$. Since $B_j'\subset B_j$ and $\Prob(B_j)=h$, this implies the third assertion in \eqref{eq:linear-scm-sbj}.

Next, note that $A_{Ir}\subset B_I$ for $r\in\mathcal C_I$, which follows from $A_{jr}\subset B_j'\subset B_j$ and $A_{sr}\subset B_s'\subset B_s$. Conversely, let $\omega\in B_I\setminus\bigcup_{r\in\mathcal C_I}A_{Ir}$. We claim that either $\omega\in (B_j\setminus B_j')\cup(B_s\setminus B_s')$, or $\omega\in A_{jr}\cap A_{st}$ for some $r\in\An(j)$ and $t\in\An(s)$ with $r\ne t$. Indeed, suppose that $\omega\in B_j'\cap B_s'$. Then there exist $r\in\An(j)$ and $t\in\An(s)$ with $\omega\in A_{jr}\cap A_{st}$. If $r=t$, then $r\in\mathcal C_I$ and $\omega\in A_{jr}\cap A_{sr}=A_{Ir}$, which is impossible. Hence, $r\ne t$. Consequently, by the second and third assertion in \eqref{eq:linear-scm-sbj}, $\Prob\big(B_I\setminus\bigcup_{r\in\mathcal C_I}A_{Ir}\big)=o(h)$. Since, moreover, any two of the events $A_{Ir}$, $r\in\mathcal C_I$, intersect with probability $O(h^2)$ by the second assertion in \eqref{eq:linear-scm-sbj} (note that $A_{Ir}\cap A_{It}\subset A_{jr}\cap A_{jt}$ for $r\ne t$, with $r,t\in\mathcal C_I\subset\An(j)$), and since $f$ is bounded, we obtain
\begin{align}
\label{eq:rep1}
    \frac1h\Exp[f(\bm V)\bm 1(B_I)]
    =
    \frac1h\Exp\Big[f(\bm V)\bm 1\Big(\bigcup_{r\in\mathcal C_I}A_{Ir}\Big)\Big]+o(1)
    =
    \sum_{r\in\mathcal C_I}\frac1h\Exp[f(\bm V)\bm 1(A_{Ir})]+o(1).
\end{align}

We next argue why one may replace $\bm V$ by $\bm W^{(r)}$ in \eqref{eq:rep1}. To this end, fix $r\in\mathcal C_I$. For $m\in\Des(r)$, we have $X_m\ge\beta_{r\to m}\varepsilon_r$, such that, on $A_{Ir}\subset A_{jr}$, $X_m\ge c\,x_j$ with $c=\min_{m\in\Des(r)}\beta_{r\to m}/\beta_{r\to j}>0$. Hence, on $A_{Ir}$, $\max_{m\in\Des(r)}V_m\le\max_{m\in\Des(r)}\{1-F_m(c\,x_j)\}=o(1)$. Since $\bm V$ and $\bm W^{(r)}$ differ only in the coordinates $m\in\Des(r)$, Lipschitz continuity of $f$, with Lipschitz constant $L$ with respect to the maximum norm, yields
\begin{align*}
    \frac1h\Exp\big[|f(\bm V)-f(\bm W^{(r)})|\bm 1(A_{Ir})\big]
    &\le \frac Lh\Exp\Big[\max_{m\in\Des(r)}V_m\,\bm 1(A_{Ir})\Big]
    \\&\le L\max_{m\in\Des(r)}\{1-F_m(c\,x_j)\}\cdot\frac{\Prob(A_{jr})}{h}=o(1),
\end{align*}
where we used the first assertion in \eqref{eq:linear-scm-sbj}.

Finally, $\bm W^{(r)}$ is a function of the innovations $\varepsilon_t$, $t\ne r$, and hence independent of $A_{Ir}$. Moreover, $A_{Ir}=\{\varepsilon_r>\max(x_j/\beta_{r\to j},x_s/\beta_{r\to s})\}$, such that $\Prob(A_{Ir})=\Prob(A_{jr})\wedge\Prob(A_{sr})$, and hence $h^{-1}\Prob(A_{Ir})\to w_{Ir}$ by the first assertion in \eqref{eq:linear-scm-sbj}. Therefore, for every $r \in \mathcal C_I$,
\begin{align*}
    \frac1h\Exp[f(\bm V)\bm 1(A_{Ir})]
    &=
    \frac1h\Exp[f(\bm W^{(r)})\bm 1(A_{Ir})]+o(1)
    \\&=
    \frac{\Prob(A_{Ir})}{h}\Exp[f(\bm W^{(r)})]+o(1)
    \longrightarrow
    w_{Ir}\Exp[f(\bm W^{(r)})],
\end{align*}
which, together with \eqref{eq:rep1}, proves \eqref{eq:pt-conv1}.
\end{proof}

\begin{proof}[Proof of Lemma~\ref{lem:linear-scm-condition-D}]
Fix $(j,\ell)\in\Pc_2$, let $D_{j\ell}$ be the copula of $(V_j, V_\ell)$, and fix $v \in (0,1)$. By Lemma~\ref{lem:linear-scm-multivariate} and Proposition~\ref{prop:ancestral-mixture}, the conditional tail limit CDF $H_{j\ell}$ is continuous on $(0,1)$. Hence, by Lemma~\ref{lem:gamma-representation-bivariate},
\[
    \partial_1^+D_{j\ell}(0,v)
    =
    H_{j\ell}(v)
    =
    p_{j\ell}+(1-p_{j\ell})v.
\]
It is therefore sufficient to show that $\partial_1^+D_{j\ell}(u,v)$ exists for all
sufficiently small $u>0$, say $u\in(0,\eps_v)$, and that
\[
    \lim_{u\downarrow0}\partial_1^+D_{j\ell}(u,v)
    =
    p_{j\ell}+(1-p_{j\ell})v.
\]
For that purpose, first note that $X_j$ admits a density $f_j$, since the ancestral representation in \eqref{eq:linear-scm-ancestral} shows that it is
a sum of independent random variables with densities. Consequently, the measure
\[
    B\mapsto\Prob(X_j\in B,V_\ell\le v)
\]
is absolutely continuous with respect to Lebesgue measure and admits a density $q_v$, with $0\le q_v\le f_j$ almost everywhere. We can therefore write
\[
    D_{j\ell}(u,v)
    =
    \Prob(V_j\le u,V_\ell\le v)
    =
    \Prob(X_j\ge x_j(u),V_\ell\le v)
    =
    \int_{x_j(u)}^\infty q_v(x)\,\diff x,
\]
where $x_j(u)=F_j^{-1}(1-u)$.

We will show below that $f_j$ and $q_v$ admit eventually continuous versions, with
$f_j$ eventually strictly positive. Consequently, $x_j(u)$ is continuously differentiable for all sufficiently small $u>0$, and differentiation of the preceding integral yields
\[
    \partial_1^+D_{j\ell}(u,v)
    =
    \frac{q_v(x_j(u))}{f_j(x_j(u))}.
\]
We will further show that
\begin{align}
    \frac{q_v(x)}{f_j(x)}
    \longrightarrow
    p_{j\ell}+(1-p_{j\ell})v,
    \qquad x\to\infty.
    \label{eq:linear-scm-conditional-density-limit2}
\end{align}
Since $x_j(u)\to\infty$ as $u\downarrow0$, this implies
\[
    \lim_{u\downarrow0}
    \partial_1^+D_{j\ell}(u,v)
    =
    p_{j\ell}+(1-p_{j\ell})v,
\]
as required for Condition~\ref{cond:continuity2-bivariate}.

It remains to establish eventual continuity of $f_j$ and $q_v$, as well as
\eqref{eq:linear-scm-conditional-density-limit2}.
In order to treat both densities simultaneously, we allow for $v\in(0,1]$ in what follows, and note that $q_1=f_j$, since $V_\ell\le1$. For $r\in \An(j)$, let 
\[
    Y_r=\beta_{r\to j}\varepsilon_r,
    \qquad
    S_r=\sum_{t\in \An(j)\setminus\{r\}}Y_t,
\]
such that $X_j=Y_r+S_r$, with $Y_r$ and $S_r$ independent. The density $g_r$ of $Y_r$ is given by
\[
    g_r(x)
    =
    \frac1{\beta_{r\to j}}
    f_\varepsilon\left(
        \frac{x}{\beta_{r\to j}}
    \right),
    \qquad x>0.
\]
Regular variation implies that
\begin{align}
    g_r(x)
    &\sim
    \beta_{r\to j}^{\alpha}f_\varepsilon(x),
    \qquad x\to\infty.
    \label{eq:linear-scm-density-tail}
\end{align}
Moreover, by Theorem~1.1 of \citet{BinghamGoldieOmey2006}, the convolution
of two regularly varying densities is asymptotically equivalent to their sum.
Applying this result repeatedly, we obtain
\begin{align}
    f_j(x)
    &\sim
    \sum_{r\in \An(j)}g_r(x)
    \sim
    c_jf_\varepsilon(x),
    \qquad x\to\infty,
    \label{eq:linear-scm-marginal-density}
\end{align}
with $c_j$ from \eqref{eq:weights-lsem}. In particular, $f_j$ is eventually strictly positive.

We next establish eventual continuity of $q_v$, $v\in(0,1]$. Let $m=|\An(j)|$ and fix an arbitrary ordering of $\An(j)$. Choose $M>0$ (depending on $v$) sufficiently
large such that all densities $g_r$, $r\in\An(j)$, are continuous on $[M,\infty)$ and
\begin{align}
\label{eq:beta-quotient}
    \frac{\beta_{r\to\ell}}{\beta_{r\to j}}M
    \ge z_v,
    \qquad r\in\An(j)\cap\An(\ell),
\end{align}
where $z_v=F_\ell^{-1}(1-v)$ for $v<1$ and $z_1=0$. For $r\in\An(j)$, let
\[
    C_r=\prod_{t<r}\bm 1(Y_t\le M),
    \qquad
    h_r(y)=g_r(y)\bm 1(y>M),
\]
where $t<r$ refers to the chosen ordering. Since $X_j=\sum_{r\in\An(j)}Y_r$,
the events $\{Y_r>M,C_r=1\}$ are pairwise disjoint, and their union contains
$\{X_j>mM\}$. Consequently, for every Borel set $B\subset(mM,\infty)$,
\[
\Prob(X_j\in B,V_\ell\le v)
=
\sum_{r\in\An(j)}
    \Exp\big[
        \bm 1(X_j\in B) \bm 1(V_\ell\le v)
        \bm 1(Y_r>M)C_r
    \big].
\]
If $r\in\An(\ell)$, then on $\{Y_r>M\}$, by \eqref{eq:beta-quotient},
\[
    X_\ell
    \ge
    \frac{\beta_{r\to\ell}}{\beta_{r\to j}}Y_r
    \ge z_v,
\]
and hence $\bm 1(V_\ell\le v)=1$. Therefore,  for all $r\in \An(j)$
\[
\bm 1(Y_r>M)\bm 1(V_\ell\le v) 
=
\bm 1(Y_r>M) J_r, \qquad J_r=
    \begin{cases}
        1, & r\in\An(\ell),\\
        \bm 1(V_\ell\le v), & r\notin\An(\ell).
    \end{cases}
\]
Note that $J_r=1$ if $v=1$. Since $C_r$, $J_r$, and $S_r$ depend only
on innovations other than $\varepsilon_r$, they are jointly independent of $Y_r$.
Integrating out $Y_r$, we obtain
\begin{align*}
    \Prob(X_j\in B,V_\ell\le v)
    &=
    \sum_{r\in\An(j)}
    \Exp\big[
        \bm 1(Y_r+S_r\in B)
        \bm 1(Y_r>M)C_rJ_r
    \big]\\
    &=
    \int_B
    \sum_{r\in\An(j)}
    \Exp\big[
        h_r(x-S_r)C_rJ_r
    \big]\,\diff x.
\end{align*}
Thus,
\begin{align}
    q_v(x)
    =
    \sum_{r\in\An(j)}
    \Exp\big[
        h_r(x-S_r)C_rJ_r
    \big],
    \qquad x>mM,
    \label{eq:linear-scm-q-continuous}
\end{align}
defines a version of $q_v$.

Each $h_r$ is bounded and continuous except possibly at $M$. If $m\ge2$, then $S_r$
has an absolutely continuous distribution, so that
$\Prob(x-S_r=M)=0
$
for every fixed $x>mM$. Dominated convergence therefore shows that the right-hand side of
\eqref{eq:linear-scm-q-continuous}
is continuous in $x$. If $m=1$, continuity follows directly from $h_r(x)=g_r(x)$ for $x>M$.

It remains to establish
\eqref{eq:linear-scm-conditional-density-limit2}.
For $K>0$, $r\in \An(j)$ and $v\in(0,1]$, let $q_{v,r,K}$ denote the density of the measure
\[
    B\mapsto
    \Prob(X_j\in B,V_\ell\le v,S_r\le K).
\]
On the event $\{X_j>K+M\}\cap\{S_r\le K\}$, we have $Y_r=X_j-S_r>M$, and hence $\bm 1(V_\ell\le v)=J_r$, as shown above. Since $S_r$ and $J_r$ are independent of $Y_r$, integrating out $Y_r$ as in the derivation of \eqref{eq:linear-scm-q-continuous} shows that
\[
    q_{v,r,K}(x)
    =
    \Exp\big[
        g_r(x-S_r)\bm 1(S_r\le K)J_r
    \big],
    \qquad x>K+M,
\]
defines a version of $q_{v,r,K}$, which is continuous by continuity of $g_r$ on $[M,\infty)$ and dominated convergence.
By \eqref{eq:linear-scm-density-tail} and the uniform convergence theorem for regularly varying functions,
\[
    \lim_{x \to \infty }\frac{g_r(x-s)}{f_\varepsilon(x)}
    =
    \beta_{r\to j}^{\alpha}
\]
uniformly for $s\in[0,K]$ \citep[Proposition~0.5]{resnick1987extreme}.
Consequently, for every $v\in(0,1]$,
\begin{align}
    \lim_{x \to \infty }\frac{q_{v,r,K}(x)}{f_\varepsilon(x)}
    &=
    \beta_{r\to j}^{\alpha}
    \Exp[\bm 1(S_r\le K)J_r],
    \label{eq:linear-scm-truncated-density}
\end{align}
where $J_r=1$ for $v=1$.

Next, observe that, for distinct $r,t\in \An(j)$,
\[
    S_r+S_t
    = X_j + (X_j - Y_r -Y_t)
    =
    X_j+\sum_{q\in \An(j)\setminus\{r,t\}}Y_q
    \ge X_j.
\]
Consequently, the events $\{S_r\le K\}$, $r\in \An(j)$, are mutually exclusive on
$\{X_j>2K\}$. For every Borel set $B\subset(2K,\infty)$, we therefore have
\begin{align*}
    &\Prob(X_j\in B,V_\ell\le v)
    -
    \sum_{r\in \An(j)}
    \Prob(X_j\in B,V_\ell\le v,S_r\le K)
    \\
    &\quad=
    \Prob\big(
        X_j\in B,V_\ell\le v,
        S_r>K\text{ for all }r\in \An(j)
    \big)
    \\
    &\quad\le
    \Prob\big(
        X_j\in B,
        S_r>K\text{ for all }r\in \An(j)
    \big),
\end{align*}
and the right-hand side equals the left-hand side with $v=1$. Passing to densities, we obtain
\begin{align}
\label{eq:density-inequality}
    0\le
    q_v(x)-\sum_{r\in\An(j)}q_{v,r,K}(x)
    \le
    q_1(x)-\sum_{r\in\An(j)}q_{1,r,K}(x)
\end{align}
for Lebesgue-almost every $x>2K$. Since all densities admit eventually
continuous versions, this inequality holds for all sufficiently large $x$.

Combining \eqref{eq:density-inequality} with $q_1=f_j$, \eqref{eq:linear-scm-marginal-density} and \eqref{eq:linear-scm-truncated-density},
applied with $v$ and with $1$ in place of $v$, we obtain, for every $K>0$,
\begin{align*}
    \sum_{r\in\An(j)}
    \beta_{r\to j}^{\alpha}
    \Exp[\bm 1(S_r\le K)J_r]
    \le
    \liminf_{x\to\infty}
    \frac{q_v(x)}{f_\varepsilon(x)}
    &\le
    \limsup_{x\to\infty}
    \frac{q_v(x)}{f_\varepsilon(x)}\\
    &\le
    \sum_{r\in\An(j)}
    \beta_{r\to j}^{\alpha}
    \big\{
        \Exp[\bm 1(S_r\le K)J_r]
        +\Prob(S_r>K)
    \big\}.
\end{align*}
Letting $K\to\infty$, both bounds
converge to
\begin{align*}
    \sum_{r\in\An(j)}
    \beta_{r\to j}^{\alpha}\Exp[J_r]
    =
    \sum_{r\in\An(j)\cap\An(\ell)}
    \beta_{r\to j}^{\alpha}
    +
    v\sum_{r\in\An(j)\setminus\An(\ell)}
    \beta_{r\to j}^{\alpha}
    =
    c_j\{p_{j\ell}+(1-p_{j\ell})v\}.
\end{align*}
Consequently,
\[
    \frac{q_v(x)}{f_\varepsilon(x)}
    \longrightarrow
    c_j\{p_{j\ell}+(1-p_{j\ell})v\}.
\]
Together with
\eqref{eq:linear-scm-marginal-density},
this establishes
\eqref{eq:linear-scm-conditional-density-limit2}
and completes the proof.
\end{proof}

\begin{proof}[Proof of Lemma \ref{lem:linear-scm-bias}]
Fix $(j,\ell)\in\Pc_2$ and $\delta\in(0,\alpha)$ with $\delta\le\min\{\beta,1\}$. Recall the notation $w_{jr}$, $x_j=x_j(h)$, $B_j=B_j(h)$,
$A_{jr}=A_{jr}(h)$, and $B_j'=B_j'(h)$ from \eqref{eq:weights-lsem}, \eqref{eq:xj-bj}, and the proof of Lemma~\ref{lem:linear-scm-multivariate} and, as there, we omit the argument $h$.
We first establish the quantitative bounds, for $r\in\An(j)$ and as $h\downarrow0$,
\begin{align}
\label{eq:linear-scm-bias-sbj}
    h^{-1}\Prob(A_{jr})
    =
    w_{jr}+O(h^{\delta/\alpha}),
    \qquad
    \Prob(B_j\setminus B_j')
    =
    O(h^{1+\delta/\alpha});
\end{align}
see the first and third assertion in \eqref{eq:linear-scm-sbj} for the non-quantitative versions.
For $r\in\An(j)$, let
\[
    Y_r=\beta_{r\to j}\varepsilon_r,
    \qquad
    S_r=\sum_{t\in\An(j)\setminus\{r\}}Y_t,
\]
so that $X_j=Y_r+S_r$.
By assumption,
\begin{align}
    \Prob(Y_r>x)
    &=
    C\beta_{r\to j}^{\alpha}x^{-\alpha}
    +O(x^{-\alpha-\beta}),
    \qquad x\to\infty.
    \label{eq:linear-scm-bias-innovation}
\end{align}
In particular, since $\delta<\alpha$,
$\Exp[S_r^\delta]<\infty$.

We claim that
\begin{align}
    &\Prob\big(
        X_j>x,\,
        Y_r\le x\text{ for all }r\in\An(j)
    \big)
    =
    O(x^{-\alpha-\delta}).
    \label{eq:linear-scm-bias-remainder}
\end{align}
If $|\An(j)|=1$, this probability is zero. Otherwise, let $m=|\An(j)|\ge2$ and
$a=1-1/m \in [1/2,1)$. On the event in \eqref{eq:linear-scm-bias-remainder},
at least one $Y_r$ exceeds $x/m$. For this $r$, either $S_r\le ax$ or
$S_r>ax$. Consequently, since $Y_r = X_j - S_r$, the probability
in \eqref{eq:linear-scm-bias-remainder} is bounded by
\begin{align} \label{eq:two-sums}
    &\sum_{r\in\An(j)}
    \Prob(x-S_r<Y_r\le x,S_r\le ax)
    +
    \sum_{r\in\An(j)}
    \Prob(Y_r>x/m,S_r>ax).
\end{align}
For the first sum, independence of $Y_r$ and $S_r$ yields
\begin{align*}
    &\Prob(x-S_r<Y_r\le x,S_r\le ax)=
    \Exp\Big[
        \big\{
            \bar F_{Y_r}(x-S_r)
            -\bar F_{Y_r}(x)
        \big\}
        \bm 1(S_r\le ax)
    \Big].
\end{align*}
By \eqref{eq:linear-scm-bias-innovation},
uniformly for $s\in[0,ax]$,
\begin{align*}
    &\bar F_{Y_r}(x-s)-\bar F_{Y_r}(x)=
    C\beta_{r\to j}^{\alpha}x^{-\alpha}
    \big\{(1-s/x)^{-\alpha}-1\big\}
    +O(x^{-\alpha-\beta}).
\end{align*}

To bound the first term, observe that $y\mapsto \varphi(y) := (1-y)^{-\alpha}-1$ is convex on $[0,a]$ and vanishes at zero.
Hence, for $y\in[0,a]$,
\begin{align*}
    (1-y)^{-\alpha}-1
    &\le \Big(1-\frac ya\Big) \varphi(0) + \frac ya \varphi(a)  =
    \big\{ (1-a)^{-\alpha}-1 \big\} \frac ya
    \le 
    \big\{ (1-a)^{-\alpha}-1 \big\} \Big(\frac ya\Big)^\delta,
\end{align*}
where the last inequality uses $0<\delta\le1$. Thus, uniformly for $s\in[0,ax]$,
\begin{align*}
    0
    &\le
    \bar F_{Y_r}(x-s)-\bar F_{Y_r}(x)
    \le
    C_1x^{-\alpha-\delta}s^\delta
    +C_2x^{-\alpha-\beta}
\end{align*}
for suitable finite constants $C_1,C_2$. It follows that
\begin{align*}
    &\Prob(x-S_r<Y_r\le x,S_r\le ax)\le
    C_1x^{-\alpha-\delta}\Exp[S_r^\delta]
    +C_2x^{-\alpha-\beta}=
    O(x^{-\alpha-\delta}),
\end{align*}
where we used $\delta\le\beta$. For the second sum in \eqref{eq:two-sums}, independence of $Y_r$ and $S_r$ gives
\begin{align*}
    &\Prob(Y_r>x/m,S_r>ax)=
    \Prob(Y_r>x/m)\Prob(S_r>ax).
\end{align*}
Both probabilities are $O(x^{-\alpha})$, the latter by the union bound
\[
    \Prob(S_r>y)
    \le
    \sum_{t\in\An(j)\setminus\{r\}}
    \Prob\Big(
        Y_t>\frac{y}{m-1}
    \Big)
    =
    O(y^{-\alpha}).
\]
Consequently,
$
    \Prob(Y_r>x/m,S_r>ax)
    =
    O(x^{-2\alpha})
    =
    O(x^{-\alpha-\delta}),
$
since $\delta<\alpha$. This establishes
\eqref{eq:linear-scm-bias-remainder}.

Next, independence implies that the probability of two or more distinct $Y_r$
exceeding $x$ is $O(x^{-2\alpha})$. Therefore, defining $N_j(x) = |\{r \in \An(j) : Y_r > x\}|$, \eqref{eq:linear-scm-bias-remainder} and \eqref{eq:linear-scm-bias-innovation} yield
\begin{align}
    \Prob(X_j>x)
    &=
    \Prob(X_j>x, N_j(x) = 0) + \Prob(N_j(x) = 1) + \Prob(N_j(x) \ge 2)
    \notag\\ &=
    \Prob(N_j(x) = 1) + O(x^{-\alpha-\delta}) 
    \notag\\ &=
    \sum_{r\in\An(j)}\Prob(Y_r>x)
    +O(x^{-\alpha-\delta})\notag\\
    &=
    Cc_jx^{-\alpha}
    +O(x^{-\alpha-\delta}),
    \label{eq:linear-scm-bias-marginal}
\end{align}
with $c_j$ from \eqref{eq:weights-lsem}.
In particular, since $\bar F_j(x_j)=h$, we obtain that
\begin{align}
\label{eq:xjh-equivalence}
    x_j\sim(Cc_j/h)^{1/\alpha},
    \qquad h\downarrow0.
\end{align}
Combining
\eqref{eq:linear-scm-bias-innovation}
and \eqref{eq:linear-scm-bias-marginal},
we obtain
\begin{align*}
    \frac{\Prob(A_{jr})}{h}
    =
    \frac{\Prob(Y_r>x_j)}
         {\Prob(X_j>x_j)}
    = w_{jr} + O((x_j)^{-\delta}) = 
    w_{jr}+O(h^{\delta/\alpha}),
\end{align*}
which proves the first assertion in \eqref{eq:linear-scm-bias-sbj}.
Moreover, substituting $x=x_j$
in \eqref{eq:linear-scm-bias-remainder}
yields the second assertion in \eqref{eq:linear-scm-bias-sbj}.

We now repeat the argument for proving \eqref{eq:pt-conv1} in the proof of Lemma~\ref{lem:linear-scm-multivariate}, with $s=j$ and the Lipschitz function $f(\bm v)=1-v_\ell$, while keeping track of the rates. Since $\delta<\alpha$, we have $O(h)=O(h^{\delta/\alpha})$. First, since $0\le f\le1$, the second assertions in \eqref{eq:linear-scm-bias-sbj} and \eqref{eq:linear-scm-sbj} imply that \eqref{eq:rep1} holds with $o(1)$ replaced by $O(h^{\delta/\alpha})$. Second, if $r\notin\An(\ell)$, then $f(\bm V)=f(\bm W^{(r)})$. If $r\in\An(\ell)$, then, on $A_{jr}=\{Y_r>x_j\}$, we have $X_\ell\ge(\beta_{r\to\ell}/\beta_{r\to j})Y_r\ge(\beta_{r\to\ell}/\beta_{r\to j})x_j$, and hence
\[
    |f(\bm V)-f(\bm W^{(r)})|
    =
    V_\ell
    \le
    \bar F_\ell\Big(\frac{\beta_{r\to\ell}}{\beta_{r\to j}}x_j\Big)
    =
    O(h),
\]
by \eqref{eq:xjh-equivalence} and since $\bar F_\ell(x)=O(x^{-\alpha})$ by the union bound. As $\Prob(A_{jr})/h=O(1)$, replacing $\bm V$ by $\bm W^{(r)}$ in \eqref{eq:rep1} thus costs $O(h)$. Finally, by independence of $\bm W^{(r)}$ and $A_{jr}$ and the first assertion in \eqref{eq:linear-scm-bias-sbj},
\begin{align*}
    \frac1h
    \Exp\big[
        f(\bm V)\bm 1(B_j)
    \big]
    &=
    \sum_{r\in\An(j)}
    \frac{\Prob(A_{jr})}{h}\Exp[f(\bm W^{(r)})]
    +O(h^{\delta/\alpha})
    \\&=
    \sum_{r\in\An(j)}
    w_{jr}\Exp[f(\bm W^{(r)})]
    +O(h^{\delta/\alpha})
    =\Gamma_{j\ell}+O(h^{\delta/\alpha}),
\end{align*}
where we have used Proposition~\ref{prop:ancestral-mixture} and Lemma~\ref{lem:linear-scm-multivariate} in the last step.

Finally, Fubini's theorem gives
\begin{align*}
    \frac1h\int_0^1D_{j\ell}(h,v)\,\diff v
    &=
    \frac1h
    \Exp\big[
        (1-V_\ell)\bm 1(V_j\le h)
    \big]=
    \frac1h
    \Exp\big[
        f(\bm V)\bm 1(B_j)
    \big].
\end{align*}
Consequently,
\[
    \frac1h\int_0^1D_{j\ell}(h,v)\,\diff v
    =
    \Gamma_{j\ell}+O(h^{\delta/\alpha}).
\]
Setting $h=1/t$ proves
$
    B_{j\ell}(t)
    =
    O(t^{-\delta/\alpha}),
$
and the sufficient condition on $k$
follows immediately.
\end{proof}

\begin{proof}[Proof of Lemma~\ref{lem:max-linear-multivariate}]
As in the proof of Lemma~\ref{lem:linear-scm-multivariate}, let $I=\{j,s\}$ with $j,s\in[d]$, where possibly $j=s$, let $\mathcal C_I=\An(j)\cap\An(s)$ and $w_{Ir}=w_{jr}\wedge w_{sr}$. We need to show that $\nu_{I,h}\weak\nu_I$ with
\[
\nu_I=\sum_{r\in\mathcal C_I}w_{Ir}\,\mu_r,
\]
which reduces to the claimed formulas for $\nu_j$ and $\nu_{\{j,s\}}$ in the cases $j=s$ and $j\neq s$, respectively.

Let $\lambda(u)=-\log(1-u)$ for $u\in[0,1)$. Since the copula of $\bm X$ can be written as $C(\bm u)=\exp\{-L(-\log u_1,\dots,-\log u_d)\}$ with $L$ from \eqref{eq:max-linear-stdf}, we have, for every $M\subset[d]$ and $u_m\in(0,1)$, $m\in M$,
\begin{align}
\label{eq:max-linear-survival}
\Prob(V_m>u_m, \  \forall m\in M)
=
\exp\Big\{-\sum_{r\in[d]}\max_{m\in M}w_{mr}\lambda(u_m)\Big\},
\end{align}
with the convention $\max_\emptyset=0$. For $S\subset[d]$ and $\bm v=(v_m)_{m\in S}\in(0,1)^S$, let
\[
Q_{S,\bm v} = \{\bm z \in [0,1]^d: z_m > v_m \text{ for all } m \in S\},
\]
with $Q_{\emptyset,\bm v}=[0,1]^d$. We first show that, for all such $S$ and $\bm v$,
\begin{align}
\label{eq:max-linear-orthant-limit}
\lim_{h \downarrow 0} \nu_{I,h}(Q_{S,\bm v}) = \nu_I(Q_{S,\bm v}).
\end{align}
If $S\cap I\neq\emptyset$, both sides of \eqref{eq:max-linear-orthant-limit} vanish for all $h<\min_{m\in S\cap I}v_m$: indeed, $\nu_{I,h}$ is concentrated on $\{\bm z:z_j\le h,z_s\le h\}$, and $\nu_I$ is concentrated on $\{\bm z:z_j=z_s=0\}$, since $W^{(r)}_j=W^{(r)}_s=0$ for $r\in\mathcal C_I$. Hence, we may assume that $S\subset[d]\setminus I$ in what follows.
Let $\mathcal U_S=\{r\in[d]:S\cap\Des(r)=\emptyset\}$ denote the set of nodes that have no descendants in $S$. Further, write $\pi_S=\Prob(V_m>v_m, \ \forall m\in S)$, with $\pi_\emptyset=1$. Since $v_m>0$ for all $m\in S$ and $W_m^{(r)}=0$ for $m\in\Des(r)$, we have $\mu_r(Q_{S,\bm v})=\Prob(W_m^{(r)} > v_m ,\ \forall m \in S)=\bm 1(r\in\mathcal U_S)\,\pi_S$, and hence, since $w_{Ir}=0$ for $r\notin\mathcal C_I$,
\begin{align}
\label{eq:max-linear-nuI-Q}
\nu_I(Q_{S,\bm v})
=
\pi_S\sum_{r\in\mathcal U_S}w_{Ir}.
\end{align}
For $r\in[d]$, let $\Lambda_r=\max_{m\in S}w_{mr}\lambda(v_m)$. Since $v_m>0$ for $m\in S$, we have $\Lambda_r=0$ if and only if $w_{mr}=0$ for all $m\in S$, that is, if and only if $r\in\mathcal U_S$. Hence, by \eqref{eq:max-linear-survival} with $M=S$,
\[
\pi_S
=
\exp\Big(-\sum_{r\in[d]}\Lambda_r\Big)
=
\exp\Big(-\sum_{r\notin\mathcal U_S}\Lambda_r\Big).
\]
Moreover, if $r\notin\mathcal U_S$, then $\Lambda_r>0$, and hence $\max\{w_{jr}\lambda(h),w_{sr}\lambda(h)\}\le\Lambda_r$ for all sufficiently small $h$, since $\lambda(h)\to0$ as $h\downarrow0$. Write $\tau_j=\sum_{r\in\mathcal U_S}w_{jr}$, $\tau_s=\sum_{r\in\mathcal U_S}w_{sr}$ and $\tau_{\vee}=\sum_{r\in\mathcal U_S}(w_{jr}\vee w_{sr})$. Then, for all sufficiently small $h$, \eqref{eq:max-linear-survival} with $M=S\cup\{j\}$ and $M=S\cup\{j,s\}$ yields
\begin{align*}
\Prob(V_j>h,\, V_m>v_m, \ \forall m\in S)
&=
\exp\Big(-\sum_{r\notin\mathcal U_S}\Lambda_r-\sum_{r\in\mathcal U_S}w_{jr}\lambda(h)\Big)
=
\pi_S\, e^{-\lambda(h)\,\tau_j},
\\
\Prob(V_j>h,\, V_s>h,\, V_m>v_m, \  \forall m\in S)
&=
\exp\Big(-\sum_{r\notin\mathcal U_S}\Lambda_r-\sum_{r\in\mathcal U_S}(w_{jr}\vee w_{sr})\lambda(h)\Big)
=
\pi_S\, e^{-\lambda(h)\,\tau_{\vee}},
\end{align*}
and, likewise, $\Prob(V_s>h,\, V_m>v_m, \ \forall m\in S)=\pi_S\, e^{-\lambda(h)\,\tau_s}$. By inclusion--exclusion, we therefore obtain
\begin{align*}
h \nu_{I,h}(Q_{S,\bm v})
&=
\Prob\big(V_j\le h, V_s\le h, V_m > v_m , \ \forall m \in S\big)
\\&=
\pi_S
- \Prob(V_j>h,\, V_m>v_m, \ \forall m\in S)
\\&\hspace{1cm}
- \Prob(V_s>h,\, V_m>v_m , \ \forall m\in S)
+ \Prob(V_j>h,\, V_s>h,\, V_m>v_m, \ \forall m\in S)
\\&=
\pi_S
\Big\{1 - e^{-\lambda(h)\tau_j} - e^{-\lambda(h)\tau_s} + e^{-\lambda(h)\tau_{\vee}}\Big\}.
\end{align*}
Since $\lambda(h)/h\to1$, we have $(1-e^{-\lambda(h)\tau})/h\to\tau$ as $h\downarrow0$ for every $\tau\ge0$. Hence,
\[
\lim_{h \downarrow 0} \nu_{I,h}(Q_{S,\bm v})
=
\pi_S\,(\tau_j+\tau_s-\tau_{\vee})
=
\pi_S \sum_{r \in \mathcal U_S} w_{Ir},
\]
where we used $x+y-x\vee y=x\wedge y$. By \eqref{eq:max-linear-nuI-Q}, this proves \eqref{eq:max-linear-orthant-limit}.

It remains to show that \eqref{eq:max-linear-orthant-limit} implies $\nu_{I,h} \weak \nu_I$. Choosing $S=\emptyset$ in \eqref{eq:max-linear-orthant-limit} shows that the total masses converge. If $\nu_I([0,1]^d)=0$, weak convergence to the zero measure is then immediate. Otherwise, it suffices to show weak convergence of the normalized probability measures, for which we apply \citet[Theorem~2.2]{Bil99}. Let $\mathcal R$ denote the class of rectangles $R=\prod_{m \in [d]} J_m$, where each $J_m$ is of the form $[0,1]$, $[0,b]$, $(a,1]$ or $(a,b]$ with $0<a<b<1$. Since $\bm 1_{[0,b]} = 1 - \bm 1_{(b,1]}$ and $\bm 1_{(a,b]} = \bm 1_{(a,1]}-\bm 1_{(b,1]}$, the indicator of every $R \in \mathcal R$ is a finite signed linear combination of indicators of sets of the form $Q_{S,\bm v}$; hence \eqref{eq:max-linear-orthant-limit} implies $\nu_{I,h}(R) \to \nu_I(R)$ for all $R \in \mathcal R$. The class $\mathcal R\cup\{\emptyset\}$ is closed under finite intersections, and every relatively open subset of $[0,1]^d$ is a countable union of sets from $\mathcal R$ with rational endpoints. Hence, \citet[Theorem~2.2]{Bil99} yields $\nu_{I,h} \weak \nu_I$, which completes the proof.
\end{proof}

\begin{proof}[Proof of Lemma~\ref{lem:max-linear-condition-D}]
For $(j,\ell)\in\Pc_2$, the copula $C_{j\ell}$ is an extreme value copula with Pickands dependence function $A_{j\ell}$ from \eqref{eq:max-linear-pickands}. The assertion hence follows from Example~\ref{ex:evc-D}.
\end{proof}

\begin{proof}[Proof of Lemma~\ref{lem:max-linear-bias}]
Fix $(j,\ell)\in\Pc_2$, recall the Pickands dependence function $A_{j\ell}$ from \eqref{eq:max-linear-pickands}, and recall $p_{j\ell}=\sum_{r\in\An(j)\cap\An(\ell)}w_{jr}$ from Proposition~\ref{prop:ancestral-mixture}. We show that $A_{j\ell}$ is affine in a left neighbourhood of $1$, more precisely, that
\begin{align}
\label{eq:max-linear-pickands-affine}
A_{j\ell}(t) = 1 - p_{j\ell}(1-t), \qquad t \in [t_{j\ell}, 1], \qquad
t_{j\ell}
:=
\max_{r \in \An(\ell)} \frac{w_{jr}}{w_{jr} + w_{\ell r}} < 1.
\end{align}
Indeed, for $r \in \An(\ell)$, we have $w_{\ell r}>0$ and hence $\max\{w_{jr}(1-t), w_{\ell r}t\} = w_{\ell r}t$ for all $t \ge w_{jr}/(w_{jr} + w_{\ell r})$. For $r \notin \An(\ell)$, we have $w_{\ell r}=0$, such that $\max\{w_{jr}(1-t), w_{\ell r}t\} = w_{jr}(1-t)$ for all $t \in [0,1]$. Consequently, for $t \in [t_{j\ell},1]$,
\[
A_{j\ell}(t)
=
t \sum_{r \in \An(\ell)} w_{\ell r} + (1-t)\sum_{r \notin \An(\ell)} w_{jr}
=
t + (1-t)(1-p_{j\ell})
=
1 - p_{j\ell}(1-t),
\]
where we have used that $\sum_{r \in \An(\ell)} w_{\ell r}=\sum_{r \in [d]} w_{\ell r}=1$ and
\[
\sum_{r \notin \An(\ell)} w_{jr} = 1 - \sum_{r \in \An(\ell)} w_{jr} = 1 - \sum_{r \in \An(j) \cap \An(\ell)} w_{jr} = 1-p_{j\ell},
\]
which in turn follows from $w_{jr}=0$ for $r \notin \An(j)$. This proves \eqref{eq:max-linear-pickands-affine}, and in particular $A_{j\ell}'(1-)=p_{j\ell}$. The assertion then follows from Example~\ref{ex:evc-bias}, since $t_{j\ell}<1$.
\end{proof}

\subsection{Proofs for Section~\ref{sec:graph-recovery}}

\begin{proof}[Proof of Proposition~\ref{prop:graph-recovery-studentized}]
Recall $\Gb_{n,j\ell}$ from \eqref{eq:weak-CTC-multivariate},
and let $\bm G=(G_{j\ell})_{(j,\ell)\in\Pc_2}$ be a centered Gaussian vector with covariance matrix $\Sigma$ from \eqref{eq:covariance-identity}.
By Theorem~\ref{theo:CTC-multivariate} and the assumed bias condition,
\[
    \Big(
    \sqrt{k}
    (\widehat\Gamma_{j\ell}-\Gamma_{j\ell})
    \Big)_{(j,\ell)\in\Pc_2}
    \weak \bm G.
\]
Moreover, by Theorem~\ref{theo:bootstrap-multivariate}, we have $\widehat\sigma_{j\ell}=\sigma_{j\ell}+o_\Prob(1)$ for all $(j,\ell)\in\Pc_2$, and hence
\[
    \widehat\delta_n
    =
    \delta + o_\Prob(1), \qquad \delta:=
    c\max_{(j,\ell)\in\Pc_2}\sigma_{j\ell}>0.
\]
Consequently,
\[
    \widehat s_{j\ell}
    =
    s_{j\ell} + o_\Prob(1), \qquad s_{j\ell}:=\max\{\sigma_{j\ell},\delta\}>0
\]
for all $(j,\ell)\in\Pc_2$.
By Slutsky's lemma, we obtain
\[
    \max_{(j,\ell)\in\Pc_2}
    \frac{\sqrt{k}
    |\widehat\Gamma_{j\ell}-\Gamma_{j\ell}|}
    {\widehat s_{j\ell}}
    \weak
    \max_{(j,\ell)\in\Pc_2}
    \frac{|G_{j\ell}|}{s_{j\ell}}.
\]
Similarly, Theorem~\ref{theo:bootstrap-multivariate} and the conditional continuous mapping theorem imply
\[
    \max_{(j,\ell)\in\Pc_2}
    \frac{|\widehat\Gb_{n,j\ell}|}
    {\widehat s_{j\ell}}
    \cweak{\xi}
    \max_{(j,\ell)\in\Pc_2}
    \frac{|G_{j\ell}|}{s_{j\ell}}.
\]
Since at least one marginal variance is strictly positive, the CDF of
$\max_{(j,\ell)\in\Pc_2}|G_{j\ell}|/s_{j\ell}$ is continuous and strictly increasing on $(0,\infty)$. Consequently, $\widehat q_{1-\alpha}$ converges
in probability to its $(1-\alpha)$-quantile, denoted by $q_{1-\alpha}$.
It follows that
\begin{align}
\label{eq:simultaneous-coverage-studentized}
    \lim_{n \to \infty} \Prob\Big(
    \max_{(j,\ell)\in\Pc_2}
    \frac{\sqrt{k}
    |\widehat\Gamma_{j\ell}-\Gamma_{j\ell}|}
    {\widehat s_{j\ell}}
    \le\widehat q_{1-\alpha}
    \Big)
    = 1-\alpha.
\end{align}
On the event inside the probability, we have
\[
    \widehat\Gamma_{j\ell}
    -\frac{\widehat q_{1-\alpha}
    \widehat s_{j\ell}}{\sqrt{k}}
    \le \Gamma_{j\ell}
    \le
    \widehat\Gamma_{j\ell}
    +\frac{\widehat q_{1-\alpha}
    \widehat s_{j\ell}}{\sqrt{k}}
\]
for all $(j,\ell)\in\Pc_2$.
By the definitions of the confidence graphs, this implies
$
    \widehat E_{\gamma,\alpha}^L
    \subset E_\gamma
    \subset\widehat E_{\gamma,\alpha}^U.
$
Hence, \eqref{eq:graphrecovery-studentized} follows from
\eqref{eq:simultaneous-coverage-studentized}.

Next, for every $(j,\ell)\in E_\gamma$, we have $\Gamma_{j\ell}>\gamma$.
Since
\[
    \widehat\Gamma_{j\ell}
    =\Gamma_{j\ell}+o_\Prob(1),
    \qquad
    \frac{\widehat q_{1-\alpha}
    \widehat s_{j\ell}}{\sqrt{k}}
    =o_\Prob(1),
\]
it follows that
\[
    \lim_{n\to\infty}\Prob\big(
    (j,\ell)\in\widehat E_{\gamma,\alpha}^L
    \big)=1.
\]
Since $d$ is fixed, this implies
\[
    \lim_{n\to\infty}\Prob\big(
    E_\gamma\subset\widehat E_{\gamma,\alpha}^L
    \big)=1.
\]
Together with \eqref{eq:graphrecovery-studentized}, this yields
\[
    \liminf_{n\to\infty}
    \Prob\big(
    \widehat E_{\gamma,\alpha}^L=E_\gamma
    \big)\ge1-\alpha.
\]

Finally, if $\Gamma_{j\ell}\ne\gamma$ for all
$(j,\ell)\in\Pc_2$, then
\[
    \widehat\Gamma_{j\ell}
    \pm\frac{\widehat q_{1-\alpha}
    \widehat s_{j\ell}}{\sqrt{k}}
    =\Gamma_{j\ell}+o_\Prob(1)
\]
for every $(j,\ell)\in\Pc_2$. Since $d$ is fixed, both confidence graphs
then recover the true edge set with probability tending to one, which proves the final assertion.
\end{proof}

\subsection{Proofs for Section~\ref{sec:general-tail-quantities}}

\begin{proof}[Proof of Theorem~\ref{theo:unknown-margins-bivariate}]
Throughout, we write $t_n=k/n$.
Define a version of $\widehat G_n$ where $\widehat V_{i2}$ has been replaced by the unobservable $V_{i2}=1-F_2(X_{i2})$, that is, 
\[
\gghat^{\circ}(x) = \frac1k \sum_{i=1}^n g(V_{i2}) \bm 1\big( \widehat V_{i1} \le t_nx\big).
\]
We will show below that
\begin{align} \label{eq:hat-circ-general}
\sup_{x \in [0,T]} \Big| \ggprocess(x) 
- 
\sqrt{k}\big\{ \gghat^\circ(x) - \ggmean(x) \big\}  \Big| = o_\Prob(1).
\end{align}
Next, define a version where also $\widehat V_{i 1}$ has been replaced by $ V_{i 1}$:
\begin{align*}
    \gghat^{\star}(x) &= \frac1k \sum_{i=1}^n g(V_{i2}) \bm 1\big(  V_{i1} \le t_nx\big)
\end{align*}
Write $\hat x_{n} := t_n^{-1} G_{n1}^{-1}(t_nx)$ where $G_{n1}(u) = n^{-1}\sum_{i=1}^n \bm 1(V_{i1} \le u)$. Here and below, we suppress the dependence of
$\hat x_{n}$ on $x$. 
By the definitions of $\widehat V_{i1}$ and the generalized inverse, the numbers of summands in $\gghat^\circ(x)$ and $\gghat^\star(\hat x_n)$ are given by
$\lfloor kx(n+1)/n\rfloor$ and $\lceil kx\rceil$, respectively, for $x>0$.
For sufficiently large $n$, these differ by at most one, uniformly over $x\in[0,T]$.
Consequently,
\begin{align} \label{eq:hat-circ-general2}
    \sup_{x\in[0,T]}
    \big|
        \sqrt k \{ \gghat^\circ(x) - \ggmean(x) \} - \sqrt k \{\gghat^\star(\hat x_n) - \ggmean(x)\}
    \big|
    \le
    \frac{\|g\|_\infty}{\sqrt k} = o_\Prob(1).
\end{align}
Next, we may decompose
\begin{align} \label{eq:decomp1-general}
\sqrt k \{\gghat^\star(\hat x_n) - \ggmean(x)\}
=
\Bb_{n1}(x) + \Bb_{n2}(x), 
\end{align}
where
\[
\Bb_{n1}(x) = \sqrt k \big\{  \gghat^\star(\hat x_{n}) - \ggmean(\hat x_{n}) \big\} , \qquad
\Bb_{n2}(x)  =  \sqrt k \big\{ \ggmean(\hat x_{n}) - \ggmean(x)\big\}.
\]
Also, let 
\begin{align*}
    \Bb^\star_{n1}(x) &:= \frac1{\sqrt n} \sum_{i=1}^n \big\{ \psi_{n,x}(V_{i1},V_{i2}) - \Exp[\psi_{n,x}(V_{i1},V_{i2})] \big\},
    \\
    \Bb^\star_{n2}(x) &:= \frac1{\sqrt n} \sum_{i=1}^n \big\{ \zeta_{n,x}(V_{i1}) - \Exp[\zeta_{n,x}(V_{i1})] \big\},
\end{align*}
where
\[
\psi_{n,x}(v_1,v_2) = t_n^{-1/2} g(v_2) \bm 1\big(v_{1} \le  t_nx \big),
\qquad 
\zeta_{n,x}(v_{1}) = - t_n^{-1/2} \, \ggcmean(t_nx) \, \bm 1\big(v_{1} \le  t_nx \big).
\]
We will show below that
\begin{align} \label{eq:decomp1-term1-general}
    \sup_{x \in [0,T]} \Big| \Bb_{n1}(x) 
    -
    \Bb^\star_{n1}(x)
    \Big| =  o_\Prob(1), 
\end{align}
and that 
\begin{align} \label{eq:decomp1-term2-general}
    \sup_{x \in [0,T]}\Big| \Bb_{n2} (x) 
    - 
    \Bb^\star_{n2}(x) \Big| = o_\Prob(1).
\end{align}
The assertion in \eqref{eq:linear-bivariate-process-g} then follows from \eqref{eq:hat-circ-general}, \eqref{eq:hat-circ-general2}, \eqref{eq:decomp1-general}, \eqref{eq:decomp1-term1-general} and \eqref{eq:decomp1-term2-general}, observing that
\[
\phi_{n,x}(v_1,v_2)
=
\psi_{n,x}(v_1,v_2) + \zeta_{n,x}(v_1)
\]
and that $\Exp[\phi_{n,x}(V_{1},V_{2})]=0$, since $\Exp[g(V_2)\bm 1(V_1\le t_nx)]=t_nx\,\ggcmean(t_nx)$.

It remains to show \eqref{eq:hat-circ-general},
\eqref{eq:decomp1-term1-general} and
\eqref{eq:decomp1-term2-general}. For later use, we start by showing that
\begin{align}
\label{eq:quantile-expansion-general}
\sup_{x\in[0,T]}
\left|
\sqrt{k}(\hat x_n-x)+\alpha_{n1}(x)
\right|
=
o_{\Prob}(1),
\end{align}
where $\alpha_{n1}$ is the tail empirical process defined by
\begin{align}
\label{eq:tail-empirical-process}
\alpha_{n1}(x)
=
\frac{1}{\sqrt{k}}
\sum_{i=1}^n
\left\{
\bm 1(V_{i1}\le t_nx)-t_nx
\right\},
\qquad x\ge0.
\end{align}
This is a standard consequence of Vervaat's inversion lemma
\citep[Lemma~1]{Vervaat1972}; for completeness, we provide a short
direct argument.
Proposition~\ref{prop:known-margins-bivariate-2param}, applied with \(g\equiv1\) and $a=0$,
implies that
\[
\big(\alpha_{n1}(x)\big)_{x\in[0,2T]}
\rightsquigarrow
\big(W(x)\big)_{x\in[0,2T]}
\qquad\text{in }\ell^\infty([0,2T]),
\]
where \(W\) is a standard Brownian motion. In particular,
\(\alpha_{n1}\) is asymptotically tight and asymptotically uniformly
equicontinuous.
Moreover, with
\[
\widetilde G_{n1}(x)
=
t_n^{-1}G_{n1}(t_nx)
=
x+\frac{1}{\sqrt{k}}\alpha_{n1}(x),
\]
we can write $\hat x_n=\widetilde G_{n1}^{-1}(x)$. Since the \(V_{i1}\) have a continuous distribution, \(G_{n1}\) has, almost surely, jumps of size \(1/n\), and hence \(\widetilde G_{n1}\) has jumps of size \(1/k\). By the defining
property of the generalized inverse, for every
\(x\in(0,t_n^{-1}]\),
\[
\widetilde G_{n1}(\hat x_n-)
<
x
\le
\widetilde G_{n1}(\hat x_n).
\]
It follows that
\[
0
\le
\widetilde G_{n1}(\hat x_n)-x
\le
\frac1k.
\]
Consequently,
\begin{align}
\label{eq:bahadur}
\left|
\sqrt{k}(\hat x_n-x)+\alpha_{n1}(\hat x_n)
\right|
&=
\sqrt{k}
\left|
\hat x_n-x+\widetilde G_{n1}(\hat x_n)-\hat x_n
\right|
\nonumber\\
&=
\sqrt{k}
\left|
\widetilde G_{n1}(\hat x_n)-x
\right|
\le
\frac1{\sqrt{k}}.
\end{align}
The same inequality holds for \(x=0\), since then \(\hat x_n=0\) and \(\alpha_{n1}(0)=0\). We next verify that
\begin{align}
\label{eq:hatx-cons}
\sup_{x\in[0,T]}|\hat x_n-x|
=
O_{\Prob}(k^{-1/2})
=
o_{\Prob}(1).
\end{align}
To this end, note that
\[
\widetilde G_{n1}(2T)
=
2T+\frac{1}{\sqrt{k}}\alpha_{n1}(2T) = 2T + o_\Prob(1),
\]
and hence
$\Prob\bigl(\widetilde G_{n1}(2T)>T\bigr)\to1$.
On the event \(\{\widetilde G_{n1}(2T)>T\}\), monotonicity of the
generalized inverse implies
$
\hat x_n
=
\widetilde G_{n1}^{-1}(x)
\le 2T$ 
for all $x\in[0,T]$.
Therefore,
\begin{align}
\label{eq:2tt}
\Prob\Big(
\sup_{x\in[0,T]}\hat x_n\le2T
\Big)\to1.
\end{align}
On this event, \eqref{eq:bahadur} yields
\begin{align*}
\sqrt{k}\sup_{x\in[0,T]}|\hat x_n-x|
&\le
\sup_{x\in[0,T]}|\alpha_{n1}(\hat x_n)|
+\frac1{\sqrt{k}}
\\
&\le
\sup_{y\in[0,2T]}|\alpha_{n1}(y)|
+\frac1{\sqrt{k}}
=
O_{\Prob}(1),
\end{align*}
where the final equality follows from asymptotic tightness of
\(\alpha_{n1}\) in \(\ell^\infty([0,2T])\). Since the event under
consideration has probability tending to one, this proves
\eqref{eq:hatx-cons}.

Finally, asymptotic uniform equicontinuity of \(\alpha_{n1}\), together
with \eqref{eq:hatx-cons}, implies that
\[
\sup_{x\in[0,T]}
\left|
\alpha_{n1}(\hat x_n)-\alpha_{n1}(x)
\right|
=
o_{\Prob}(1).
\]
Indeed, for every fixed \(\delta>0\), on the event
\[
\sup_{x\in[0,T]}|\hat x_n-x|<\delta
\quad\text{and}\quad
\sup_{x\in[0,T]}\hat x_n\le2T,
\]
we have
\[
\sup_{x\in[0,T]}
\left|
\alpha_{n1}(\hat x_n)-\alpha_{n1}(x)
\right|
\le 
\sup_{\substack{y,z\in[0,2T]\\|y-z|<\delta}}
|\alpha_{n1}(y)-\alpha_{n1}(z)|,
\]
and the claim follows by first letting \(n\to\infty\) and then
\(\delta\downarrow0\). Combining this with \eqref{eq:bahadur} proves
\eqref{eq:quantile-expansion-general}.

\medskip
\noindent
Proof of \eqref{eq:hat-circ-general}: let $L=L_g$ denote the Lipschitz constant of $g$. It follows that, if $n$ is sufficiently large such that $kT <n$,
\begin{align*}
\sup_{x \in [0,T]}|\gghat(x) - \gghat^\circ(x)|
&\le
\frac{1}k \sum_{i=1}^n | g(\widehat V_{i2}) - g(V_{i2}) | \bm 1\Big( \widehat V_{i1} \le \frac{kT}{n} \Big)
\\&\le 
\frac{L}k \sum_{i=1}^n \sup_{v \in [0,1]} \Big( | G_{n2}(v) - v| + \frac1{n+1}\Big)\bm 1\Big( \widehat V_{i1} \le \frac{kT}{n} \Big)
\\&=
O_\Prob(n^{-1/2}) \cdot \frac1k \sum_{i=1}^n  \bm 1\Big( \widehat V_{i1} \le \frac{kT}{n} \Big)
= O_\Prob(n^{-1/2})=o_\Prob(k^{-1/2}),
\end{align*}
where we have used that
$
    \sup_{v\in[0,1]}|G_{n2}(v)-v|
    =O_\Prob(n^{-1/2}),
$
with $G_{n2}(v)=n^{-1}\sum_{i=1}^n
\bm 1(V_{i2}\le v)$, together with
$
    \max_{i\in[n]}|\widehat V_{i2}-V_{i2}|
    \le
    \sup_{v\in[0,1]}|G_{n2}(v)-v|
    +1/(n+1)
$
and
$
    \sum_{i=1}^n
    \bm 1(\widehat V_{i1}\le kT/n)
    =
    \lfloor(n+1)kT/n\rfloor
    =O(k)
$
with probability one.

\medskip
\noindent
Proof of \eqref{eq:decomp1-term1-general}: 
It follows from the definition of $\Bb^\star_{n1}$ and Proposition~\ref{prop:known-margins-bivariate-2param} that $\Bb^\star_{n1}$, considered as a map from the underlying probability space to $\ell^\infty([0,2T])$, is asymptotically tight \citep[Theorem 1.5.4]{VanWel23}. Since the limit process in Proposition~\ref{prop:known-margins-bivariate-2param} has continuous sample paths almost surely, it follows from \citet[Addendum 1.5.8]{VanWel23} that, for every $\eps>0$,
\[
\lim_{\delta\to 0}\limsup_{n\to\infty} \Prob\Big( w_{\Bb^\star_{n1}}(\delta) >\eps\Big) = 0, \qquad 
w_{\Bb^\star_{n1}}(\delta)  := \sup_{\substack{x,y\in[0,2T]:\\|x-y|<\delta}} \big| \Bb^\star_{n1}(x) - \Bb^\star_{n1}(y) |.
\]
Now, for arbitrary $\eps>0$ and $\delta \in (0,T)$, 
\begin{align*}
    \Prob\Big( \sup_{x \in [0,T]} \big| \Bb^\star_{n1}(\hat x_n) - \Bb^\star_{n1}(x) | >\eps\Big) 
    &\le 
    \Prob\Big( \sup_{x \in [0,T]}\big| \Bb^\star_{n1}(\hat x_n) - \Bb^\star_{n1}(x) | >\eps, \sup_{x \in [0,T]}|\hat x_n - x| < \delta \Big) 
    \\& \hspace{5cm} + \Prob\Big( \sup_{x \in [0,T]}|\hat x_n - x| \ge \delta  \Big)
    \\&
    \le 
    \Prob\Big(w_{\Bb^\star_{n1}}(\delta)  >\eps\Big)
    +  \Prob\Big( \sup_{x \in [0,T]} |\hat x_n - x|  \ge \delta  \Big).
\end{align*}
Since $\sup_{x \in [0,T]} | \hat x_n - x | = o_\Prob(1)$ as shown in \eqref{eq:hatx-cons}, for every fixed $\delta \in (0,T)$,
\[
\lim_{n\to\infty}
\Prob\big(\sup_{x \in [0,T]}|\hat x_n-x|\ge\delta\big)
=
0.
\]
Consequently,
\[
\limsup_{n\to\infty}
\Prob\Big(\sup_{x \in [0,T]}
\big|\Bb^\star_{n1}(\hat x_n)-\Bb^\star_{n1}(x)\big|>\eps
\Big)
\le
\limsup_{n\to\infty}
\Prob\big(w_{\Bb^\star_{n1}}(\delta)>\eps\big),
\]
which can be made arbitrary small by decreasing $\delta$. This completes the proof of \eqref{eq:decomp1-term1-general}.

\medskip
\noindent
Proof of \eqref{eq:decomp1-term2-general}:  throughout, recalling $t_n=k/n \to 0$, let
\[
\Bb^h_{n2}(x) = \sqrt k (\hat x_n - x) \ggcmean(t_n x).
\]
The proof of \eqref{eq:decomp1-term2-general} is finished once we show that
\begin{align} \label{eq:tbs2}
    \sup_{x \in [0,T]}\Big| \Bb_{n2} (x) - \Bb^h_{n2}(x) \Big| = o_\Prob(1),
    \qquad
    \sup_{x \in [0,T]}\Big| \Bb^h_{n2} (x) - \Bb^\star_{n2}(x) \Big| = o_\Prob(1)
\end{align}

Since $g$ is Lipschitz, it is absolutely continuous, which allows to define, for $z \in [0,1]$,
\[
\Psi(z) := \int_0^1 g'(v) D(z,v) \,\diff v.
\]
Integration by parts yields
$
\Exp[
g(V_2)\bm 1(V_1\le z)
]
=
zg(1)-\Psi(z),
$
which in turn implies that $\ggmean$ and $\ggcmean$ from \eqref{eq:gn-g}  and \eqref{eq:gn-Ag} satisfy
\begin{align}
\label{eq:gmean-A}
\ggmean(x)
=
xg(1)-t_n^{-1}\Psi(t_nx),
\qquad
x\in[0,t_n^{-1}],
\end{align}
and
\begin{align}
\label{eq:gcmean-A}
\ggcmean(z)
=
g(1)-\frac{\Psi(z)}{z},
\qquad z>0.
\end{align}

Since $D$ is a bivariate copula, for every $v \in (0,1)$, the partial derivative $\partial_1D(z,v)$ exists for Lebesgue-almost every $z \in (0,1)$. Thus, by Fubini's theorem, for almost every $z\in(0,1)$, $\partial_1D(z,v)$ exists for almost every $v\in(0,1)$.
Consequently, for almost every $z$, by dominated convergence 
\[
\frac{\Psi(z+h) - \Psi(z)}h
= 
\int_0^1 \frac{D(z+h,v) - D(z,v)}h  g'(v) \, \diff v
\longrightarrow \int_0^1 \partial_1 D(z,v) g'(v) \, \diff v =: Q(z)
\]
as $h\to 0$.
Hence, $\Psi'(z) = Q(z)$ Lebesgue-almost everywhere, say for all $z\in\mathcal Z$, where $\mathcal Z\subset(0,1)$ has full Lebesgue measure.
We extend the definition of $Q$ to all $z\in[0,1]$ by setting $Q(z) = \int_0^1  q_z(v) g'(v) \, \diff v$, where  $q_z(v) = \partial_1^+ D(z,v)$ whenever the latter exists and $q_z(v)=0$ else. By construction, the new definition coincides with the old one for $z\in\mathcal Z$; in particular, the preceding equality continues to hold. Moreover, Condition~\ref{cond:continuity2-bivariate} implies that, for almost every $v\in(0,1)$, $q_z(v) \to \partial_1^+ D(0,v)$ for $z \downarrow 0$, such that $Q(z) \to Q(0) = \int_0^1 \partial_1^+ D(0,v)g'(v)\, \diff v$ for $z \downarrow 0$ by dominated convergence. Similarly, dominated convergence yields $\Psi(z) / z \to Q(0)$ for $z \downarrow 0$. 
With the convention $\Psi(0) / 0 := Q(0)$, we therefore obtain
\begin{align}
\label{eq:uniform-Q-A}
\Delta_n := \sup_{s,x \in [0,2T]} \Big| Q(t_ns)  - \frac{\Psi(t_nx)}{t_nx} \Big| = o(1), \qquad n \to \infty.
\end{align}

For $x,y\in[0,2T]$, Equations
\eqref{eq:gmean-A} and \eqref{eq:gcmean-A} give
\begin{align*}
\ggmean(y)-\ggmean(x)
-
(y-x)\ggcmean(t_nx)
&= t_n^{-1} \big\{ \Psi(t_nx) - \Psi(t_ny) \big\} + (y-x) \frac{\Psi(t_nx)}{t_nx}.
\end{align*}
Since $\Psi$ is absolutely continuous (by the Lipschitz property of $D$), the fundamental theorem of calculus allows to write
\begin{align*} 
t_n^{-1} \big\{ \Psi(t_nx) - \Psi(t_ny) \big\} 
&=
- t_n^{-1} \int_{t_nx}^{t_ny} Q(z)  \,\diff z
=
- \int_{x}^{y} Q(t_ns)  \,\diff s,
\end{align*}
such that 
\begin{align*}
\ggmean(y)-\ggmean(x) - (y-x)\ggcmean(t_nx)
&=
-
\int_{x}^{y} \left\{Q(t_ns)  - 
\frac{\Psi(t_nx)}{t_nx}\right\} \,\diff s.
\end{align*}
Hence, for all $x,y\in[0,2T]$,
\[
\left|
\ggmean(y)-\ggmean(x)
-
(y-x)\ggcmean(t_nx)
\right|
\le
|y-x|\Delta_n.
\]

We now prove the  left convergence in \eqref{eq:tbs2}. By definition, we have
\[
\Bb_{n2} (x) - \Bb^h_{n2}(x)
= 
\sqrt k \big\{\ggmean( \hat x_n) - \ggmean(x) - (\hat x_n - x) \ggcmean(t_n x) \big\}
\]
On the event
$
\sup_{x\in[0,T]}\hat x_n\le2T,
$
whose probability tends to one as shown in \eqref{eq:2tt}, we obtain that 
\begin{align*}
\sup_{x\in[0,T]}
\left|
\Bb_{n2}(x)-\Bb^h_{n2}(x)
\right|
&\le
\Delta_n
\sup_{x\in[0,T]}
\left|
\sqrt{k}(\hat x_n-x)
\right| =
o_{\Prob}(1),
\end{align*}
because \(\Delta_n\to0\) by \eqref{eq:uniform-Q-A} and $
\sqrt k \sup_{x\in[0,T]}|\hat x_n-x|
=
O_{\Prob}(1)$
by \eqref{eq:hatx-cons}.

Finally, we prove the right convergence in \eqref{eq:tbs2}. By definition of \(\zeta_{n,x}\), we can write
\begin{align*}
\Bb^\star_{n2}(x)
=
-\frac{\ggcmean(t_nx)}{\sqrt{k}}
\sum_{i=1}^n
\big\{
\bm 1(V_{i1}\le t_nx)-t_nx
\big\}
=
-\ggcmean(t_nx)\alpha_{n1}(x).
\end{align*}
Consequently,
\begin{align*}
\sup_{x\in[0,T]}
\left|
\Bb^h_{n2}(x)-\Bb^\star_{n2}(x)
\right|
&\le
\sup_{x\in[0,T]}
|\ggcmean(t_nx)|
\sup_{x\in[0,T]}
\left|
\sqrt{k}(\hat x_n-x)+\alpha_{n1}(x)
\right|
\\
&\le
\|g\|_\infty
\sup_{x\in[0,T]}
\left|
\sqrt{k}(\hat x_n-x)+\alpha_{n1}(x)
\right|
=
o_{\Prob}(1),
\end{align*}
where we used \(|\ggcmean(z)|\le\|g\|_\infty\) and \eqref{eq:quantile-expansion-general}. This completes the proof of \eqref{eq:tbs2}.
\end{proof}

\begin{proof}[Proof of Corollary~\ref{cor:unknown-margins-bivariate}] 
Recall $t_n=k/n$. Let $\Hb_n^\star(x,a)$ denote the two-parameter process from
Proposition~\ref{prop:known-margins-bivariate-2param}. Since $\Exp[\phi_{n,x}(V_{i1},V_{i2})]=0$ for all $x\ge 0$, 
we have
\[
    \frac1{\sqrt n}
    \sum_{i=1}^n
    \phi_{n,x}(V_{i1},V_{i2})
    =
    \Hb_n^\star\bigl(x,\ggcmean(t_nx)\bigr),
    \qquad x\in[0,T].
\]
Recall that Condition~\ref{cond:continuity2-bivariate} implies Condition~\ref{cond:conditional-weak-limit-bivariate}. Hence, by Proposition~\ref{prop:known-margins-bivariate-2param}, with $\abound = \|g\|_\infty+1$, we have $\Hb_n^\star
    \weak
    \Hb^\star$ in $\ell^\infty\big(
        [0,T]\times[-\abound,\abound]
    \big)$, where $\Hb^\star$ has almost surely continuous sample paths.
Consequently, by \citet[Addendum~1.5.8]{VanWel23}, $\Hb_n^\star$ is asymptotically uniformly
equicontinuous in probability on
$[0,T]\times[-\abound,\abound]$. Since
\[
    \sup_{x\in[0,T]}
    \big|
        \ggcmean(t_nx)-\mu_g
    \big|
    =o(1), \qquad \mu_g := \Exp_{Z\sim H}[g(Z)],
\]
by Condition~\ref{cond:conditional-weak-limit-bivariate},
it follows that
\[
    \sup_{x\in[0,T]}
    \Big|
        \Hb_n^\star\bigl(x,\ggcmean(t_nx)\bigr)
        -
        \Hb_n^\star(x,\mu_g)
    \Big|
    =
    o_\Prob(1).
\]
Moreover, by the continuous mapping theorem,
\[
    \Big(
        \Hb_n^\star(x,\mu_g)
    \Big)_{x\in[0,T]}
    \weak
    \Big(
        \Hb^\star(x,\mu_g)
    \Big)_{x\in[0,T]}
\]
in $\ell^\infty([0,T])$. The latter is a tight, centered Gaussian
process with covariance function 
\[
\Cov\big[
        \Hb^\star(x,\mu_g),
        \Hb^\star(y,\mu_g)
    \big] = (x \wedge y) \Exp_{Z\sim H}[ (g(Z) - \mu_g)^2] = (x\wedge y)
    \Var_{Z\sim H}[g(Z)].
    \]
The assertion now follows from
Theorem~\ref{theo:unknown-margins-bivariate} and Slutsky's lemma.
\end{proof}

\begin{proof}[Proof of Theorem~\ref{theo:bootstrap}] 
Fix $\abound>2\|g\|_\infty$ and write $\mathcal T=[0,T]\times[-\abound,\abound]$. Let $t_n = k/n$.
For $(x,a)\in\mathcal T$, define
\begin{align*}
    f_{ni}(x,a)
    &=
    \frac1{\sqrt{k}}
    \{g(V_{i2})-a\}
    \bm 1(V_{i1}\le t_nx),
    \\
    \widehat f_{ni}(x,a)
    &=
    \frac1{\sqrt{k}}
    \{g(\widehat V_{i2})-a\}
    \bm 1(\widehat V_{i1}\le t_nx),
\end{align*}
and let
\[
\widehat\Bb_n(x,a):=
        \sum_{i=1}^n
        \xi_i\, \widehat f_{ni}(x,a).
\]
Since 
$
\widehat \Hb_n(x) = \widehat\Bb_n(x, \widehat a_n(x))
$
where $\widehat a_n(x) = \gghat(x) / x$ if $x>0$ and $\widehat a_n(0)=0$, the proof is finished once we show that
\begin{equation}
\label{eq:two-param-bootstrap}
    \Big(
        \widehat\Bb_n(x,a)
    \Big)_{(x,a)\in\mathcal T}
    \cweak{\xi}
    \Big(
        \Hb^\star(x,a)
   \Big)_{(x,a)\in\mathcal T} \qquad \text{ in } \quad \ell^\infty(\mathcal T),
\end{equation}
and that 
\begin{align}
\label{eq:hb-bb-approx}
    \sup_{x\in[0,T]}
    \left|
        \widehat\Bb_n(x,\widehat a_n(x))
        -
        \widehat\Bb_n(x, \mu_g)
    \right|
    =
    o_{\Prob_\xi}(1), \qquad \mu_g:=\Exp_{Z\sim H}[g(Z)],
\end{align}
in probability. Indeed, note that
\[
    \Cov\big[
        \Hb^\star(x,\mu_g),
        \Hb^\star(y,\mu_g)
    \big]
    =
    (x\wedge y)\Var_{Z\sim H}[g(Z)],
\]
which is the covariance function of the Gaussian limit process $\gglimitprocess$ from Corollary~\ref{cor:unknown-margins-bivariate}.

For the proof of \eqref{eq:two-param-bootstrap}, we apply Theorem~11.18 in \citet{Kosorok2008} with 
\[
    \mu_{ni}(x,a)=\Exp[f_{ni}(x,a)]
    \qquad\text{and}\qquad
    \widehat\mu_{ni}(x,a)
    =
    f_{ni}(x,a)-\widehat f_{ni}(x,a);
\]
note that $\widehat f_{ni}(x,a) = f_{ni}(x,a)-\widehat\mu_{ni}(x,a)$. 
We need to verify the conditions of Theorem~11.18. 

First, the assumptions on the two-parameter triangular array $(f_{ni})$ follow from other proofs. Indeed,  note 
that the proof of Proposition~\ref{prop:known-margins-bivariate-2param} is based on verifying the conditions of Theorem~11.20 in \citet{Kosorok2008}. In turn, Theorem~11.20 is proved by showing that these conditions imply those of Theorem~11.16. Hence, the two-parameter triangular array $(f_{ni})$ satisfies the conditions of Theorem~11.16 in \citet{Kosorok2008}. 

It remains to verify Conditions~(G)--(I) of Theorem~11.18.
For that purpose, we start by showing that
\begin{equation}
\label{eq:two-param-L2}
    \sup_{(x,a)\in\mathcal T}
    \sum_{i=1}^n
    \big\{
        \widehat f_{ni}(x,a)-f_{ni}(x,a)
    \big\}^2
    =
    o_\Prob(1).
\end{equation}
Indeed,
\begin{align*}
    \sqrt{k}
    \big\{
        \widehat f_{ni}(x,a)-f_{ni}(x,a)
    \big\}
    ={}&
    \big\{
        g(\widehat V_{i2})-g(V_{i2})
    \big\}
    \bm 1(\widehat V_{i1}\le t_nx)
    \\
    &+
    \{g(V_{i2})-a\}
    \big\{
        \bm 1(\widehat V_{i1}\le t_nx)
        -
        \bm 1(V_{i1}\le t_nx)
    \big\}.
\end{align*}
We consider both terms on the right separately. First, since $g$ is Lipschitz-continuous and the empirical marginal
CDF is uniformly consistent, we have
$
    \max_{1\le i\le n}
    \big|
        g(\widehat V_{i2})-g(V_{i2})
    \big|
    =o_\Prob(1).
$
Moreover, $\sum_{i=1}^n
\bm 1(\widehat V_{i1}\le t_nx)
\le (n+1)kT/n$
for all $x \in [0,T]$,
and hence
\[
    \sup_{x\in[0,T]}
    \frac1k
    \sum_{i=1}^n
    \big|
        g(\widehat V_{i2})-g(V_{i2})
    \big|^2
    \bm 1\big(\widehat V_{i1}\le t_nx\big)
    =o_\Prob(1).
\]
Second,  boundedness of $g$  gives
\begin{align*}
    &\sup_{(x,a)\in \mathcal T}
    \frac1k
    \sum_{i=1}^n
    \big|
        g(V_{i2})-a
    \big|^2
    \big|
        \bm 1\big(\widehat V_{i1}\le t_nx\big) - \bm 1\big(V_{i1}\le t_nx\big)
    \big|
    \\ &\hspace{4cm}\le 
    (\| g\|_\infty+\abound)^2 
    \sup_{x\in[0,T]}
    \frac1k \sum_{i=1}^n
    \big|
        \bm 1\big(\widehat V_{i1}\le t_nx\big) - \bm 1\big(V_{i1}\le t_nx\big)
    \big|.
\end{align*}
Since the pseudo-observations $\widehat V_{i1}$
preserve the ordering of the $V_{i1}$,
the sets $\{i:\widehat V_{i1}\le t_nx\}$ and $\{i:V_{i1}\le t_nx\}$
are nested for every $x\in[0,T]$.
Consequently,
\begin{align*}
\sum_{i=1}^n
    \big|
        \bm 1(\widehat V_{i1}\le t_nx)
        -
        \bm 1(V_{i1}\le t_nx)
    \big|
&=
    \Big|
        \sum_{i=1}^n
        \bm 1(\widehat V_{i1}\le t_nx)
        -
        \sum_{i=1}^n
        \bm 1(V_{i1}\le t_nx)
    \Big|\\
&=
    \Big|
        \Big\lfloor\frac{n+1}{n}kx\Big\rfloor
        -\sum_{i=1}^n
        \bm 1\big(V_{i1}\le t_nx\big)
    \Big|
\\&= 
    \Big| \Big\lfloor\frac{n+1}{n}kx\Big\rfloor - kx - \sqrt k \alpha_{n1}(x)
    \Big|
    \le
    1+\frac{kT}{n}+\sqrt{k}|\alpha_{n1}(x)|,
\end{align*}
where $\alpha_{n1}$ is the tail empirical process
from \eqref{eq:tail-empirical-process}.
It follows that 
\begin{align*}
    &\sup_{(x,a) \in \mathcal T}
    \frac1k
    \sum_{i=1}^n
    \big|
        g(V_{i2})-a
    \big|^2
    \big|
        \bm 1\big(\widehat V_{i1}\le t_nx\big) - \bm 1\big(V_{i1}\le t_nx\big)
    \big|
    \\ &\hspace{4cm}\le 
    (\|g\|_\infty+\abound)^2
    \Big\{
        \frac1k+\frac Tn
        +\frac1{\sqrt{k}}
        \sup_{x\in[0,T]}|\alpha_{n1}(x)|
    \Big\}
   = O_\Prob(k^{-1/2}).
\end{align*}
Combining the two bounds proves \eqref{eq:two-param-L2}.

Moreover,
\[
    \sup_{(x,a)\in\mathcal T}|\mu_{ni}(x,a)|
    \le
    \frac{\|g\|_\infty+\abound}{\sqrt{k}}\,
    t_nT = (\|g\|_\infty+\abound)T \frac{\sqrt{k}}{n}.
\]
It follows from \eqref{eq:two-param-L2} and $(u+v)^2\le 2u^2+2v^2$ that
\begin{align*}
    &\sup_{(x,a)\in\mathcal T}
    \sum_{i=1}^n
    \big|
        \widehat\mu_{ni}(x,a)-\mu_{ni}(x,a)
    \big|^2
    \\
    &\qquad\le
    2
    \sup_{(x,a)\in\mathcal T}
    \sum_{i=1}^n
    \big|
        \widehat f_{ni}(x,a)-f_{ni}(x,a)
    \big|^2
    +
    2n
    \sup_{(x,a)\in\mathcal T}
    |\mu_{ni}(x,a)|^2
    =
    o_\Prob(1),
\end{align*}
which is Condition~(G).

Conditions~(H) and~(I) are also immediate in the present setting.
Indeed, for each fixed realization, the arrays $(f_{ni}(x,a))$ and $(\widehat f_{ni}(x,a))$ are evaluations of the same two-parameter VC-subgraph class at two different collections of points. The entropy bound established in Proposition~\ref{prop:known-margins-bivariate-2param}, and hence the manageability bound obtained from it in the proof of Theorem~11.20, is uniform in the evaluation points. Thus both arrays are manageable. By the permanence property for sums of manageable arrays, their difference $(\widehat\mu_{ni})=(f_{ni}-\widehat f_{ni})$ is manageable with the sum of the corresponding envelopes.

Furthermore, $\widehat\mu_{ni}$ is separable, and hence AMS by
Lemma~11.15 in \citet{Kosorok2008}. Indeed, one may use the same
countable set $S$
as in the proof of Proposition~\ref{prop:known-margins-bivariate-2param}:
for every realization, the dependence on $x$ is right-continuous and
the dependence on $a$ is continuous.

Finally, for a sufficiently large constant $C<\infty$,
\[
    \widehat F_{ni}
    =
    \frac{C}{\sqrt{k}}
    \left\{
        \bm 1(V_{i1}\le t_nT)
        +
        \bm 1(\widehat V_{i1}\le t_nT)
    \right\}
\]
is an envelope for $\widehat\mu_{ni}$. Since $k^{-1}\sum_{i=1}^n
    \bm 1(V_{i1}\le t_nT) = T + o_\Prob(1)$ and $k^{-1}\sum_{i=1}^n
    \bm 1(\widehat V_{i1}\le t_nT) = T + o_{a.s.}(1)$,
we have
\[
    L\vee
    \sum_{i=1}^n\widehat F_{ni}^2
    = L + o_\Prob(1)
\]
for any $L>4C^2 T$.
This verifies Condition~(I), and completes the proof of \eqref{eq:two-param-bootstrap}.

It remains to show \eqref{eq:hb-bb-approx}. Note that 
$|\widehat a_n(x)|\le2\|g\|_\infty$ uniformly in $x\in[0,T]$, so that
$(x,\widehat a_n(x))\in\mathcal T$. Moreover, we will show below that
\begin{equation}
\label{eq:chat-metric}
    \sup_{x\in[0,T]}
    x|\widehat a_n(x)-\mu_g|^2
    =
    o_\Prob(1).
\end{equation}
The claimed convergence in \eqref{eq:hb-bb-approx} then follows from
asymptotic equicontinuity of $\widehat\Bb_n$ with respect to the canonical
semimetric $\rho$ of $\Hb^\star$, which satisfies
\[
\rho^2((x,a),(x,b))= \Var(\Hb^\star(x,a) - \Hb^\star(x,b)) = x(a-b)^2.
\]
Hence, \eqref{eq:chat-metric} says exactly that
\[
\sup_{x\in[0,T]}
\rho\big((x,\widehat a_n(x)),(x,\mu_g)\big)=o_\Prob(1).
\]
Next, \eqref{eq:two-param-bootstrap} implies unconditional weak convergence
$\widehat\Bb_n\weak\Hb^\star$ in $\ell^\infty(\mathcal T)$; see \citet{BucKoj19}. Since
$\Hb^\star$ is tight and Gaussian, its sample paths are uniformly
$\rho$-continuous \citep[Example~1.5.10]{VanWel23}, and hence
$\widehat\Bb_n$ is asymptotically uniformly $\rho$-equicontinuous in
probability \citep[Addendum~1.5.8]{VanWel23}. Writing
$D_n=\sup_{x\in[0,T]}|\widehat\Bb_n(x,\widehat a_n(x))-\widehat\Bb_n(x,\mu_g)|$,
we obtain, for all $\eps,\delta>0$,
\[
\Prob(D_n>\eps)
\le
\Prob\Big(\sup_{x\in[0,T]}\rho\big((x,\widehat a_n(x)),(x,\mu_g)\big)\ge\delta\Big)
+\Prob\Big(\sup_{u,v\in\mathcal T:\,\rho(u,v)<\delta}|\widehat\Bb_n(u)-\widehat\Bb_n(v)|>\eps\Big),
\]
such that $D_n=o_\Prob(1)$. Finally, by Markov's inequality,
$\Prob\big(\Prob_\xi(D_n>\eps)>\kappa\big)\le\kappa^{-1}\Prob(D_n>\eps)\to0$
for all $\eps,\kappa>0$, which is \eqref{eq:hb-bb-approx}.

To verify \eqref{eq:chat-metric}, write $a_n(x)=\ggcmean(t_nx)$ and fix $\eta \in (0,T)$. By definition,  $a_n(x)$ is bounded by $\|g\|_\infty$, uniformly in $x \ge 0$. Moreover,  
$\widehat a_n$ is bounded by $2 \|g\|_\infty$; indeed, for $x>0$,
\[
    |\widehat a_n(x)|
    \le
    \frac{\|g\|_\infty}{kx}
    \sum_{i=1}^n
    \bm 1(\widehat V_{i1}\le t_nx)
    \le
    \frac{\|g\|_\infty}{kx}
    \frac{n+1}{n}kx
    \le 
    \frac{n+1}{n}\|g\|_\infty
    \le 2\|g\|_\infty.
\]
Consequently, for sufficiently large $n$,
\[
    \sup_{x \in (0,\eta)}
    x\big|
        \widehat a_n(x)-a_n(x)
    \big|^2
    \le 9 \|g\|_\infty^2  \eta.
\]
On the other hand, for $x\in [\eta,T]$, we can write 
$
    x \big|
        \widehat a_n(x)-a_n(x)
    \big|^2
    =
    x^{-1}\big| \gghat(x)-\ggmean(x) \big|^2 
    =
    x^{-1}|\ggprocess(x)|^2/ k.
$
By Corollary~\ref{cor:unknown-margins-bivariate}, we have $\sup_{x\in[0,T]}|\ggprocess(x)|^2=O_\Prob(1)$, 
and thus
\[
    \sup_{x\in [\eta,T]} x
    \big|
        \widehat a_n(x)-a_n(x)
    \big|^2
    \le \frac1{\eta k} \sup_{x\in[0,T]}|\ggprocess(x)|^2
    =
    O_\Prob(k^{-1}).
\]
Letting first $n\to\infty$ and then $\eta\downarrow0$ proves
\[
\sup_{x\in[0,T]}
    x|\widehat a_n(x)-a_n(x)|^2
    =
    o_\Prob(1)
\]
and  the claim in \eqref{eq:chat-metric} then follows from
$
\sup_{x\in[0,T]}|a_n(x)-\mu_g|=o(1)
$ by Condition~\ref{cond:conditional-weak-limit-bivariate}.
\end{proof}

\section{Auxiliary results}
\label{sec:auxiliary}

\begin{proposition}
\label{prop:known-margins-bivariate-2param}
Let $D$ be a bivariate copula, and suppose that $(V_1,V_2)\sim D$ satisfies Condition~\ref{cond:conditional-weak-limit-bivariate}. Let $(V_{i1},V_{i2})_{i=1}^n$ be $n$ IID observations with $(V_{i1},V_{i2})\sim D$, let $g:[0,1]\to\mathbb R$ be continuous, and let $k\in[n]$ with $k\to\infty$ and $k=o(n)$ as $n\to\infty$.
For $x\ge0$ and $a\in\R$, define
\begin{align*}
    f_{n,x,a}(v_1,v_2)
    =
    \sqrt{\frac nk}
    \big\{g(v_2)-a\big\}
    \bm 1\Big(v_1\le\frac{kx}{n}\Big)
\end{align*}
and
\begin{align*}
    \Hb_n^\star(x,a)
    =
    \frac1{\sqrt n}
    \sum_{i=1}^n
    \Big\{
        f_{n,x,a}(V_{i1},V_{i2})
        -
        \Exp[f_{n,x,a}(V_{i1},V_{i2})]
    \Big\}.
\end{align*}
Then, for every $T>0$ and $\abound>0$,
\begin{align}
\label{eq:weak_conv_delta-2param}
    \big(
        \Hb_n^\star(x,a)
    \big)_{(x,a)\in[0,T]\times[-\abound,\abound]}
    \weak
    \big(
        \Hb^\star(x,a)
    \big)_{(x,a)\in[0,T]\times[-\abound,\abound]}
\end{align}
in $\ell^\infty([0,T]\times[-\abound,\abound])$, where $\Hb^\star$ is a
tight, centered Gaussian process with covariance function
\[
    \Cov\big[
        \Hb^\star(x,a),
        \Hb^\star(y,b)
    \big]
    =
    (x\wedge y)
    \Exp_{Z\sim H}
    \big[
        \{g(Z)-a\}\{g(Z)-b\}
    \big].
\]
The process $\Hb^\star$ admits a version whose sample paths are
continuous on $[0,\infty)\times[-\abound,\abound]$.
\end{proposition}

\begin{proof}[Proof of Proposition~\ref{prop:known-margins-bivariate-2param}]
Fix $T,\abound>0$, write $t_n=k/n$, and let $\mathcal T=[0,T]\times[-\abound,\abound]$.
For $n\in\N$, let
\[
    \mathcal F_{n,T}
    =
    \{
        f_{n,x,a}:(x,a)\in\mathcal T
    \}.
\]
A corresponding envelope is
\[
    F_{n,T}(v_1,v_2)
    =
    (\|g\|_\infty+\abound)t_n^{-1/2}
    \bm 1\big(v_1\le t_n T\big).
\]
We apply Theorem~11.20 in \citet{Kosorok2008}, with the index set
appearing therein replaced by $\mathcal T$. 
We verify its conditions in turn, which are as follows:
\begin{itemize}
    \item[(1)]  The sequence $\mathcal F_{n,T}$  satisfies the bounded uniform entropy integral condition
    \begin{align*}
        \limsup_{n\to\infty} \sup_Q \int_0^1 \sqrt{ \log N(\eps \|F_{n,T}\|_{Q,2}, \mathcal  F_{n,T}, L_2(Q))}\, {\mathrm d}\eps<\infty,
    \end{align*}
    where the inner supremum ranges over all probability measures $Q$ with finite support contained in $[0,1]^2$ for which $\|F_{n,T}\|_{Q, 2}^2=\int F_{n,T}^2\, \diff Q >0$;
    
    \item[(2)] The sequence $\mathcal F_{n,T}$ is almost measurable Suslin (AMS), i.e., for all $n\in \N$, there exists a Suslin topological space $S_n\subset \mathcal T$ with Borel sets $\mathcal B_n$ such that 
    \begin{itemize}
        \item[(a)] $\mathbb P^*\left(\sup_{ (x,a) \in \mathcal T} \inf_{(y,b)\in S_n}|f_{n, x,a}(V_1, V_2)-f_{n,y,b}(V_1, V_2)|>0 \right)=0$; and
        \item[(b)] the functions $f_{n, \cdot, \cdot}:[0,1]^2\times S_n \to \mathbb R$ are $\mathcal B|_{[0,1]^2}\times \mathcal B_n$ measurable for all $n$;
    \end{itemize}
    
    \item[(3)] $B(x,a ; y,b) =\lim_{n\to\infty}\Exp[\Hb_n^\star(x,a)\Hb_n^\star(y,b)]$ exists for every $(x,a), (y,b) \in \mathcal T$;
    
     \item[(4)] $\lim_{n\to\infty} \Exp F_{n,T}^2(V_1, V_2)< \infty$;
     
    \item[(5)] $\lim_{n\to\infty} \Exp F_{n,T}^2(V_1, V_2)\bm 1(F_{n,T}(V_1, V_2)>\eps \sqrt n)=0$ for all $\eps>0$;
    
    \item[(6)] The limit $\rho(x,a; y,b)\coloneqq \lim_{n\to\infty} \rho_n(x,a;  y,b)$ exists for every $(x,a),(y,b)\in\mathcal T$, where $\rho_n(x,a;y,b)=\sqrt{\Exp|f_{n,x,a}(V_1, V_2)-f_{n, y,b}(V_1, V_2)|^2}$. Furthermore, for all deterministic sequences $(x_n,a_n)_n, (y_n,b_n)_n$ in $\mathcal T$, if $\rho(x_n,a_n; y_n, b_n)\to 0$ as $n\to \infty$ then $\rho_n(x_n,a_n; y_n,b_n)\to 0$.
    \end{itemize}

\medskip

\noindent
\textit{Verification of (1):} Let
\[
    \mathcal G
    =
    \big\{
        v_2\mapsto g(v_2)-a:
        a\in[-\abound,\abound]
    \big\}
, \qquad
    \mathcal I_{n,T}
    =
    \left\{
        v_1\mapsto
        \bm 1(v_1\le t_nx):
        x\in[0,T]
    \right\}.
\]
Both classes are VC-subgraph classes with VC-index bounded uniformly
in $n$. Indeed, the subgraphs associated with $\mathcal G$ are nested
in the parameter $a$, while $\mathcal I_{n,T}$ is a one-dimensional
threshold class.
Moreover, writing  $h_a(v_2)=g(v_2)-a$ and $I_x(v_1)=\bm 1(v_1 \le t_n x)$,
we have
\begin{align*}
    &\big\{
        (v_1,v_2,z):
        z<h_a(v_2)I_x(v_1)
    \big\}
    =
    \big(
        \{I_x=1\}\cap\{z<h_a(v_2)\}
    \big)
    \cup
    \big(
        \{I_x=0\}\cap\{z<0\}
    \big).
\end{align*}
Hence the subgraph class associated with
\[
    \big\{
        (v_1,v_2)\mapsto h_a(v_2)I_x(v_1):
        (x,a)\in\mathcal T
    \big\}
\]
is obtained from VC classes by a finite number of complements,
intersections and unions. By the permanence properties of VC classes,
it is therefore itself a VC class, with VC-index bounded uniformly in
$n$. Multiplication by the positive constant $\sqrt{n/k}$ does not
affect the VC-subgraph property. Thus $\mathcal F_{n,T}$ is a
VC-subgraph class with uniformly bounded VC-index, and Condition~(1)
follows from Lemma~11.21 in \citet{Kosorok2008}.

\medskip

\noindent
\textit{Verification of (2):}
Let
$
    S_n
    =
    \big((\Q\cap[0,T])\cup\{T\}\big)
    \times
    \big((\Q\cap[-\abound,\abound])\cup\{-\abound,\abound\}\big).
$
For fixed $n$, $(x,a)\in\mathcal T$, and $\omega\in\Omega$, choose
$(x_m,a_m)\in S_n$ such that $
    x_m\downarrow x$ and $
    a_m\to a$.
By right-continuity of the indicator as a function of $x$ and
continuity in $a$,
\[
    f_{n,x_m,a_m}
    (V_{i1}(\omega),V_{i2}(\omega))
    \longrightarrow
    f_{n,x,a}
    (V_{i1}(\omega),V_{i2}(\omega))
\]
for every $i\in[n]$. Since $n$ is fixed,
\[
    \inf_{(y,b)\in S_n}
    \sum_{i=1}^n
    \left\{
        f_{n,y,b}(V_{i1}(\omega),V_{i2}(\omega))
        -
        f_{n,x,a}(V_{i1}(\omega),V_{i2}(\omega))
    \right\}^2
    =0.
\]
Thus the triangular array is separable in the sense of
\citet[Page~219]{Kosorok2008}, and Condition~(2) follows from
Lemma~11.15 therein.

\medskip

\noindent
\textit{Verification of (3):}
For $(x,a),(y,b)\in\mathcal T$,
\begin{align*}
    \Exp\big[
        \Hb_n^\star(x,a)\Hb_n^\star(y,b)
    \big]
    &=
    \Cov\big[
        f_{n,x,a}(V_1,V_2),
        f_{n,y,b}(V_1,V_2)
    \big]
    \\
    &=
    B_n(x,a;y,b)
    -
    \Exp[f_{n,x,a}(V_1,V_2)]
    \Exp[f_{n,y,b}(V_1,V_2)],
\end{align*}
where
\[
    B_n(x,a;y,b)
    =
    \Exp\big[
        f_{n,x,a}(V_1,V_2)
        f_{n,y,b}(V_1,V_2)
    \big].
\]
First,
\begin{align*}
    \big|
        \Exp[f_{n,x,a}(V_1,V_2)]
    \big|
    \le
    t_n^{-1/2}
    \Exp\left[
        |g(V_2)-a|
        \bm 1(V_1 \le t_n x)
    \right]
    &\le
    (\|g\|_\infty + \abound) x t_n^{1/2}
    \\&\le
    (\|g\|_\infty + \abound) T t_n^{1/2}
    =
    o(1),
\end{align*}
uniformly in $(x,a)\in\mathcal T$. Hence the product of the two
expectations converges to zero uniformly.
For the treatment of the product moment term, let $m_1=\Exp_{Z\sim H}[g(Z)]$ and $m_2=\Exp_{Z\sim H}[g^2(Z)]$.
For $u>0$, define
\[
    m_{1,n}(u)
    =
    \Exp\left[
        g(V_2)
        \,\middle|\,
        V_1\le t_n u
    \right],
    \qquad
    m_{2,n}(u)
    =
    \Exp\left[
        g^2(V_2)
        \,\middle|\,
        V_1\le t_n u
    \right].
\]
Condition~\ref{cond:conditional-weak-limit-bivariate} and the bounded
continuity of $g$ and $g^2$ imply
\begin{equation}
\label{eq:uniform-moments-a}
    \sup_{u\in(0,T]}
    |m_{1,n}(u)-m_1|
    =o(1),
    \qquad
    \sup_{u\in(0,T]}
    |m_{2,n}(u)-m_2|
    =o(1), \qquad n \to \infty.
\end{equation}
Let $u=x\wedge y$. If $u=0$, then $B_n(x,a;y,b)=0$. If $u>0$, then
\begin{align}
    B_n(x,a;y,b)
    &= \nonumber
    t_n^{-1}
    \Exp\left[
        \{g(V_2)-a\}\{g(V_2)-b\}
        \bm 1\big(V_1\le t_n u\big)
    \right]
    \\
    &= \nonumber
    u
    \Exp\left[
        \{g(V_2)-a\}\{g(V_2)-b\}
        \,\middle|\,
        V_1\le t_nu
    \right]
    \\
    &= \label{eq:def_bn}
    u
    \left\{
        m_{2,n}(u)
        -
        (a+b)m_{1,n}(u)
        +
        ab
    \right\}.
\end{align}
By \eqref{eq:uniform-moments-a},
\begin{align}
\label{eq:uniform-Bn}
    \sup_{(x,a), (y,b) \in \mathcal T}
    \left|
        B_n(x,a;y,b)
        -
        B(x,a;y,b)
    \right|
    \longrightarrow0,
\end{align}
where
\begin{align*}
    B(x,a;y,b)
    &=
    (x\wedge y)
    \left\{
        m_2-(a+b)m_1+ab
    \right\}
    \nonumber\\
    &=
    (x\wedge y)
    \Exp_{Z\sim H}
    \big[
        \{g(Z)-a\}\{g(Z)-b\}
    \big].
\end{align*}
Overall,
$
    \lim_{n\to\infty}
    \Exp\big[
        \Hb_n^\star(x,a)\Hb_n^\star(y,b)
    \big]
    =
    B(x,a;y,b),
$
which verifies Condition~(3).

\medskip

\noindent
\textit{Verification of (4):}
We have $
    \Exp[
        F_{n,T}^2(V_1,V_2)
    ]
    =
    (\|g\|_\infty  +\abound)^2 t_n^{-1}
    \Prob(V_1 \le t_n T)
=
    (\|g\|_\infty  +\abound)^2 T.
$
Thus Condition~(4) holds.

\medskip

\noindent
\textit{Verification of (5):}
Since
$
    F_{n,T}(v_1,v_2)
    \le
    (\|g\|_\infty  +\abound) t_n^{-1/2},
$
the inequality $
    F_{n,T}(V_1,V_2)>\eps\sqrt n$
can hold only if
$
    (\|g\|_\infty  +\abound) /{\sqrt k}>\eps.
$
Since $k\to\infty$, this is impossible for sufficiently large $n$.
Consequently,
\[
    \Exp\left[
        F_{n,T}^2(V_1,V_2)
        \bm 1\big(
            F_{n,T}(V_1,V_2)>\eps\sqrt n
        \big)
    \right]
    =0
\]
for all sufficiently large $n$, which proves Condition~(5).

\medskip

\noindent
\textit{Verification of (6):}
For $(x,a),(y,b)\in\mathcal T$, recalling the definition  of $B_n$ in \eqref{eq:def_bn}, we have
\begin{align*}
    \rho_n^2\big(x,a; y,b)
    &=
    \Exp
    \left|
        f_{n,x,a}(V_1,V_2)
        -
        f_{n,y,b}(V_1,V_2)
    \right|^2
\\&=
    B_n(x,a;x,a)
    -
    2B_n(x,a;y,b)
    +
    B_n(y,b;y,b).
\end{align*}
Therefore, by \eqref{eq:uniform-Bn},
\begin{equation}
\label{eq:uniform-rho}
    \sup_{(x,a),(y,b)\in\mathcal T}
    \left|
        \rho_n^2\big(x,a;y,b\big)
        -
        \rho^2\big(x,a;y,b\big)
    \right|
    \longrightarrow0,
\end{equation}
where
\begin{align*}
    \rho^2\big(x,a;y,b\big)
    &=
    x\Exp_{Z\sim H}\big[\{g(Z)-a\}^2\big]
    -2(x\wedge y)
    \Exp_{Z\sim H}
    \big[
        \{g(Z)-a\}\{g(Z)-b\}
    \big]
    \\
    &\hspace{6cm}
    +y\Exp_{Z\sim H}\big[\{g(Z)-b\}^2\big].
\end{align*}
Consequently, if $(x_n,a_n),(y_n,b_n)\in\mathcal T$
are deterministic sequences such that
$
    \rho\big(x_n,a_n;y_n,b_n\big) = o(1),
$
then, by \eqref{eq:uniform-rho},
\begin{align*}
    \rho_n^2\big(x_n,a_n;y_n,b_n\big)
    &\le
    \sup_{(x,a),(y,b)\in\mathcal T}
    \left|
        \rho_n^2\big(x,a;y,b\big)
        -
        \rho^2\big(x,a;y,b\big)
    \right|
    +
    \rho^2\big(x_n,a_n;y_n,b_n\big) = o(1).
\end{align*}
This verifies Condition~(6).

\medskip

All conditions of Theorem~11.20 in \citet{Kosorok2008} are therefore
satisfied, and we obtain the claimed weak convergence result in \eqref{eq:weak_conv_delta-2param}.  Moreover, Theorem~11.20 also implies that the limiting Gaussian process has sample paths that are almost surely uniformly continuous with respect to its canonical semimetric $\rho$. Since the explicit expression for $\rho$ obtained above is continuous in its arguments, $\rho$-continuity implies continuity with respect to the usual topology on $\mathcal T$. Since $T>0$ was arbitrary, the same conclusion holds on $[0,\infty)\times[-\abound,\abound]$, by considering $T=N\in\N$ and noting that continuous versions on nested index sets coincide almost surely on the smaller set.
\end{proof}

\section{Additional simulation results}
\label{sec:additional-simulations}

\subsection{On the choice of \texorpdfstring{$k$}{k} in bivariate RMLMs}\label{sec:additional-simulations_doublegrid}

The goal of this section is to study the behaviour of the mean squared error (MSE) of the two estimators $\widehat \Gamma_{12}$ and $\widehat \Gamma_{21}$ in the bivariate RMLM setup introduced in Section~\ref{sec:simulations}. We use $k=\max(\lfloor c\, n^{\kappa}\rfloor,10)$ and evaluate the performance of the two estimators in a two-dimensional $(c,\kappa)$-grid, where $c\in\{0.25, 0.5, 0.75, 1\}$ and $\kappa\in\{0.3, 0.4, 0.5, 0.6\}$. We report the results for the sample size $n=50000$ and $R=1000$ replicates.

\begin{figure}[htbp]
    \centering
    \includegraphics[width=.85\linewidth]{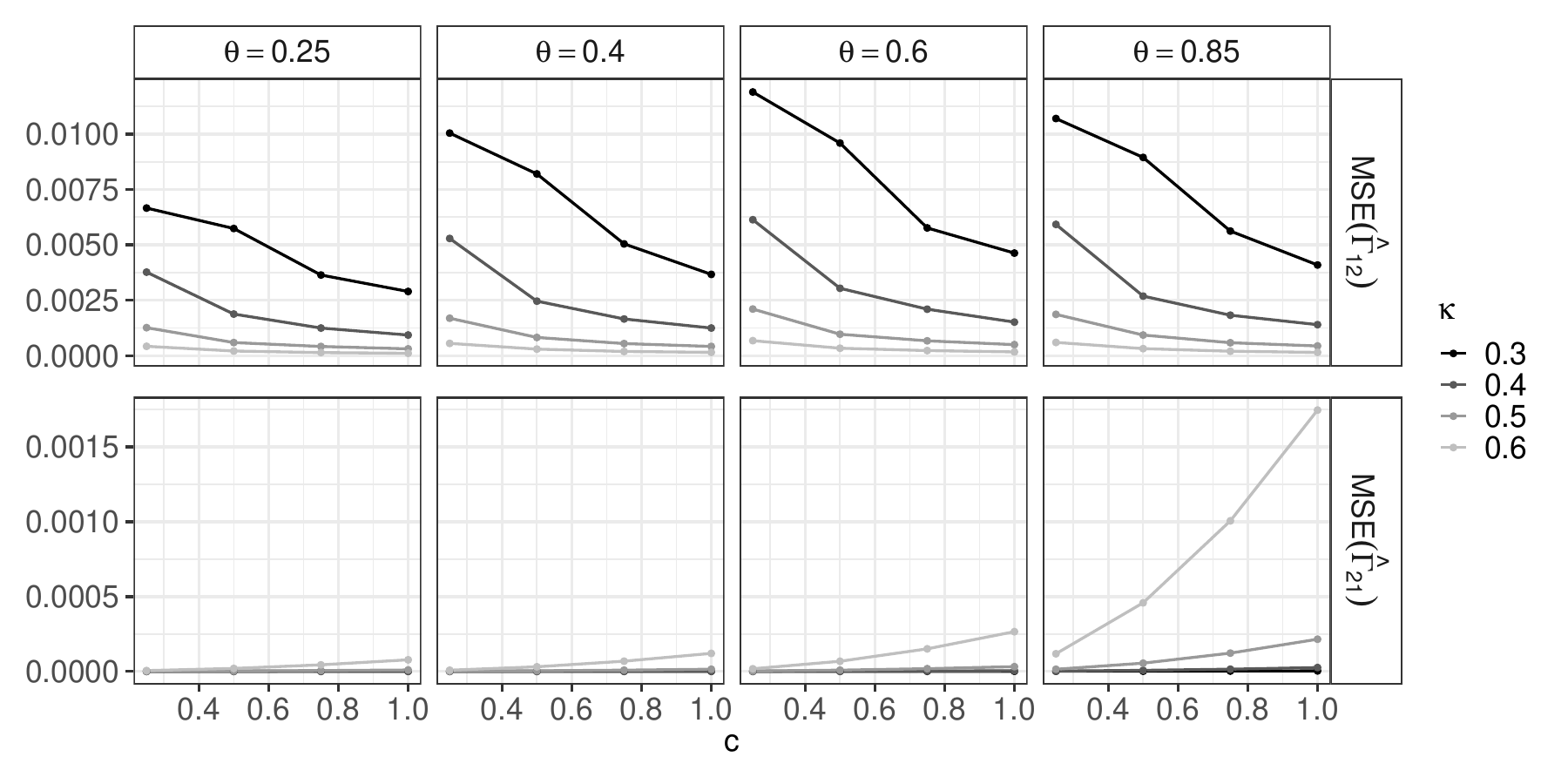} 
    \caption{Bivariate RMLM: MSE of $\widehat\Gamma_{12}$ (top) and $\widehat\Gamma_{21}$ (bottom) based on
    $R=1000$ replications for different values of
    $\theta$, $n=50000$, and $k=\max(\lfloor c n^{\kappa} \rfloor,10)$ with $c\in\{0.25, 0.5, 0.75, 1\}$ and $\kappa\in\{0.3, 0.4, 0.5, 0.6\}$.}
   \label{fig:MSE_gamma_double_grid}
\end{figure}

Figure~\ref{fig:MSE_gamma_double_grid} shows a contrasting effect of the scalar $c$ and the exponent $\kappa$ in the two directions. On the one hand, larger $c$ and $\kappa$ yield a low MSE when estimating $\Gamma_{12}$, whereas they lead to higher MSE when estimating $\Gamma_{21}$. In general, the MSE for $\widehat \Gamma_{21}$ is much lower than for $\widehat \Gamma_{12}$. On the other hand, smaller $c$ and $\kappa$ offer a low MSE for $\widehat \Gamma_{21}$ and a high MSE for $\widehat \Gamma_{21}$. The exponents $\kappa\in\{0.3, 0.5, 0.6\}$ appear to be the most affected by such phenomenon, pointing to $\kappa=0.4$ as offering a balanced trade-off between $\Gamma_{12}$ and $\Gamma_{21}$.

\subsection{On CTC estimation}

\begin{figure}[htbp]
    \centering
    \includegraphics[width=.85\linewidth]{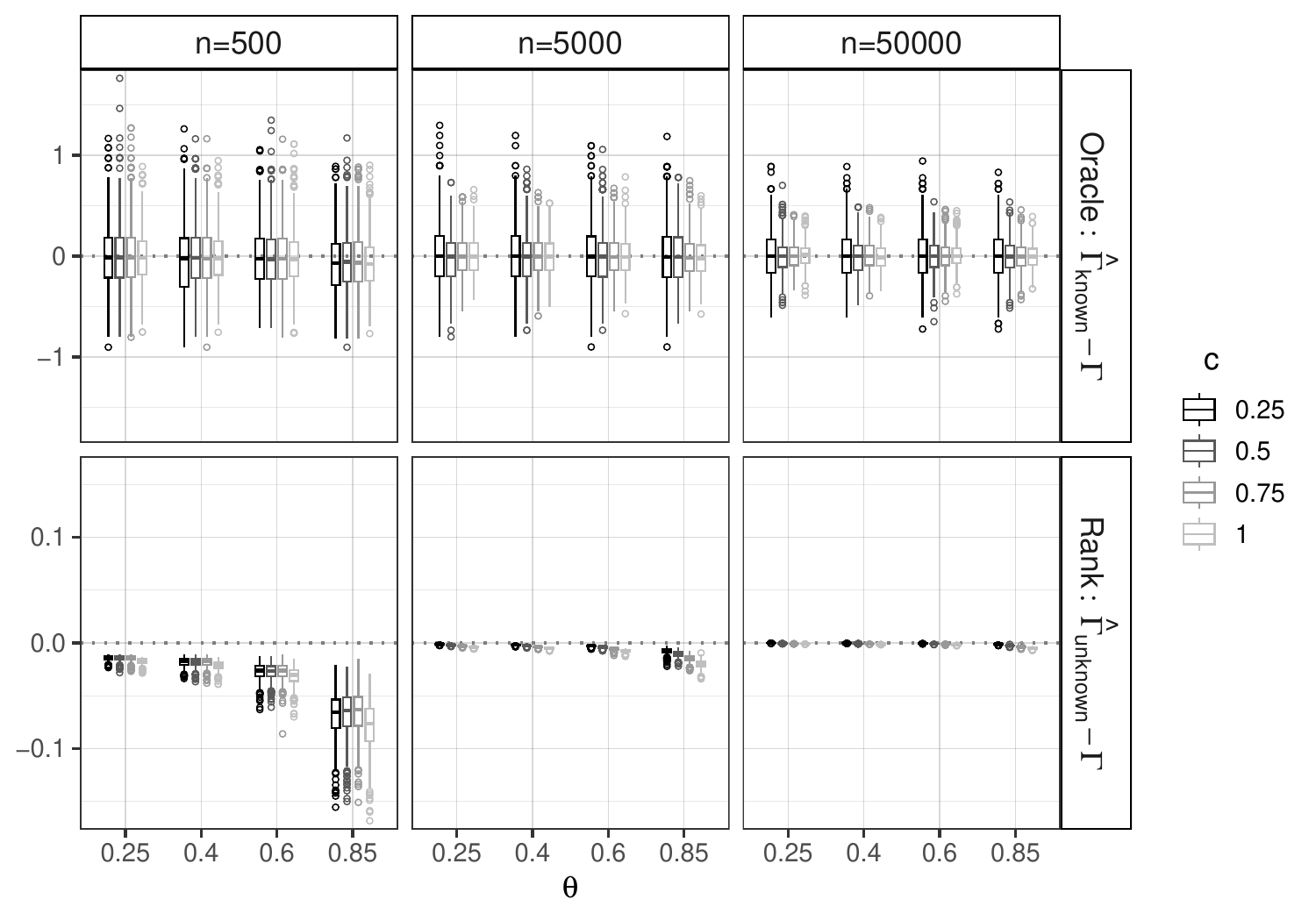} 
    \caption{Bivariate RMLM: Boxplots of the estimation
    errors $\widehat\Gamma_{21}^*-\Gamma_{21}$ (top) and $\widehat\Gamma_{21}-\Gamma_{21}$ (bottom) based on
    $R=1000$ replications for different values of
    $\theta$, $n$, and $k$.}
   \label{fig:MLgamma21}
\end{figure}

Figure~\ref{fig:MLgamma21} complements the simulation results
in Section~\ref{sec:simulations} by displaying the estimation
errors for $\Gamma_{21}=1$ under the same parameter configurations.
Consistent with the degenerate limiting distribution, the rank-based
estimator exhibits almost no variability, in sharp contrast to
the oracle estimator. However, it shows a negative finite-sample
bias that becomes more pronounced for larger values of $\theta$
and $k$. This bias generally decreases with increasing sample
size, illustrating the importance of threshold selection even
when the asymptotic variance vanishes.

\subsection{On variance estimation}
\label{app:variance21}

\begin{figure}[htbp]
    \centering
    \includegraphics[width=.85\linewidth]
        {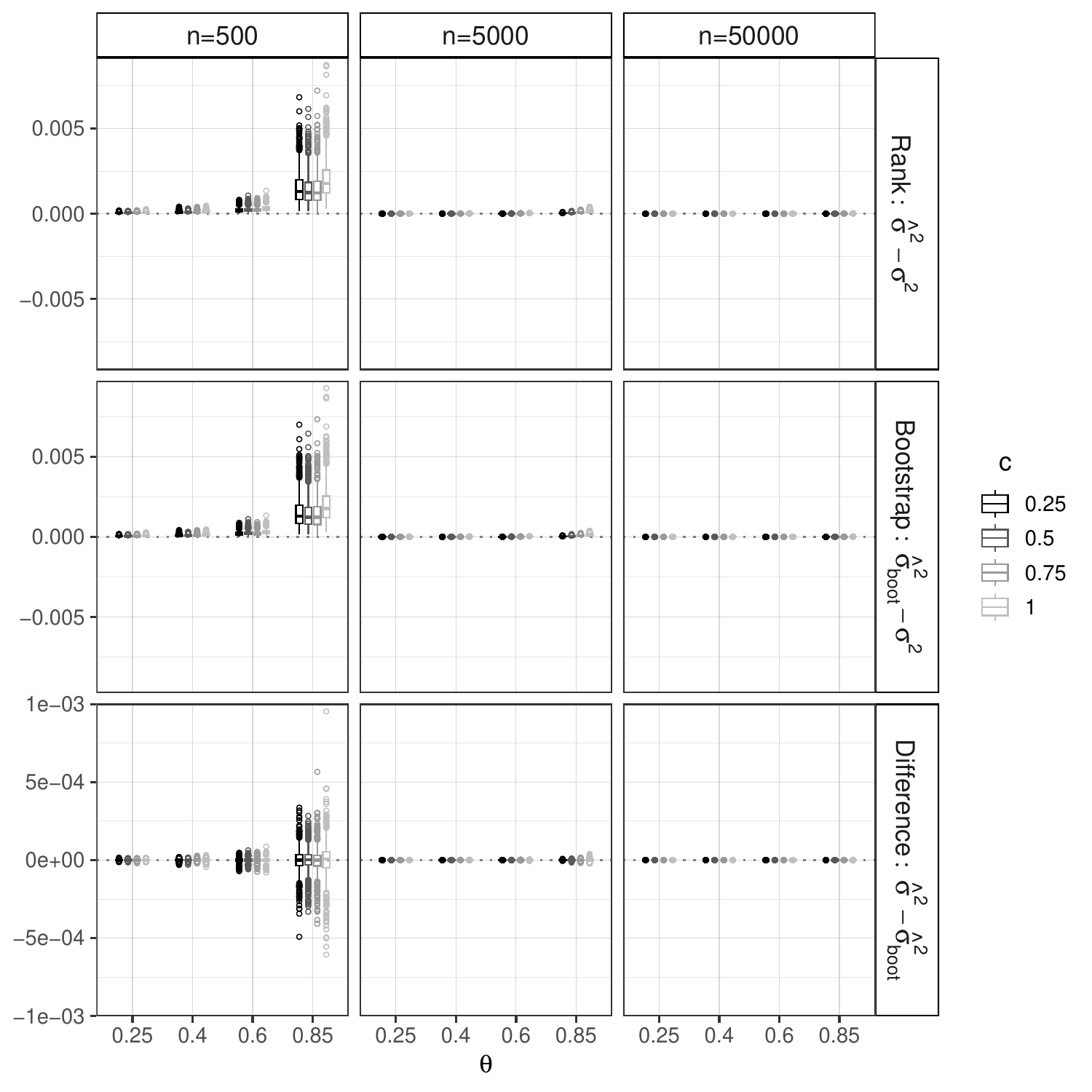}
     \caption{Bivariate RMLM, direction $(2,1)$:
    Boxplots of the estimation errors of the plug-in
    variance estimator (top) and the multiplier-bootstrap
    variance estimator (middle), together with their
    differences (bottom).
    Results are based on $R=1000$ simulation replications
    with $1000$ multiplier replicates each.}
    \label{fig:ML21vars}
\end{figure}

Figure~\ref{fig:ML21vars} complements the results of
Section~\ref{sec:variance-simulations} by considering
variance estimation in the direction $(2,1)$, where
$\sigma_{21}^2=0$.
Both estimators exhibit a positive finite-sample bias,
which is most pronounced for large values of $a_{11}$
and $k$. This bias decreases as the sample size increases
or the number of exceedances decreases.
Despite the degenerate limiting variance, the plug-in
and bootstrap estimators show broadly similar behaviour,
with their differences generally centred around zero.

\subsection{On the choice of \texorpdfstring{$k$}{k} in graph recovery}
\label{app:k_graph}

We carry out a simulation study on the choice of $k$ used to estimate the CTC matrix $(\widehat\Gamma_{j\ell})_{(j,\ell)\in\Pc_2}$ of a $d$-dimensional linear SEM under the setup of Section~\ref{subsec:sims-graph}. More specifically, the goal is to recover the confidence graphs from $\widehat \Gamma$. We again use $k=\max(\lfloor c\, n^{\kappa}\rfloor,10)$ and evaluate the performance of the joint coverage metric on a two-dimensional $(c,\kappa)$-grid, where $c\in\{0.25, 0.5, 0.75, 1\}$ and $\kappa\in\{0.3, 0.4, 0.5, 0.6\}$. We report the results for the sample size and $R=100$ random DAGs replicates for each combination of $d\in\{5,10,20\}$, $p\in\{0.1, 0.2\}$ and $n\in\{500, 5000, 50000\}$. 

\begin{figure}[htbp]
    \centering
    \includegraphics[width=.85\linewidth]{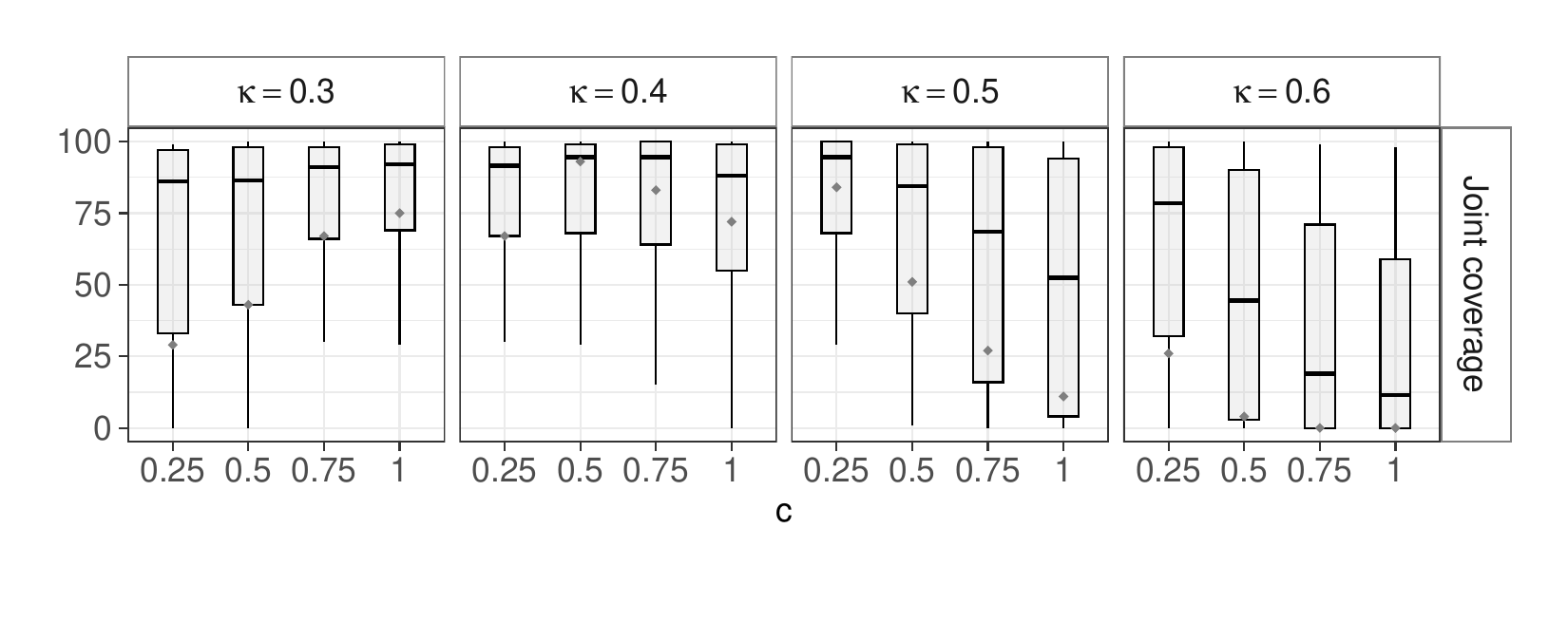} 
    \caption{Boxplots of the joint coverage metric aggregated over different sample sizes $n$, dimensions $d$ and edge probabilities $p$. The metric is computed for $R=100$ replications. The gray diamonds show the joint coverage for $(d, p, n)=(20, 0.2, 50000)$.}
   \label{fig:k_coverage}
\end{figure}

Figure~\ref{fig:k_coverage} contains the results of the analysis, with each boxplot bar summarising the joint coverage over all $(d, p, n)$ configurations for a given $(c,\kappa)$. The second panel, corresponding to the exponent $\kappa=0.4$, yields a relatively stable performance across different scalars $c$. In particular, the values $(c,\kappa)=(0.5, 0.4)$ help achieve a relatively good joint coverage of complex graph configuration $(d,p)=(20, 0.2)$ for large sample sizes. We therefore choose $\kappa=0.4$ for the simulations of Section~\ref{subsec:sims-graph}.

\subsection{On graph recovery}
\label{app:graph}

\begin{figure}[htbp]
    \centering
    \includegraphics[width=.85\linewidth]{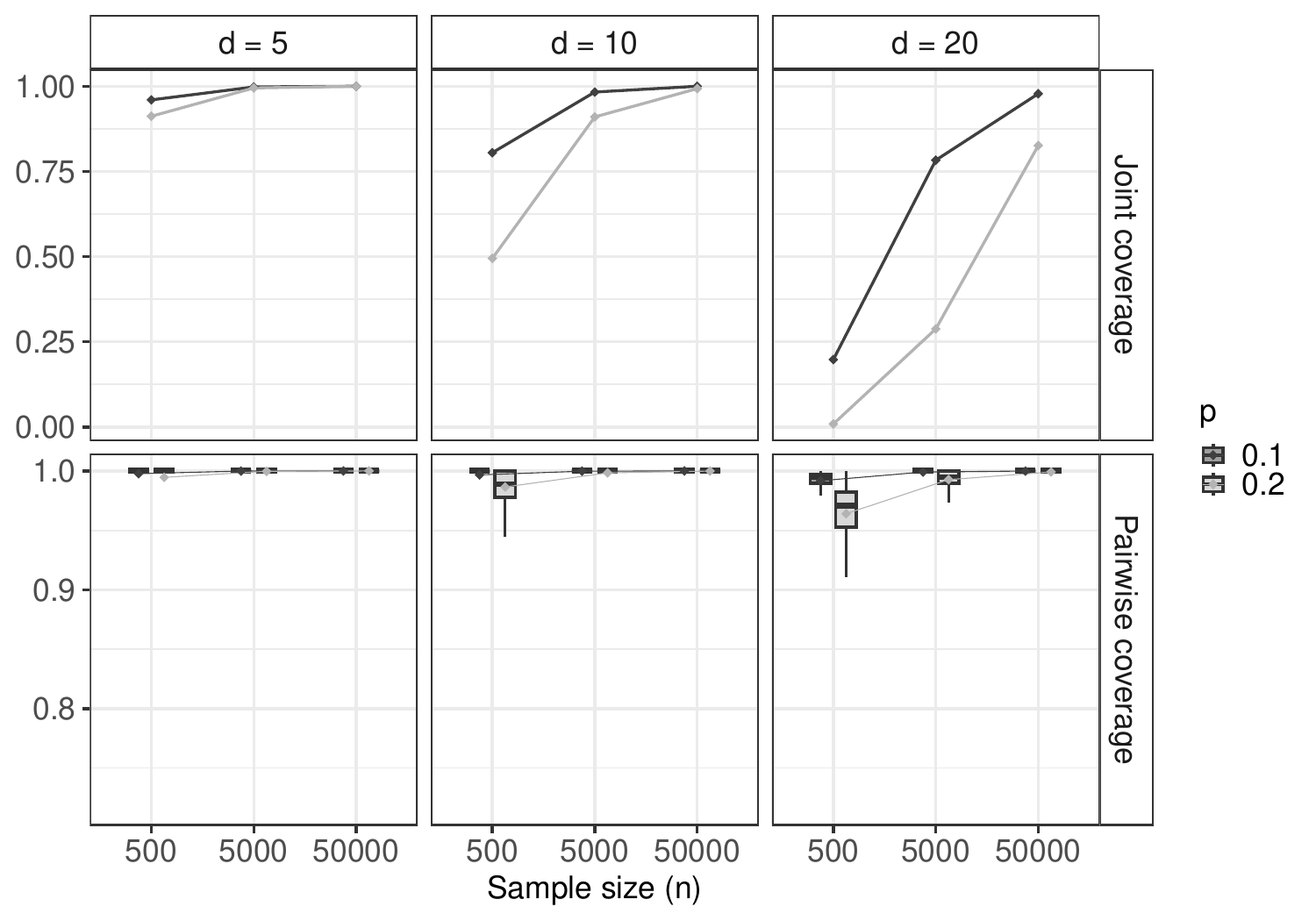} 
    \caption{Empirical joint coverage (top panel) and pairwise coverage (bottom panel) of the lower and upper graphs, with respective edge sets $\widehat E^L_{\gamma,\alpha}$ and $\widehat E^U_{\gamma,\alpha}$ for $R=1000$ replications. The diamonds indicate the averages.}
   \label{fig:coverage}
\end{figure}

Figure~\ref{fig:coverage} provides a visual illustration of the coverage metrics from Table~\ref{tab:graph-recovery}, showing the difficulty~in obtaining joint coverage for complex graphs with $d=20$. The effect of the latter is also visible from the pairwise coverage, translating into higher uncertainty in particular for smaller sample sizes.

\begin{figure}[htbp]
    \centering
    \includegraphics[width=.85\linewidth]
        {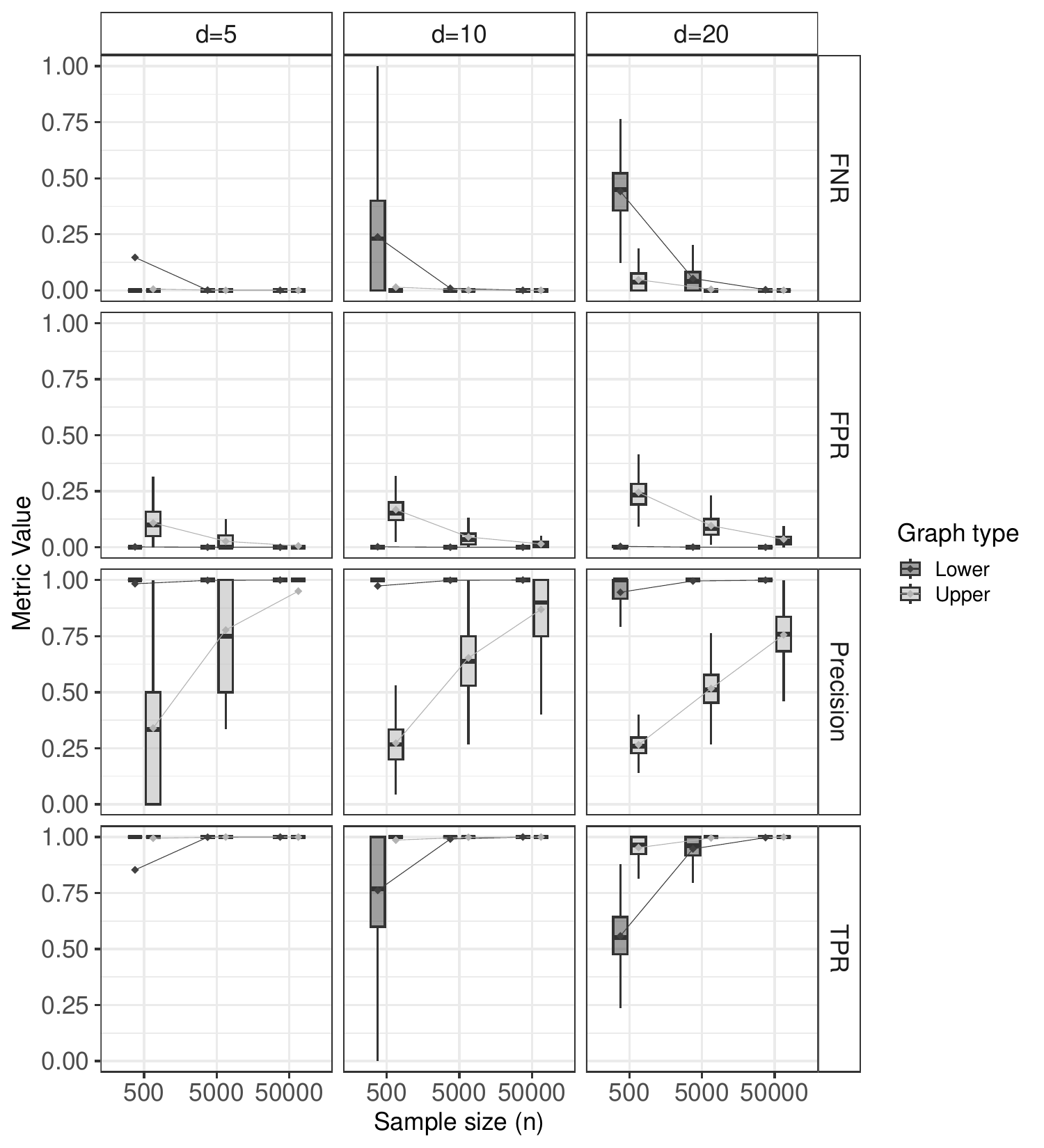}
    \caption{Distributions of the recovery metrics for the lower and upper
    confidence graphs
    across $R=1000$ simulation replications for
    $p=0.1$, different dimensions $d$ and sample
    sizes~$n$.} 
    \label{fig:LSEM-p01}
\end{figure}

\begin{figure}[htbp]
    \centering
    \includegraphics[width=.85\linewidth]
        {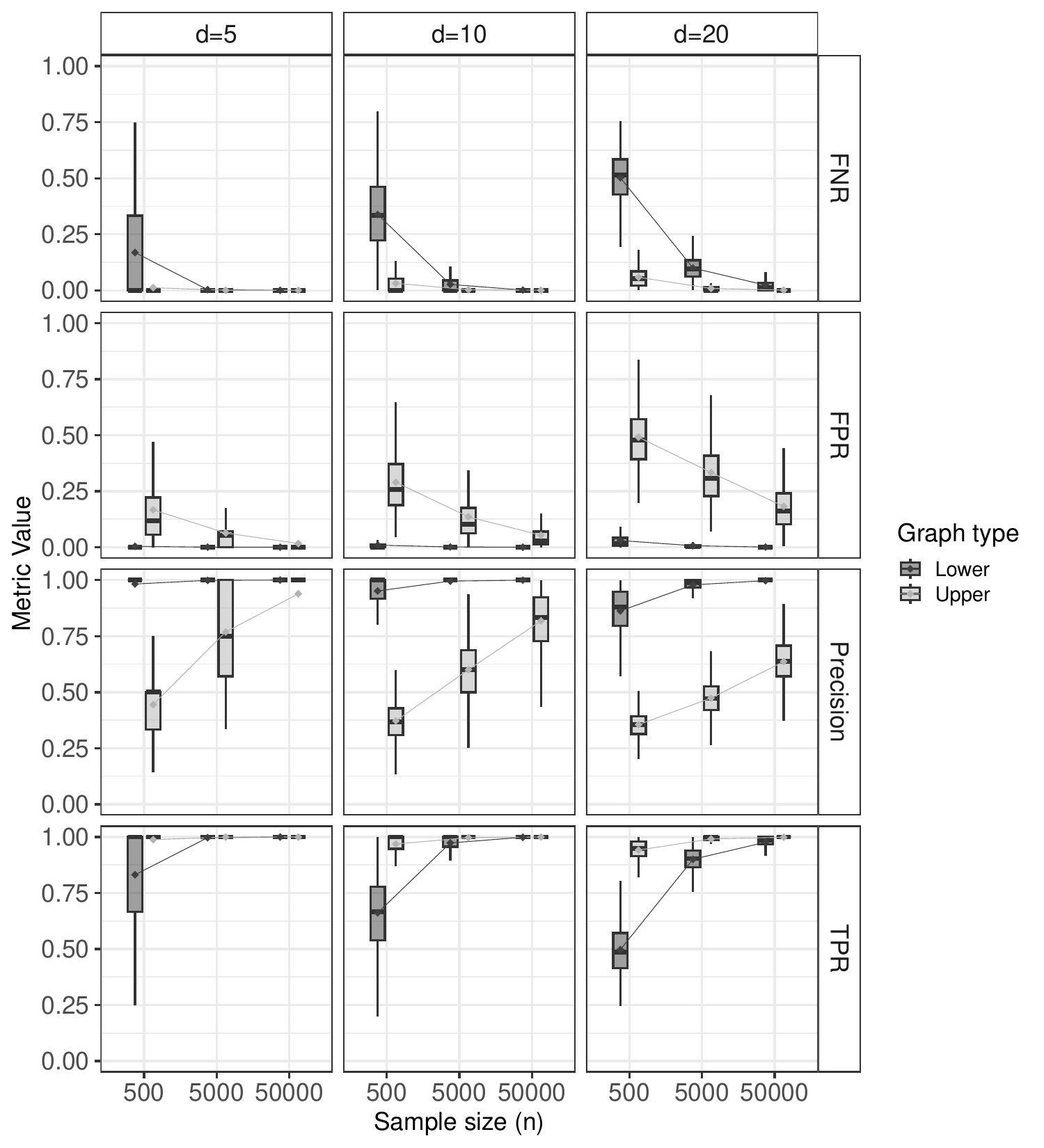}
        \caption{Distributions of the recovery metrics for the lower and upper
    confidence graphs
    across $R=1000$ simulation replications for
    $p=0.2$, different dimensions $d$ and sample
    sizes~$n$.}
    \label{fig:LSEM-p02}
\end{figure}

Figures~\ref{fig:LSEM-p01} and~\ref{fig:LSEM-p02}
complement Table~\ref{tab:graph-recovery} by showing
the distributions of the individual recovery metrics
across simulation replications for $p=0.1$ and $p=0.2$,
respectively.
The results confirm the generally conservative behaviour
of the lower confidence graphs, which exhibit high
precision and low false positive rates, but can have
substantial false negative rates in smaller samples.
Conversely, the upper confidence graphs achieve high
true positive rates at the expense of lower precision
and higher false positive rates.
The boxplots also reveal considerable variability
across replications, particularly for small sample sizes
and larger dimensions.
As the sample size increases, the distributions generally
concentrate around their ideal values, although the
recovery problem remains more challenging for larger
and denser graphs.

\clearpage

\bibliographystyle{apalike}
\bibliography{biblio}

\end{document}